\documentclass{article}%
\usepackage{amsmath}
\usepackage{amsfonts}
\usepackage{amssymb}
\usepackage{graphicx}
\usepackage[top=2.5cm, bottom=3cm, left=2.5cm, right=2.5cm, headheight=0.5cm,
footskip=60pt]{geometry}
\usepackage{color}%
\providecommand{\U}[1]{\protect\rule{.1in}{.1in}}
\makeatletter
\newcommand{\pushright}[1]{\ifmeasuring@#1\else\omit\hfill$\displaystyle#1$\fi\ignorespaces}
\newcommand{\pushleft}[1]{\ifmeasuring@#1\else\omit$\displaystyle#1$\hfill\fi\ignorespaces}
\makeatother
\newtheorem{theorem}{Theorem}

\newtheorem{definition}{Definition}

\newtheorem{lemma}{Lemma}

\newtheorem{proposition}{Proposition}
\newtheorem{remark}{Remark}

\numberwithin{equation}{section}
\numberwithin{remark}{section}
\numberwithin{proposition}{section}
\numberwithin{definition}{section}
\numberwithin{lemma}{section}
\allowdisplaybreaks
\begin{document}

\title{\textcolor{red}{Extension of the Hoff solutions framework to cover Navier-Stokes
equations for a compressible fluid with anisotropic viscous-stress tensor}}
\author{D. Bresch\thanks{Univ. Grenoble Alpes, Univ. Savoie Mont-Blanc, CNRS, LAMA,
Chamb\'ery, France; didier.bresch@univ-smb.fr}, \hskip1cm C. Burtea
\thanks{Universit\'e de Paris and Sorbonne Universit\'e , CNRS, IMJ-PRG,
F-75006 Paris, France; cosmin.burtea@imj-prg.fr}\,}
\maketitle

\begin{abstract}

This paper deals with the Navier-Stokes system governing the evolution of a compressible barotropic fluid. We extend Hoff's intermediate regularity solutions framework \cite{hoff1995global,hoff1995strong} by relaxing the integrability needed for the initial density which is usually assumed to be $L^{\infty}$. By achieving this, we are able to take into account general fourth order symmetric viscous-stress tensors with coefficients depending smoothly on the time-space variables. More precisely, in space dimensions $d=2,3$, under periodic boundary conditions, considering a pressure law $p(\rho)=a\rho^{\gamma}$ whith $a>0$ respectively $\gamma\geq d/(4-d)$) and under the assumption that the norms of the initial data $\left(  \rho_{0}-M,u_{0}\right)  \in L^{2\gamma}\left(\mathbb{T}^{d}\right)  \times(H^{1}(\mathbb{T}^{d}))^{d}$ are sufficiently small, we are able to construct global weak solutions.  Above, $M$ denotes the total mass of the fluid while $\mathbb{T}$ with $d=2,3$ stands for periodic box. When comparing to the results known for the global weak solutions \`{a} la Leray, i.e. constructed assuming only the basic energy bounds, we obtain a relaxed condition on the range of admissible adiabatic coefficients $\gamma$.

\end{abstract}



\section{Introduction and main result}

In this paper, we study the problem of existence of global solutions with
intermediate regularity as pioneered by Hoff
\cite{hoff1995global,hoff1995strong} for the Navier-Stokes equations governing
the flow of a compressible fluid. Our aim is to extend the existence theory as
to accommodate general smooth fourth order symmetric viscous-stress tensors.
More precisely we consider the following system:
\begin{equation}
\left\{
\begin{array}
[c]{l}%
\rho_{t}+\operatorname{div}\left(  \rho u\right)  =0,\\
\left(  \rho u\right)  _{t}+\operatorname{div}\left(  \rho u\otimes u\right)
+\nabla(a\rho^{\gamma})=\operatorname{div}(\mathcal{\tilde{E}(}\nabla u)),
\end{array}
\right.  \label{system}%
\end{equation}
where $\left(  \rho,u\right)  $ represent the density and the velocity field
of the fluid. We assume that the pressure is given by $p(\rho)=a\rho^{\gamma}$
with $a>0$ and $\gamma\geq d/(4-d)$. We will restrict ourselves to the case of
$d=2,3$ space-dimensions and periodic boundary conditions. The $d$-dimensional
torus, $d\in{2,3}$ is denoted by $\mathbb{T}^{d}$. We consider a general
fourth order symmetric viscous stress tensor:
\[
\mathcal{\tilde{E}}=\left(  \tilde{\varepsilon}_{ijk\ell}\right)
_{i,j,k,\ell\in\overline{1,d}}%
\]
where we use the notations\footnote{$\overline{1,d}$ stands for the set
${1,2,\cdots,d}$} \footnote{All along this paper, we will use the Einstein
summation over repeated indices convention.}
\[
\mathcal{\tilde{E}\,(}\nabla u)=\tilde{\varepsilon}_{ijk\ell}\partial_{\ell
}u^{k},\text{ \ }\operatorname{div}\left(  \mathcal{\tilde{E}}\,\nabla
u\right)  =\partial_{j}\left(  \tilde{\varepsilon}_{ijk\ell}\partial_{\ell
}u^{k}\right)  .
\]
The system is completed with the initial data
\begin{equation}
\rho|_{t=0}=\rho_{0}\geq0,\qquad\rho u|_{t=0}=m_{0}. \label{Ini}%
\end{equation}
Most of the literature concerning compressible fluid mechanics deals with the
classical isotropic tensor
\begin{equation}
\mathcal{I=}\left(  \varepsilon_{ijk\ell}^{iso}\right)  _{i,j,k,\ell
\in\overline{1,d}} \label{iso1}%
\end{equation}
which is given by%
\begin{equation}
\varepsilon_{ijk\ell}^{iso}=\varepsilon_{k\ell ij}^{iso}=\left\{
\begin{array}
[c]{l}%
\mu\text{ if }\left(  i,j\right)  =\left(  k,\ell\right)  \text{ and }%
i\not =j,\\
2\mu+\lambda\text{ if }\left(  i,j\right)  =\left(  k,\ell\right)  \text{ and
}i=j,\\
0\text{ otherwise,}%
\end{array}
\right.  \label{iso2}%
\end{equation}
and $\mu,\lambda>0$ are given constants. This implies in particular that one
has
\[
-\operatorname{div}(\mathcal{I\nabla}u)=-\mu\Delta u-(\mu+\lambda
)\nabla\operatorname{div}u.\,
\]
In the following, we will recall some well known results concerning the
existence of solutions for compressible Navier-Stokes equations for both
isotropic viscous tensor and for anisotropic viscous tensors. In order to
understand why we are require to work within the Hoff-solutions framework and,
in particular, why it is necessary to relax the integrability of density to
$L^{p}$ with $p<\infty$, we will briefly pass in review different notions of
solutions (strong solutions, critical spaces, global weak solutions \`{a} la
Leray, intermediate regularity \`{a}~la Hoff). In this paper, we do not
discuss density dependent viscosity for compressible Navier-Stokes equations.
\bigskip

\noindent1) \textit{A short review of know results for isotropic stress
tensors.} The study of system $\left(  \text{\ref{system}}\right)
-\eqref{Ini}$ with a given pressure law $s\mapsto p(s)$ in the case of
isotropic stress tensors goes back to the work of J. Nash
\cite{nash1962probleme} where the author shows the existence of local-in-time
strong solutions in H\"{o}lder spaces. Then local strong existence for initial
data in Sobolev spaces was investigated by Solonnikov
\cite{solonnikov1980solvability} in the $80^{\prime}s$ while the first global
result is due to Matsumura and Nishida \cite{matsumura1980initial} where they
prove the existence of global-in-time solutions in $3D$ if the initial data
are sufficiently close to equilibrium in $H^{3}$. In the $2000s$, R. Danchin
\cite{danchin2000global} constructed global small solutions in the so-called
critical spaces.
Smooth solutions for the IBVP with Dirichlet boundary conditions in dimension $d\ge2$ are known to blow-up since the work of Vaigant \cite{vaigant1994example} where the author constructs an explicit solution for the NSC system and shows that the $L^{\infty}$-norm of the density blows-up in finite time.
Very recently, \cite{merle2019implosion} F. Merle, P. Rapha\"{e}l, I.
Rodnianski and J.Szeftel prove that in $3D$, for small $\gamma\leq1+2/\sqrt
{3}$ there exists local smooth solutions which explode in finite time : the
$L^{\infty}$-norms of the density and the velocity blow-up. Thus, in some
sense, the smallness condition, which express the fact that the initial
configuration is sufficiently close to an constant equilibrium state, is
necessary in order to insure global well-posedness. Let us mention the recent result of R. Danchin and P.B. Mucha \cite{DanchinMucha2019} where the authors construct global solutions in the two-dimensional case requiring only that the divergence of the velocity field should be small.

Another category of results regarding the solvability of $\left(
\text{\ref{system}}\right)  -\eqref{Ini}$ concerns the so-called
weak-solutions \`{a}~la Leray: solutions in the sense of distributions
satisfying the energy inequality for which one can guarantee their global
existence for arbitrary large initial data. Of course, few things are known
regarding the uniqueness of these solutions. We mention the, by now classical
results of P.L. Lions \cite{lions1996mathematical2}, E. Feireisl et al.
\cite{feireisl2001existence}. Recently, the first author and P.--E. Jabin
extended these two results in order to cover on one hand some anisotropic
stress tensors \cite{bresch2018global} and on the other hand more general
pressure functions \cite{bresch2018global,bresch2021compressible} that could
not be treated within the Lions-Feireisl theory. We will return back and
comment a bit more on these results in the context of anisotropy.

A third category of results concerns an intermediate regularity functional
framework which was pioneered in the works of D. Hoff
\cite{hoff1995global,hoff1995strong,hoff2002dynamics,hoff2008lagrangean} (that
we will called solutions \`{a} la Hoff) and B. Desjardins
\cite{desjardins1997regularity}. By intermediate regularity we mean of course
between the regularity needed to construct strong solutions and weak solutions
\`{a} la Leray (see \cite{lions1996mathematical2}). \emph{These solutions are
interesting since they allow to work with discontinuous densities while
granting some extra regularity for the velocity field which turns out to be
sufficiently regular in order to generate a log-Lipschitz flow}. These
solutions were used by D. Hoff to study the dynamics of a surface of
discontinuities initially present in the density, see
\cite{hoff2002dynamics,hoff2008lagrangean} and found applications in the
context of multifluid flows, see the work of the first author
and X. Huang \cite{bresch2011multi}. Since the present work deals with these
kind of solutions, we will take the time to give more details. In
\cite{hoff1995global,hoff1995strong}, for the case of isotropic stress
tensors, D. Hoff introduced and studied the properties of two energy-type
functionals
\begin{equation}
A_{1}\left(  t\right)  =\frac{\mu\sigma\left(  t\right)  }{2}\sum_{i=1}%
^{d}\sum_{k=1}^{d}\int\left\vert \partial_{k}u^{i}\left(  t\right)
\right\vert ^{2}+\frac{\left(  \mu+\lambda\right)  \sigma\left(  t\right)
}{2}\int\left\vert \operatorname{div}u\left(  t\right)  \right\vert ^{2}%
+\int_{0}^{t}\int\sigma\rho\left\vert \dot{u}\right\vert ^{2}
\label{Hoff1_Intro}%
\end{equation}
and
\begin{equation}
A_{2}\left(  t\right)  =\sigma^{1+d}\left(  t\right)  \int\frac{\rho\left(
t\right)  \left\vert \dot{u}\left(  t\right)  \right\vert ^{2}}{2}+\mu
\sum_{i=1}^{d}\sum_{k=1}^{d}\int_{0}^{t}\int\sigma^{1+d}\left\vert
\partial_{k}\dot{u}^{i}\right\vert ^{2}+\left(  \mu+\lambda\right)  \int
_{0}^{t}\int\sigma^{1+d}\left\vert \operatorname{div}\dot{u}\right\vert ^{2}
\label{Hoff2_Intro}%
\end{equation}
where
\[
\dot{u}=u_{t}+u\cdot\nabla u\text{ and }\sigma\left(  t\right)  =\min\left\{
1,t\right\}  .
\]
These functionals naturally appear: The first one when multiplying the
momentum equation with $\sigma\dot{u}$ and integrating while the other one
appears when applying $\partial_{t}+\operatorname{div}\left(  u\cdot\right)  $
to the momentum equation and multiplying with $\sigma^{1+d}\dot{u}$. D. Hoff
shows that $A_{1},A_{2}$ can be controlled globally in time if the initial
data have suitably small energy and $\rho_{0}$\emph{ is close to a constant in
}$L^{\infty}$. The fact that these two functionals can be controlled translate
some fine smoothing properties due to the diffusion: it turns out that $u$ is
H\"{o}lder continuous in time-space, far from $t=0$ and that
$\operatorname{curl}u$ and the effective flux%
\[
F=\left(  2\mu+\lambda\right)  \operatorname{div}u-p\left(  \rho\right)
\]
are $H^{1}$ in space for a.e. $t>0$. In particular, this later properties
render mathematically clear the fact that discontinuities in the density are
advected by the flow but in such a way that the so called effective flux, i.e.
$F$ stays fairly smooth. This property enjoyed by the effective flux, known
and exploited in the $1d$ case in \cite{HS85,hoff1987global}, also turns out
to be crucial when showing the stability of sequences of weak-solutions. In
order to give a meaning to $A_{1}$, $A_{2}$ very little extra information is
need when comparing to the energy level, by which we mean $\rho_{0}u_{0}%
^{2}\in L^{1}$, $\rho_{0}\in L^{\gamma}$, which essentially is that $\rho
_{0}\in L^{\infty}$ and $u_{0}\in L^{2^{n}}$ (in the whole space case). If
more information is available for the initial data, modified versions of the
two functionals can be used: for instance if $u_{0}\in H^{1}$, one can control%
\begin{equation}
\tilde{A}_{1}\left(  t\right)  =\frac{\mu}{2}\sum_{i=1}^{d}\sum_{k=1}^{d}%
\int\left\vert \partial_{k}u^{i}\left(  t\right)  \right\vert ^{2}+\frac
{\mu+\lambda}{2}\int\left\vert \operatorname{div}u\left(  t\right)
\right\vert ^{2}+\int_{0}^{t}\int\rho\left\vert \dot{u}\right\vert ^{2},
\label{Hoff11_Intro}%
\end{equation}
respectively
\begin{equation}
\tilde{A}_{2}\left(  t\right)  =\sigma\left(  t\right)  \int\frac{\rho\left(
t\right)  \left\vert \dot{u}\left(  t\right)  \right\vert ^{2}}{2} +\mu
\sum_{i=1}^{d}\sum_{k=1}^{d}\int_{0}^{t}\int\sigma\left\vert \partial_{k}%
\dot{u}^{i}\right\vert ^{2} +\left(  \mu+\lambda\right)  \int_{0}^{t}%
\int\sigma\left\vert \operatorname{div}\dot{u}\right\vert ^{2},
\label{Hoff22_Intro}%
\end{equation}
which of course express the fact that due to the extra information the
solution is better behaved close to the initial time layer $t=0$. We also
mention the related but independent work of B. Desjardins
\cite{desjardins1997regularity} where the author obtains local in time results
showing that is possible to control a function which is essentially equivalent
to $\tilde{A}_{1}$. \emph{We mention that in all the above cited papers, the
assumption that }$\rho_{0}\in L^{\infty}$\emph{ turns out to be crucial}.
\emph{The fact that one can propagate control of the }$L^{\infty}$\emph{-norm
of the density heavily depends on the algebraic structure of the isotropic
Navier-Stokes system throughout the so called-effective flux }$F=\left(
2\mu+\lambda\right)  \operatorname{div}u-p\left(  \rho\right)  $\emph{ defined
above. }

\bigskip

\noindent2) \textit{The case of anisotropic stress tensors.} In this case, the
mathematical results are in short supply. Let us mention that in the context
of strong solutions \cite{matsumura1980initial}, \cite{danchin2000global}
where the results are proved by maximal regularity results, at least if the
stress tensor is "close enough to the isotropic tensor" then there should
virtually be little change needed in order to accommodate these kind of
solutions. However, as explained above, when dealing with classical solutions
the density is a continuous function thus excluding many interesting
situations in applications (for example, mixtures of fluids).

The first paper providing a result in this direction has been obtained by the
first author and P.--E. Jabin in \cite{bresch2018global} and concerns the
existence of global weak solutions \`{a}~la Leray with an anisotropic
diffusion of the form:
\begin{equation}
-\mathrm{div}(\tilde{A}(t)\nabla u)-(\mu+\lambda)\nabla\mathrm{div}u
\label{BJ1}%
\end{equation}
where
\[
\tilde{A}(t)=\mu\mathrm{Id}+ A(t)\text{ with }\mu>0.
\]
The result proved in \cite{bresch2018global} states that there exists an
universal constant $c>0$ such that if\footnote{$d$ stands for the space
dimension.}
\[
\Vert A(t)\Vert_{L^{\infty}}\leq c\left(  \frac{2\mu}{d}+\lambda\right)
\]
and if
\begin{equation}
\gamma>\frac{d}{2}\left[  \left(  1+\frac{1}{d}\right)  +\sqrt{1+\frac
{1}{d^{2}}}\right]  \label{restrict_BreschJabin}%
\end{equation}
then, there exists global weak solutions \`{a} la Leray for the Navier-Stokes
system $\left(  \text{\ref{system}}\right)  -\eqref{Ini}$ with the
viscous-stress tensor given by $\left(  \text{\ref{BJ1}}\right)  $. This
result extended to the anisotropic case the global existence of weak solution
\`{a} la Leray obtained for the isotropic case in
\cite{lions1996mathematical2}, \cite{feireisl2001existence}. The result of the
first author and P.-E. Jabin is based on new estimates for the transport
equation. This result requires in a crucial manner some form of compactness in
space for
\[
\left(  2\mu+\lambda\right)  \operatorname{div}u-L(a\rho^{\gamma})
\]
where $L$ is a non-local operator of order $0$. \textit{It is at this level
that the authors use the fact that }$A$\textit{ depends only on time has been
used by the authors.} \textit{The extension of this result to space dependent
strain tensors represents a serious difficulty that remains an open problem}.
Moreover, the restriction for the adiabatic coefficient $\gamma$ given by
$\left(  \text{\ref{restrict_BreschJabin}}\right)  $ excludes most of the
physically realistic values : monoatomic gases $5/3$, ideal diatomic gases
$7/5$, viscous shallow--water $\gamma=2$.

Let us also mention our results concerning related models for compressible
fluids, on the one hand, existence of global weak solutions \`{a} la Leray for
the quasi-stationary compressible Stokes in \cite{BrBu1} where an anisotropic
diffusion $-\mathrm{div}(AD(u))$ is considered and where no smallness
assumption on the anisotropic amplitude is needed in order to develop the
existence theory and on the other hand our result regarding the stationary
compressible Navier-Stokes equations in \cite{BrBu2} where we treat a viscous
diffusion operator given by $-{\mathcal{A}}u$ (under some constraints) where
${\mathcal{A}}$ is composed by a classical constant viscous part plus an
anisotropic contribution and a possible nonlocal contribution.

\bigskip

\noindent3) \textit{Motivation to extend the framework of Hoff-type solutions
and description of our main result.} When dealing with weak solutions for
non-linear PDE systems, one of the most delicate aspects is the stability
analysis: given a sequence of weak solutions for some well-chosen approximated
systems, show that this sequence converges to a solution for the initial
system. The key ingredient in \cite{BrBu1} and \cite{BrBu2} is an identity
that we found when comparing on the one hand, the limiting energy equation and
on the other hand, the equation of the energy associated to the limit system.
In order to justify such an identity, a crucial assumption seems to be the
fact that the pressure is $L^{2}$, an apriori estimate which is ensured by
basic a-priori estimates in the case of the Stokes system or for the
stationary Navier-Stokes system. However, in the case of system $\left(
\text{\ref{system}}\right)  -\eqref{Ini}$ in the isotropic case, the best
estimate for the density is due to P. Lions who showed for global weak
solutions \`{a} la Leray that $\rho\in$ $L_{t,x}^{5\gamma/3-1}$. This makes it
impossible to write the energy equation because, loosely speaking, the
velocity cannot be used as a test function in a weak-formulation of $\left(
\text{\ref{system}}\right)  $. Thus, it seems hopeless to justify the limiting
passage as in \cite{BrBu1} and \cite{BrBu2} in the most general setting of
weak-solutions \`{a} la Leray. Obviously, one may ask if we can work in an
intermediate regularity setting. However, one learns fast that we are faced
with a serious problem when trying to propagate the $L^{\infty}$-norm for the
density. In the isotropic case, this is achieved based on the fact that the
the second Hoff functional, namely $\left(  \text{\ref{Hoff2_Intro}}\right)  $
controls, at least far from the initial time layer $t=0$, the $L^{\infty}%
$-norm of the effective viscous flux:
\begin{equation}
\left(  2\mu+\lambda\right)  \operatorname{div}u-\left(  a\rho^{\gamma}%
-\int_{\mathbb{T}^{d}}a\rho^{\gamma}\right)  =\Delta^{-1}\operatorname{div}%
(\rho\dot{u}). \label{effective_flux_iso}%
\end{equation}
Of course, the situation is not the same in the anisotropic case, where
\begin{equation}
\left(  2\mu+\lambda\right)  \operatorname{div}u-\left(  a\rho^{\gamma}%
-\int_{\mathbb{T}^{d}}a\rho^{\gamma}\right)  =\Delta^{-1}\operatorname{div}%
(\rho\dot{u})+\Delta^{-1}\operatorname{div}\operatorname{div}\left(  \left(
\mathcal{\tilde{E}-I}\right)  (\nabla u)\right)  \label{effective_flux_aniso}%
\end{equation}
and the term $\Delta^{-1}\operatorname{div}\operatorname{div}\left(  \left(
\mathcal{\tilde{E}-I}\right)  (\nabla u)\right)  $, being of the same order as
$\nabla u$ we cannot expect it to be $L^{\infty}$. Because of the lack of
algebraic structure we are led to abandon any hope of propagating an
$L^{\infty}$ bound for the density. \emph{A natural question then appears: is
it possible to bound the Hoff functionals without working in an }$L^{\infty}%
$\emph{ framework for the density? The main contribution of this paper is to
show that this is indeed the case. Of course, this fact makes it possible to
construct global weak solutions close to equilibrium for the Navier--Stokes
system in the anisotropic case in an intermediate regularity setting. }This
program requires establishing $L^{p}$-estimates for the density that are
compatible with the Hoff functionals. In order to avoid further technical
difficulties, we will assume the best information possible for the velocity,
namely $u_{0}\in\left(  H^{1}\right)  ^{d}$ such that we will rather work with
the anisotropic equivalent of the functionals defined in $\left(
\text{\ref{Hoff11_Intro}}\right)  -\left(  \text{\ref{Hoff22_Intro}}\right)  $.

\smallskip

Our result should be seen to be complementary to the work of the first author
and P.-E. Jabin \cite{bresch2018global}. In particular, this extended
intermediate regularity framework allows us to:

\begin{itemize}
\item consider viscous-stress tensors depending on time and also on the space variable;

\item consider a range of adiabatic coefficient namely%
\[
\gamma\geq\frac{d}{4-d} \hbox{ for } d\in\left\{  2,3\right\}  .
\]
\emph{ In particular, in }$2D$\emph{ we are able to treat all coefficients
that are of practical interest }$\gamma\geq1$\emph{.}

\item this method could be adapted to bounded domains with Dirichlet boundary
conditions for which the existence result of global weak solutions (\`{a} la
Leray or intermediate regularity) with anisotropic tensors remains open.
\end{itemize}

Note that the range for the coefficient $\gamma$ is larger than the one in
\cite{bresch2018global} namely $\left(  \text{\ref{restrict_BreschJabin}%
}\right)  $ and we cover strain tensors which may depend on the space
variable. Of course, the price to pay is that the initial conditions are
supposed to be close to equilibrium and that we require the initial velocity
field to be in $(H^{1}(\mathbb{T}^{d}))^{d}$.

\medskip

\noindent\textit{Assumptions and notations.} We will rather write
$\mathcal{\tilde{E}}=\mathcal{I+}\mathcal{E}$ with $\mathcal{I}$, the usual
isotropic tensor $\left(  \text{\ref{iso1}}\right)  $-$\left(
\text{\ref{iso2}}\right)  $ and where $\mathcal{E}$ measures in some sense the
anisotropic perturbation. With these new notations
and setting the constant $a=1$ in the pressure term\footnote{Of
course, this does not affect the generality of the result, the choice $a=1$ is
only a matter of simplifying the computations.}, the system $\left(
\text{\ref{system}}\right)  $ becomes
\begin{equation}
\left\{
\begin{array}
[c]{l}%
\rho_{t}+\operatorname{div}\left(  \rho u\right)  =0,\\
\left(  \rho u\right)  _{t}+\operatorname{div}\left(  \rho u\otimes u\right)
+\nabla\rho^{\gamma}=\mu\Delta u+\left(  \mu+\lambda\right)  \nabla
\operatorname{div}u+\operatorname{div}(\mathcal{E(}\nabla u)).
\end{array}
\right.  \label{system_pert}%
\end{equation}
We suppose that $\mu,\lambda\in\mathbb{R}$ such that%
\begin{equation}
\mu>0\text{ and }\mu+\lambda\geq0. \label{conditie_mu_lambda}%
\end{equation}
We will assume that $\mathcal{E}$ $\mathcal{=}\left(  \varepsilon_{ijk\ell
}\right)  _{i,j,k,\ell\in\overline{1,d}}$ verifies the following properties:

\begin{itemize}
\item For all $i,j,k,\ell\in\overline{1,d}$ we assume the following symmetry
property:
\begin{equation}
\varepsilon_{ijk\ell}=\varepsilon_{k\ell ij}. \tag{$\textrm{H1}$}\label{Hyp1}%
\end{equation}
The latter property ensures that\footnote{All along this
paper, we will use the Einstein summation over repeated indices convention.}%
\[
\varepsilon_{ijk\ell}a_{ij}b_{k\ell}=\varepsilon_{ijk\ell}b_{ij}a_{k\ell}.
\]

\item Strict coercivity of the diffusive part:%
\begin{equation}
\overline{\varepsilon}\left\vert a_{ij}\right\vert ^{2}\geq\varepsilon
_{ijk\ell}a_{ij}a_{k\ell}\geq-\underline{\varepsilon}\left\vert a_{ij}%
\right\vert ^{2} \tag{$\mathrm{H2}$}\label{Hyp2}%
\end{equation}
where $\underline{\varepsilon},\overline{\varepsilon}>0$ such that%
\[
0<\mu-\underline{\varepsilon}.
\]

\item Regularity: for all $i,j,k,\ell\in\overline{1,d}$, $\varepsilon
_{ijk\ell}\in W^{1,\infty}\left(  \mathbb{(}0,\infty\mathbb{)\times T}%
^{d}\right)  $ with
\begin{equation}
\text{ }\left\Vert \partial_{t}\varepsilon_{ijk\ell}\right\Vert _{L^{\infty
}((0,+\infty)\times\mathbb{T}^{d})}+\left\Vert \nabla\varepsilon_{ijk\ell
}\right\Vert _{L^{\infty}((0,+\infty)\times\mathbb{T}^{d})}<\infty
.\tag{$\mathrm{H3}$}\label{Hyp3}%
\end{equation}

\end{itemize}

\noindent\textit{Main Result.} Let us define the following:
\[
E\left(  \rho/M,u\right)  =\int_{\mathbb{T}^{d}}\bigl(H_{1}(\rho/M)+\frac
{1}{2}\rho|u|^{2}\bigr)
\]
with $0<M<+\infty$ and where
\begin{equation}
H_{1}(\rho/M)=H_{1}(\rho)-H_{1}(M)+H_{1}^{\prime}(M)(\rho-M) \label{has1}%
\end{equation}
with
\[
H_{1}(\rho)=\rho\int_{0}^{\rho}P(s)/s^{2}\,ds=\frac{\rho^{\gamma}}{\gamma-1}.
\]
Also, we introduce
\begin{equation}
H_{\ell}(\rho/M)=\rho\int_{M}^{\rho}\frac{|P(s)-P(M)|^{\ell-1}(P(s)-P(M))}%
{s^{2}}\,ds\text{ with }\ell\in\left\{  2,3\right\}  . \label{Hshurile}%
\end{equation}
We are now in the position of stating our main result:

\begin{theorem}
\label{MainTheorem} Consider $\mu,\lambda\in\mathbb{R}$ such that $\mu
>0,\mu+\lambda>0$. Let $\mathcal{\tilde{E}}=\mathcal{I+}\mathcal{E}$ with
$\mathcal{I}$ the isotropic viscous-stress tensor defined in
\eqref{iso1}--\eqref{iso2} and $\mathcal{E=}\left(  \varepsilon_{ijk\ell
}\right)  _{i,j,k,\ell\in\overline{1,d}}$ a fourth order tensor verifying the
hypothesis \eqref{Hyp1}--\eqref{Hyp3}.  Then, there exists two positive
constant $\eta$,$c_{0}$ independent of $\mu$ and $\lambda$ such that the
following holds true. Assume that
\begin{equation}
\left\Vert \mathcal{E}\right\Vert _{L^{\infty}((0,+\infty)\times\mathbb{T}%
^{d})}=\sup_{i,j,k,l\in\overline{1,d}}\left\Vert \varepsilon_{ijk\ell
}\right\Vert _{L^{\infty}((0,+\infty)\times\mathbb{T}^{d})}\leq\eta
\min\left\{  \mu,2\mu+\lambda\right\}  .\tag{$\mathrm{H4}$}\label{Hyp4}%
\end{equation}
Then, for any $\left(  \rho_{0},u_{0}\right)  \in L^{2\gamma}\left(
\mathbb{T}^{d}\right)  \times\left(  H^{1}\left(  \mathbb{T}^{d}\right)
\right)  ^{d}$ with
\[
\int_{\mathbb{T}^{d}}\rho_{0}=M,\text{ }\int_{\mathbb{T}^{d}}\rho_{0}%
u_{0}=\mathcal{P\in}\mathbb{R}^{d},
\]
such that%
\[
E\left(  \rho_{0}/M,u_{0}\right)  +\int_{\mathbb{T}^{d}}H_{2}\left(  \rho
_{0}/M\right)  +\left\Vert u_{0}\right\Vert _{(H^{1}(\mathbb{T}^{d}))^{d}}%
^{2}\leq c_{0},
\]
there exists a global weak solution $\left(  \rho,u\right)  $ for
\eqref{system}--\eqref{Ini} with
\[
(\rho-M,\rho u-\mathcal{P})\in\mathcal{C}([0,+\infty);H^{-1}({\mathbb{T}}%
^{d}))\times\mathcal{C}([0,+\infty);(H^{-1}({\mathbb{T}}^{d}))^{d})
\]
and such that for all $t\geq0$ we have:
\[%
\begin{array}
[c]{lll}%
E(\rho\left(  t\right)  /M,u\left(  t\right)  )+(\mu-\underline{\varepsilon})%
{\displaystyle\sum\limits_{i=1}^{d}}
{\displaystyle\sum\limits_{k=1}^{d}}
{\displaystyle\int_{\mathbb{T}^{d}}}
\left\vert \partial_{k}u^{i}\left(  t\right)  \right\vert ^{2}+\left(
\mu+\lambda\right)
{\displaystyle\int_{\mathbb{T}^{d}}}
\left\vert \operatorname{div}u\right\vert ^{2} & \leq & E(\rho_{0}/M,u_{0}),\\
\text{ }\frac{1}{2}\left\{  \mu%
{\displaystyle\sum\limits_{i=1}^{d}}
{\displaystyle\sum\limits_{k=1}^{d}}
{\displaystyle\int_{\mathbb{T}^{d}}}
\left\vert \partial_{k}u^{i}\right\vert ^{2}\left(  t\right)  +\left(
\mu+\lambda\right)
{\displaystyle\int_{\mathbb{T}^{d}}}
\left\vert \operatorname{div}u\right\vert ^{2}\left(  t\right)  \right\}  +%
{\displaystyle\int_{0}^{t}}
{\displaystyle\int_{\mathbb{T}^{d}}}
\rho\left\vert \dot{u}\right\vert ^{2} & \leq & Cc_{0},\\
\text{ }\sigma\left(  t\right)
{\displaystyle\int_{\mathbb{T}^{d}}}
\dfrac{\rho\left(  t\right)  \left\vert \dot{u}\left(  t\right)  \right\vert
^{2}}{2}+\mu%
{\displaystyle\sum\limits_{i=1}^{d}}
{\displaystyle\sum\limits_{k=1}^{d}}
{\displaystyle\int_{0}^{t}}
{\displaystyle\int_{\mathbb{T}^{d}}}
\sigma\left\vert \partial_{k}\dot{u}^{i}\right\vert ^{2}+\left(  \mu
+\lambda\right)
{\displaystyle\int_{0}^{t}}
{\displaystyle\int_{\mathbb{T}^{d}}}
\sigma\left\vert \operatorname{div}\dot{u}\right\vert ^{2} & \leq & Cc_{0},\\%
{\displaystyle\int_{\mathbb{T}^{d}}}
H_{2}\left(  \rho(t)/M\right)  +\sigma\left(  t\right)
{\displaystyle\int_{\mathbb{T}^{d}}}
H_{3}\left(  \rho(t)/M\right)  +%
{\displaystyle\int_{0}^{t}}
{\displaystyle\int_{\mathbb{T}^{d}}}
\left\vert P\left(  \rho\right)  -P\left(  M\right)  \right\vert ^{3}+%
{\displaystyle\int_{0}^{t}}
{\displaystyle\int_{\mathbb{T}^{d}}}
\sigma\left\vert P\left(  \rho\right)  -P\left(  M\right)  \right\vert ^{4} &
\leq & Cc_{0}.
\end{array}
\]
where $\sigma(t)=\min\{1,t\}$ while $C=C\left(  \mu,\lambda,\gamma
,M,E_{0},c_{0}\right)  $ is a constant that depends on $\mu,\lambda
,\gamma,M,E_{0}$.
\end{theorem}


\smallskip

\begin{remark}
\label{HoffRem} It is important to remark that it seems a difficult problem to
propagate the $L^{\infty}$-norm for the density as it has been done by D. Hoff
for the isotropic compressible Navier-Stokes equations with a barotropic
pressure law.$\,$
\end{remark}

\smallskip

\begin{remark}

The uniform control of the Hoff functionals that we obtain in Theorem \ref{MainTheorem} plays a crucial role in the stability part of the proof. For instance, the fact that the pressure is bounded in $L_{t,x}^{3}$ allows us to justify the equation
\[
\partial_{t}P\left(  \rho\right)  +\operatorname{div}\left(  P\left(
\rho\right)  u\right)  +\left(  \rho P^{\prime}\left(  \rho\right)  -P\left(
\rho\right)  \right)  \operatorname{div}u=0,
\]
from the mass equation
\[
\partial_{t}\rho+\operatorname{div}\left(  \rho u\right)  =0,
\]
of the limit system. Another crucial aspect is that when considering a
sequence of solutions of systems that approximate the Navier-Stokes system,
controlling the second Hoff functional allows to obtain information for the
time derivative of the velocities. As a consequence of the Aubin-Lions lemma,
we obtain that the sequence of velocities converges strongly in $L_{t,x}^{2},$
at least far from $t=0$ which is crucial in order to implement the idea from
\cite{BrBu1}.
\end{remark}

\noindent\textit{Main steps and organization of the paper.} We detail below
the main steps of the proof of Theorem \ref{MainTheorem}. Inspired
by the approximate system proposed by the first author and
P.--E. Jabin in \cite{bresch2018global}, we will consider a regularized
version of the Navier-Stokes system $\left(  \text{\ref{system}}\right)  $:%
\begin{equation}
\left\{
\begin{array}
[c]{l}%
\rho_{t}+\operatorname{div}\left(  \rho u\right)  =0,\\
\left(  \rho u\right)  _{t}+\operatorname{div}\left(  \rho u\otimes u\right)
+\nabla\rho^{\gamma}=\mu\Delta u+\left(  \mu+\lambda\right)  \nabla
\operatorname{div}u+\omega_{\delta}\ast\operatorname{div}(\mathcal{E(}%
\nabla\omega_{\delta}\ast u)),
\end{array}
\right.  \label{ANS_delta}%
\end{equation}
where
\begin{equation}
\omega_{\delta}\left(  x\right)  =\frac{1}{\delta^{d}}\omega(\frac{x}{\delta
}),
\end{equation}
with $\omega$ a smooth, nonnegative, radial function compactly supported in
the unit ball centered at the origin and with integral equal to $1$. Since
system $\left(  \text{\ref{ANS_delta}}\right)  $ can be seen as a regular
perturbation of the Navier-Stokes system for a compressible barotropic fluid,
classical results
\cite{solonnikov1980solvability,desjardins1997regularity,danchin2010solvability}
can be invoked in order to ensure the existence of a \emph{local} classical solution.

\smallskip

\begin{remark}
\label{important} \noindent\textbf{(Important remark on the anisotropy).} To
simplify the writing of the paper, we will assume in the proof that
\[
\int_{0}^{T}\int_{{\mathbb{T}}^{d}}{\mathcal{E}}(\nabla w):\nabla w\geq0
\qquad\text{ for all }T\in(0,+\infty]\text{ and }w\in L^{2}((0,T);H^{1}\left(
\mathbb{T}^{d}\right)  ).
\]
This assumption is needed in order to treat the stability of weak-solutions of
system $\left(  \text{\ref{ANS_delta}}\right)  $ part of the proof. In order
to avoid this assumption and treat the general case, it is sufficient to
consider an approximate system with diffusion given by%
\[
(\mu-\underline{\varepsilon})\Delta u+\left(  \mu+\lambda\right)
\nabla\operatorname{div}u+\underline{\varepsilon}\Delta\omega_{\delta}\ast
u+\omega_{\delta}\ast\operatorname{div}(\mathcal{E(}\nabla\omega_{\delta}\ast
u)).
\]
change the coefficients $\lambda$ and $\mu$ in the isotropic part to allow to
satisfy these assumptions.
\end{remark}

We show that these solutions have the property that the two Hoff functionals
associated are bounded independently of $\delta$. This is one of the main
contributions of this paper.

\medskip

Following exactly the same steps as in R. Danchin and P.B. Mucha
\cite{DanchinMucha2019}, shows that the local solutions of $\left(
\text{\ref{ANS_delta}}\right)  $ can be prolonged to global ones. The fact
that the Hoff functionals are independent of $\delta$ is of course crucial in
order to show that we can extract a subsequence converging to a weak-solution
\`{a} la Hoff of $\left(  \text{\ref{system}}\right)  -\eqref{Ini}$. Here, we
are faced with, let us say the classical difficulty in compressible fluid
mechanics which is to be able to identify the pressure in the limit. More
precisely,
\begin{equation}
\lim\left(  \rho_{n}\right)  ^{\gamma}=\left(  \lim\rho_{n}\right)  ^{\gamma}.
\label{pressure_ident}%
\end{equation}
Of course, when dealing with weak solutions,
the density is just a $L^p$
function for some $p<\infty$ and no gain of regularity is to be expected.
Since weak limits, in general, do not commute with nonlinear functions,
showing $\left(  \text{\ref{pressure_ident}}\right)  $ has to take into
consideration some algebraic properties of solutions of the NS system. Let us
recall that classical techniques due to P.L. Lions
\cite{lions1996mathematical2} and E. Feireisl \cite{feireisl2001compactness}
do not apply in this context, see the discussions from the introductions of
\cite{bresch2018global},\cite{BrBu1},\cite{BrBu2} for more details. Moreover,
the work by the first author and P.E. Jabin requires a relative large $\gamma$
and, maybe more importantly, as it was explained above, it is not
straightforward to extend it to heterogenous in space anisotropic tensors (the
fact that $\mathcal{E}$ can depend also on the space variable). \emph{Here, it
is crucial to extend our idea from \textrm{\cite{BrBu1}}} that we successfully
implemented in order to construct global weak solutions \`{a}~la Leray for the
Stokes-Brinkman system in \cite{BrBu1} and for the stationary NS system in
\cite{BrBu2}. In these two papers, we did not need however to impose any
restriction on the size of the initial data or the forcing terms. This is
essentially due to the fact that, in the previous cases, the pressure turns
out to be an $L_{t,x}^{2}$ function (if $\gamma$ is large for the stationary
NS system) and that the convective term behaves better in the aforementioned cases.

\bigskip

The rest of the paper unfolds as follows. In the second section, we prove the
main result: First we recall basic mass conservation and energy estimate,
secondly we extend the Hoff estimates in a $L^{p}$ framework, third we
construct a sequence of approximate solutions and then finally we show the
stability property. In an appendix, we present a tool box with Fourier
multipliers properties, Sobolev inequality and Gronwall-Bihari inequality and
finally we give the detailed computations for the Hoff functionals that we
strongly use.

\section{Proof of the main result}

\subsection{Basic mass conservation and energy estimate}

\noindent\textbf{The conservation of mass and momentum.} The simplest
\textit{a priori} estimate we have is given by the conservation of mass and
momentum:%
\begin{align}
\int_{\mathbb{T}^{d}}\rho\left(  t\right)   &  =\int_{\mathbb{T}^{d}}\rho
_{0}:=M>0,\label{mass}\\
\int_{\mathbb{T}^{d}}\rho\left(  t\right)  u\left(  t\right)   &
=\int_{\mathbb{T}^{d}}m_{0}:=\mathcal{P}\in\mathbb{R}^{d}.
\end{align}

\noindent\textbf{The energy estimate.} From the continuity equation, we can
also deduce the following equation
\begin{equation}
\partial_{t}b\left(  \rho\right)  +\operatorname{div}\left(  b\left(
\rho\right)  u\right)  +\left(  \rho b^{\prime}(\rho)-b\left(  \rho\right)
\right)  \operatorname{div}u=0, \label{renormm}%
\end{equation}
\textit{a priori} for all $b$ verifing the conditions from Proposition
\ref{Prop_ren1}. Taking $b\left(  \rho\right)  =\rho^{\alpha}$ in $\left(
\text{\ref{renormm}}\right)  $ yields
\[
\partial_{t}\rho^{\alpha}+\operatorname{div}\left(  \rho^{\alpha}u\right)
+\left(  \alpha-1\right)  \rho^{\alpha}\operatorname{div}u=0.
\]
Also, we can write that%
\begin{align*}
u\cdot\nabla P\left(  \rho\right)   &  =\operatorname{div}\left(  u(P\left(
\rho\right)  -P(M)\right)  -(P\left(  \rho\right)  -P\left(  M\right)
)\operatorname{div}u\\
&  =\operatorname{div}\left(  u(P\left(  \rho\right)  -P(M)\right)  +\frac
{d}{dt}H_{1}\left(  \rho/M\right)  +\operatorname{div}\left(  H_{1}\left(
\rho/M\right)  u\right)  ,
\end{align*}
where $H_{1}\left(  \rho/M\right)  $ has been defined in $\left(
\text{\ref{has1}}\right)  $. The function $H_{1}\left(  \rho/M\right)  $ is
more appropriate in order to study densities that are close to some constant
state. Thus we get the following energy estimate 
\begin{align*}
&  \int_{\mathbb{T}^{d}}\left(  H_{1}\left(  \rho/M\right)  +\frac{\rho |u|^{2}}{2}\right)  +\mu\int_{0}^{t}\int_{\mathbb{T}^{d}}\left\vert \nabla
u\right\vert ^{2}+\left(  \mu+\lambda\right)  \int_{0}^{t}\int_{\mathbb{T}^{d}}\left\vert \operatorname{div}u\right\vert ^{2}+\int_{0}^{t}\int_{\mathbb{T}^{d}}\varepsilon_{ijkl}\partial_{\ell}u_{\delta}^{k}\partial_{j}u_{\delta}^{i}\\
&  \leq\int_{\mathbb{T}^{d}}\left(  H_{1}\left(  \rho_{0}/M\right)  +\frac{\rho_{0}|u_{0}|^{2}}{2}\right)  :=E_{0}.
\end{align*} Note that we assume that $E_{0}$ is small in Theorem
\ref{MainTheorem}.

\subsection{Extension of Hoff's estimates in a $L^{p}$ framework}

This part is the key of the paper: Assuming the initial velocity $u_{0}\in
H^{1}({\mathbb{T}}^{d})$\label{d=2_H1} and $\rho_{0}\in L^{2\gamma
}({\mathbb{T}}^{d})$, instead of $\rho_{0}\in L^{\infty}({\mathbb{T}}^{d})$ as
in \cite{hoff1995global}, we allow more general densities
than in \cite{hoff1995global}. This $L^{p}$, $p<\infty$
framework for the density is important when considering anisotropic viscous
tensors for which it is not so straightforward to propagate $L^{\infty}%
$-information. Consider\footnote{In this section $u_{\delta}:=\omega_{\delta
}\ast u$.}\footnote{We use the convention that
\[
\sum_{i=1}^{d}\sum_{k=1}^{d}\partial_{k}u^{i}\partial_{k}u^{i}=\partial
_{k}u^{i}\partial_{k}u^{i}:=\left\vert \partial_{k}u^{i}\right\vert ^{2}.
\]
}
\begin{equation}
A_{1}\left(  t\right)  =\frac{1}{2}\mu\int_{\mathbb{T}^{d}}\left\vert
\partial_{k}u^{i}\left(  t\right)  \right\vert ^{2}+\left(  \mu+\lambda
\right)  \int_{\mathbb{T}^{d}}\left\vert \operatorname{div}u\left(  t\right)
\right\vert ^{2}+\int_{\mathbb{T}^{d}}\varepsilon_{ijkl}\partial_{\ell
}u_{\delta}^{k}\left(  t\right)  \partial_{j}u_{\delta}^{i}\left(  t\right)
+\int_{0}^{t}\int_{\mathbb{T}^{d}}\rho\left\vert \dot{u}\right\vert ^{2},
\label{A1_Hoff}%
\end{equation}
and
\begin{equation}
A_{2}\left(  t\right)  =\sigma\left(  t\right)  \int_{\mathbb{T}^{d}}%
\frac{\rho\left(  t\right)  \left\vert \dot{u}\left(  t\right)  \right\vert
^{2}}{2}+\mu\int_{0}^{t}\int_{\mathbb{T}^{d}}\sigma\left\vert \partial_{k}%
\dot{u}^{i}\left(  t\right)  \right\vert ^{2}+\left(  \mu+\lambda\right)
\int_{0}^{t}\int_{\mathbb{T}^{d}}\sigma\left\vert \operatorname{div}\dot
{u}\right\vert ^{2}+\int_{0}^{t}\int_{\mathbb{T}^{d}}\sigma\varepsilon
_{ijkl}\partial_{\ell}\dot{u}_{\delta}^{k}\partial_{j}\dot{u}_{\delta}^{i}.
\label{A2_Hoff}%
\end{equation}
Multiplying the momentum equation with
\[
\dot{u}=u_{t}+u\cdot\nabla u
\]
we obtain (see the detailed computations in the appendix) that%
\begin{align}
A_{1}\left(  t\right)  =  &  \frac{1}{2}\mu\int_{\mathbb{T}^{d}}\left\vert
\partial_{k}u_{0}^{i}\right\vert ^{2}+\left(  \mu+\lambda\right)
\int_{\mathbb{T}^{d}}\left\vert \operatorname{div}u_{0}\right\vert ^{2}%
+\int_{\mathbb{T}^{d}}\varepsilon_{ijkl}\partial_{\ell}u_{0,\delta}%
^{k}\partial_{j}u_{0,\delta}^{i}\left(  t\right) \nonumber\\
&  -\mu\int_{0}^{t}\int_{\mathbb{T}^{d}}\partial_{k}u^{i}\partial_{k}u^{\ell
}\partial_{\ell}u^{i}+\frac{\mu}{2}\int_{0}^{t}\int_{\mathbb{T}^{d}}\left\vert
\partial_{k}u^{i}\right\vert ^{2}\operatorname{div}u\nonumber\\
&  -\left(  \mu+\lambda\right)  \int_{0}^{t}\int_{\mathbb{T}^{d}%
}\operatorname{div}u\partial_{i}u^{\ell}\partial_{\ell}u^{i}+\frac{\mu
+\lambda}{2}\int_{0}^{t}\int_{\mathbb{T}^{d}}(\operatorname{div}%
u)^{3}\nonumber\\
&  +\frac{1}{2}\int_{0}^{t}\int_{\mathbb{T}^{d}}\left\{  \partial
_{t}\varepsilon_{ijkl}+\partial_{q}(\varepsilon_{ijkl}u^{q})\right\}
\partial_{j}u_{\delta}^{i}\partial_{k}u_{\delta}^{\ell}-\int_{0}^{t}%
\int_{\mathbb{T}^{d}}\varepsilon_{ijkl}\partial_{\ell}u_{\delta}^{k}%
\omega_{\delta}\ast(\partial_{j}u^{q}\partial_{q}u^{i})\nonumber\\
&  -\int_{0}^{t}\int_{\mathbb{T}^{d}}\varepsilon_{ijkl}\partial_{\ell
}u_{\delta}^{k}\left[  u^{q},\omega_{\delta}\right]  \partial_{qj}^{2}%
u^{i}\nonumber\\
&  +\int_{0}^{t}\int_{\mathbb{T}^{d}}\rho P^{\prime}\left(  \rho\right)
\partial_{\ell}u^{k}\partial_{k}u^{\ell}+\int_{0}^{t}\int_{\mathbb{T}^{d}%
}(\rho P^{\prime}\left(  \rho\right)  -P\left(  \rho\right)
)(\operatorname{div}u)^{2}+\int_{0}^{t}\int_{\mathbb{T}^{d}}\rho\dot{u}f.
\label{ref_def_A1}%
\end{align}
Applying the operator $\partial_{t}+\operatorname{div}\left(  u\cdot\right)  $
to the momentum equation we obtain (see the detailed computations in the
appendix) that:%
\begin{align}
A_{2}\left(  t\right)   &  =\int_{0}^{1}\int_{\mathbb{T}^{d}}\sigma\frac
{\rho\left\vert \dot{u}\right\vert ^{2}}{2}+\mu\int_{0}^{t}\int_{\mathbb{T}%
^{d}}\sigma\partial_{k}u^{q}\partial_{q}u^{i}\partial_{k}\dot{u}^{i}+\mu
\int_{0}^{t}\int_{\mathbb{T}^{d}}\sigma\partial_{k}u^{q}\partial_{k}%
u^{i}\partial_{q}\dot{u}^{i}-\mu\int_{0}^{t}\int_{\mathbb{T}^{d}}%
\sigma\operatorname{div}u\partial_{k}u^{i}\partial_{k}\dot{u}^{i}\nonumber\\
&  +\left(  \mu+\lambda\right)  \int_{0}^{t}\int_{\mathbb{T}^{d}}%
\sigma\partial_{\ell}u^{q}\partial_{q}u^{\ell}\operatorname{div}\dot
{u}+\left(  \mu+\lambda\right)  \int_{0}^{t}\int_{\mathbb{T}^{d}}%
\sigma\partial_{i}u^{q}\partial_{q}\dot{u}^{i}\operatorname{div}u-\left(
\mu+\lambda\right)  \int_{0}^{t}\int_{\mathbb{T}^{d}}\sigma\left\vert
\operatorname{div}u\right\vert ^{2}\operatorname{div}\dot{u}\nonumber\\
&  -\int_{0}^{t}\int_{\mathbb{T}^{d}}\sigma\left(  \partial_{t}\varepsilon
_{ijkl}+\partial_{q}\left(  u^{q}\varepsilon_{ijkl}\right)  \right)
\partial_{\ell}u_{\delta}^{k}\partial_{j}\dot{u}_{\delta}^{i}-\int_{0}^{t}%
\int_{\mathbb{T}^{d}}\sigma\varepsilon_{ijkl}(\omega_{\delta}\ast
(\partial_{\ell}u^{q}\partial_{q}u^{k}))\partial_{j}\dot{u}_{\delta}%
^{i}\nonumber\\
&  -\int_{0}^{t}\int_{\mathbb{T}^{d}}\sigma\varepsilon_{ijkl}\partial_{\ell
}u_{\delta}^{k}\omega_{\delta}\ast(\partial_{j}u^{q}\partial_{q}\dot{u}%
^{i})\nonumber\\
&  +\int_{0}^{t}\int_{\mathbb{T}^{d}}\sigma\varepsilon_{ijkl}(\left[
u^{q},\omega_{\delta}\ast\right]  \partial_{\ell q}^{2}u^{k})\partial_{j}%
\dot{u}_{\delta}^{i}+\int_{0}^{t}\int_{\mathbb{T}^{d}}\sigma\varepsilon
_{ijkl}\partial_{\ell}u_{\delta}^{k}(\left[  u^{q},\omega_{\delta}\ast\right]
\partial_{jq}^{2}\dot{u}^{i})\nonumber\\
&  -\int_{0}^{t}\int_{\mathbb{T}^{d}}\sigma\left\{  P\left(  \rho\right)
\partial_{j}u^{k}\partial_{k}\dot{u}^{j}+\left(  \rho P^{\prime}\left(
\rho\right)  -P\left(  \rho\right)  \right)  \operatorname{div}%
u\operatorname{div}\dot{u}\right\}  . \label{ref_def_A2}%
\end{align}
Let us introduce the effective flux%
\[
F=\left(  2\mu+\lambda\right)  \operatorname{div}u-\left(  P\left(
\rho\right)  -P\left(  M\right)  \right)  .
\]
The details leading to these formulae are by now classic for the isotropic
case and they were used by D. Hoff in the series of works with isotropic
viscosities
\cite{hoff1995global,hoff1995strong,hoff2002dynamics,hoff2008lagrangean}.
Here, the added value is that these estimates are adapted for the anisotropic
approximate system $\left(  \text{\ref{ANS_delta}}\right)  $. As mentioned
before, for the reader's convenience we gather and detail these computations
in the Appendix. One of the key difficulties is to recover
information for the gradient of the velocity. A quick
analysis of $A_{1}$ and $A_{2}$ reveals that we need to control
\[
\int_{0}^{t}\left\Vert \nabla u\right\Vert _{L^{3}(\mathbb{T}^{d})}^{3}\text{
and }\int_{0}^{t}\sigma\left\Vert \nabla u\right\Vert _{L^{4}(\mathbb{T}^{d}%
)}^{4}.
\]
in order to close the estimates. The classical Calder\'{o}n-Zygmund theory
ensures that for $p\in\left\{  3,4\right\}  $ one has%
\[
\left\Vert \nabla u\right\Vert _{L^{p}(\mathbb{T}^{d})}^{p}\leq C\left(
\left\Vert \operatorname{curl}u\right\Vert _{L^{p}(\mathbb{T}^{d})}%
^{p}+\left\Vert \operatorname{div}u\right\Vert _{L^{p}(\mathbb{T}^{d})}%
^{p}\right)  ,
\]
for some numerical constant $C$. We deduce that%
\begin{align*}
\left\Vert \nabla u\right\Vert _{L^{p}(\mathbb{T}^{d})}^{p}  &  \leq C\left(
\frac{1}{\mu^{p}}\left\Vert \mu\operatorname{curl}u\right\Vert _{L^{p}%
(\mathbb{T}^{d})}^{p}+\frac{1}{(2\mu+\lambda)^{p}}\left\Vert (2\mu
+\lambda)\operatorname{div}u\right\Vert _{L^{p}(\mathbb{T}^{d})}^{p}\right) \\
&  \leq C\left(  \frac{1}{\mu^{p}}\left\Vert \mu\operatorname{curl}%
u\right\Vert _{L^{p}(\mathbb{T}^{d})}^{p}+\frac{1}{(2\mu+\lambda)^{p}%
}\left\Vert F\right\Vert _{L^{p}(\mathbb{T}^{d})}^{p}+\frac{1}{(2\mu
+\lambda)^{p}}\left\Vert P\left(  \rho\right)  -P\left(  M\right)  \right\Vert
_{L^{p}(\mathbb{T}^{d})}^{p}\right)  ,
\end{align*}
from which we infer that
\begin{align}
&  \int_{0}^{t}\left\Vert \nabla u\right\Vert _{L^{3}(\mathbb{T}^{d})}%
^{3}+\int_{0}^{t}\sigma\left\Vert \nabla u\right\Vert _{L^{4}(\mathbb{T}^{d}%
)}^{4}\nonumber\\
&  \leq C\left(  \frac{1}{\mu^{3}}\int_{0}^{t}\left\Vert \mu
\operatorname{curl}u\right\Vert _{L^{3}(\mathbb{T}^{d})}^{3}+\frac{1}%
{(2\mu+\lambda)^{3}}\int_{0}^{t}\left\Vert F\right\Vert _{L^{3}(\mathbb{T}%
^{d})}^{3}+\frac{1}{(2\mu+\lambda)^{3}}\int_{0}^{t}\left\Vert P\left(
\rho\right)  -P\left(  M\right)  \right\Vert _{L^{3}(\mathbb{T}^{d})}%
^{3}\right) \nonumber\\
&  +C\left(  \frac{1}{\mu^{4}}\int_{0}^{t}\sigma\left\Vert \mu
\operatorname{curl}u\right\Vert _{L^{4}(\mathbb{T}^{d})}^{4}+\frac{1}%
{(2\mu+\lambda)^{4}}\int_{0}^{t}\sigma\left\Vert F\right\Vert _{L^{4}%
(\mathbb{T}^{d})}^{4}+\frac{1}{(2\mu+\lambda)^{4}}\int_{0}^{t}\sigma\left\Vert
P\left(  \rho\right)  -P\left(  M\right)  \right\Vert _{L^{4}(\mathbb{T}^{d}%
)}^{4}\right)  . \label{u_grad_L4_1}%
\end{align}
Thus, in order to close the estimate, we have to recover control for the density.

\begin{remark}
This is where our approach starts to diverge from Hoff's approach. In the
isotropic setting, there is an extra algebraic structure which allows to
recover an $L^{\infty}$-bound for the density. In the anisotropic case, we
have to work with weaker norms, essentially because of the failure of
homogeneous Fourier multipliers of order $0$ to map $L^{\infty}$ to
$L^{\infty}$. The idea is to try only to propagate what seems to be necessary
to show that the two functionals $A_{1}$ and $A_{2}$ are bounded:%
\[
\int_{0}^{t}\left\Vert P\left(  \rho\right)  -P\left(  M\right)  \right\Vert
_{L^{3}(\mathbb{T}^{d})}^{3},\text{ }\int_{0}^{t}\sigma\left\Vert P\left(
\rho\right)  -P\left(  M\right)  \right\Vert _{L^{4}(\mathbb{T}^{d})}^{4}.
\]

\end{remark}

\subsubsection{Bounds for the density}

In the following lines we want to obtain estimates for the density. We begin
by arranging $\left(  \text{\ref{renormm}}\right)  $ as%
\begin{align}
&  \partial_{t}b\left(  \rho\right)  +\operatorname{div}\left(  b\left(
\rho\right)  u\right)  +\left(  \rho b^{\prime}\left(  \rho\right)  -b\left(
\rho\right)  \right)  \frac{\left(  P\left(  \rho\right)  -P\left(  M\right)
\right)  }{2\mu+\lambda}\nonumber\\
&  =-\frac{1}{2\mu+\lambda}\left(  \rho b^{\prime}\left(  \rho\right)
-b\left(  \rho\right)  \right)  \left(  \left(  2\mu+\lambda\right)
\operatorname{div}u-(P\left(  \rho\right)  -P\left(  M\right)  )\right)  .
\label{renorm_9}%
\end{align}
Recall the definitions of $H_{2}\left(  \cdot/M\right)  $ and $H_{3}\left(
\cdot/M\right)  $ given in $\left(  \text{\ref{Hshurile}}\right)  $.

\noindent\textit{A $L^{3}$ control for the pressure.} Let us take
$b=H_{2}\left(  \cdot/M\right)  $ in $\left(  \text{\ref{renorm_9}}\right)  $
with
\[
\rho H_{2}^{\prime}\left(  \rho/M\right)  -H_{2}\left(  \rho/M\right)
=\left\vert P\left(  \rho\right)  -P\left(  M\right)  \right\vert \left(
P\left(  \rho\right)  -P\left(  M\right)  \right)  ,
\]
and use Young's inequality in order to obtain that%
\begin{align}
&  \int_{\mathbb{T}^{d}}H_{2}\left(  \rho(t)/M\right)  +\frac{1}%
{2(2\mu+\lambda)}\int_{0}^{t}\int_{\mathbb{T}^{d}}\left\vert P\left(
\rho\right)  -P\left(  M\right)  \right\vert ^{3}\nonumber\\
&  \leq\int_{\mathbb{T}^{d}}H_{2}\left(  \rho_{0}/M\right)  +\frac{1}%
{2\mu+\lambda}\int_{0}^{t}\int_{\mathbb{T}^{d}}\left\vert (2\mu+\lambda
)\operatorname{div}u-\left(  P\left(  \rho\right)  -P\left(  M\right)
\right)  \right\vert ^{3}. \label{bound_H2}%
\end{align}

\medskip

\noindent\textit{A $L^{4}$ control of the pressure.} Finally, take
$b=H_{3}\left(  \cdot/M\right)  $ with%
\[
\rho H_{3}^{\prime}\left(  \rho/M\right)  -H_{3}\left(  \rho/M\right)
=\left(  P\left(  \rho\right)  -P\left(  M\right)  \right)  ^{3}%
\]
in order to obtain that%
\begin{align}
\sigma\left(  t\right)   &  \int_{\mathbb{T}^{d}}H_{3}\left(  \rho
(t)/M\right)  +\frac{1}{2(2\mu+\lambda)}\int_{0}^{t}\int_{\mathbb{T}^{d}%
}\sigma\left\vert P\left(  \rho\right)  -P\left(  M\right)  \right\vert
^{4}\nonumber\\
&  \leq
\int_{0}^{1}\int_{\mathbb{T}^{d}}H_{3}\left( \rho(s)/M\right)ds+\frac
{1}{2\mu+\lambda}\int_{0}^{t}\int_{\mathbb{T}^{d}}\sigma\left\vert
(2\mu+\lambda)\operatorname{div}u-\left(  P\left(  \rho\right)  -P\left(
M\right)  \right)  \right\vert ^{4} \label{bound_H3}%
\end{align}

\noindent A simple computation gives us:%
\begin{align*}
H_{3}\left(  \rho/M\right)   &  =\frac{1}{3\gamma-1}\left(  \rho^{3\gamma
}-M^{3\gamma-1}\rho\right)  -\frac{3M^{\gamma}}{2\gamma-1}\left(
\rho^{2\gamma}-M^{2\gamma-1}\rho\right) \\
&  +\frac{3M^{2\gamma}}{\gamma-1}\left(  \rho^{\gamma}-M^{\gamma-1}%
\rho\right)  +M^{3\gamma}-\rho M^{3\gamma-1}.
\end{align*}
Then, we observe that%
\begin{align}
\int_{0}^{1}\int_{\mathbb{T}^{3}}H_{3}\left(  \rho/M\right)   &  =\frac
{1}{3\gamma-1}\int_{0}^{1}\int_{\mathbb{T}^{d}}\left(  \rho^{3\gamma
}-M^{3\gamma}\right)  -\frac{3M^{\gamma}}{2\gamma-1}\int_{0}^{1}%
\int_{\mathbb{T}^{d}}\left(  \rho^{2\gamma}-M^{2\gamma}\right)  +3M^{2\gamma
}\int_{0}^{1}\int_{\mathbb{T}^{d}}H_{1}\left(  \rho/M\right) \nonumber\\
&  \leq\frac{1}{2\gamma-1}\int_{0}^{1}\int_{\mathbb{T}^{d}}\left(
\rho^{3\gamma}-M^{3\gamma}\right)  -\frac{3M^{\gamma}}{2\gamma-1}\int_{0}%
^{1}\int_{\mathbb{T}^{d}}\left(  \rho^{2\gamma}-M^{2\gamma}\right)
+3M^{2\gamma}E_{0}\nonumber\\
&  =\frac{1}{2\gamma-1}\int_{0}^{1}\int_{\mathbb{T}^{d}}\left(  \rho^{\gamma
}-M^{\gamma}\right)  ^{3}-\frac{3\left(  \gamma-1\right)  M^{2\gamma}}%
{2\gamma-1}\int_{0}^{1}\int_{\mathbb{T}^{d}}H_{1}\left(  \rho/M\right)
+3M^{2\gamma}E_{0}\nonumber\\
&  \leq\frac{1}{2\gamma-1}\int_{0}^{1}\int_{\mathbb{T}^{d}}\left\vert P\left(
\rho\right)  -P\left(  M\right)  \right\vert ^{3}+3M^{2\gamma}E_{0},
\label{H3_comput_exp}%
\end{align}
where, we used that%
\[
\int_{0}^{1}\int_{\mathbb{T}^{d}}\left(  \rho^{3\gamma}-M^{3\gamma}\right)
=\int_{0}^{1}\int_{\mathbb{T}^{d}}\rho^{3\gamma}-M^{3\gamma}-3\gamma
M^{3\gamma-1}\left(  \rho-M\right)  \geq0,
\]
owing to the convexity of $s\rightarrow s^{3\gamma}$. Thus, since $\gamma
\geq1$ we obtain that%

\begin{equation}
\int_{0}^{1}\int_{\mathbb{T}^{3}}H_{3}\left(  \rho/M\right)  \leq3M^{2\gamma
}E_{0}+\int_{0}^{1}\int_{\mathbb{T}^{d}}\left\vert P\left(  \rho\right)
-P\left(  M\right)  \right\vert ^{3}. \label{estimate_H3_temp_court}%
\end{equation}
Using the above estimate, $\left(  \text{\ref{bound_H3}}\right)  $ and
$\left(  \text{\ref{bound_H2}}\right)  $ we obtain that%
\begin{align}
\sigma\left(  t\right)   &  \int_{\mathbb{T}^{d}}H_{3}\left(  \rho
(t)/M\right)  +\frac{1}{2(2\mu+\lambda)}\int_{0}^{t}\int_{\mathbb{T}^{d}%
}\sigma\left\vert P\left(  \rho\right)  -P\left(  M\right)  \right\vert
^{4}\nonumber\\
&  \leq3M^{2\gamma}E_{0}+\int_{\mathbb{T}^{d}}H_{2}\left(  \rho_{0}/M\right)
+\int_{0}^{1}\int_{\mathbb{T}^{d}}\left\vert P\left(  \rho)-P(M)\right)
\right\vert ^{3}+\frac{C}{2\mu+\lambda}\int_{0}^{t}\int_{\mathbb{T}^{d}}%
\sigma\left\vert F\right\vert ^{4}. \label{Bound_H3}%
\end{align}
Let us combine $\left(  \text{\ref{bound_H2}}\right)  $ with $\left(
\text{\ref{Bound_H3}}\right)  $ to deduce that:%
\begin{align}
B\left(  t\right)   &  :=\left(  1+2\left(  2\mu+\lambda\right)  \right)
\int_{\mathbb{T}^{d}}H_{2}\left(  \rho(t)/M\right)  +\sigma\left(  t\right)
\int_{\mathbb{T}^{d}}H_{3}\left(  \rho(t)/M\right) \nonumber\\
&  +\frac{1}{2(2\mu+\lambda)}\int_{0}^{t}\int_{\mathbb{T}^{d}}\left\vert
P\left(  \rho\right)  -P\left(  M\right)  \right\vert ^{3}+\frac{1}%
{2(2\mu+\lambda)}\int_{0}^{t}\int_{\mathbb{T}^{d}}\sigma\left\vert P\left(
\rho\right)  -P\left(  M\right)  \right\vert ^{4}\nonumber\\
&  \leq\left(  1+2\left(  2\mu+\lambda\right)  \right)  \int_{\mathbb{T}^{d}%
}H_{2}\left(  \rho_{0}/M\right)  +3M^{2\gamma}E_{0}\nonumber\\
&  +\frac{1+2\left(  2\mu+\lambda\right)  }{2\mu+\lambda}\int_{0}%
^{t}\left\Vert F\right\Vert _{L^{3}({\mathbb{T}}^{d})}^{3}+\frac{1}%
{2\mu+\lambda}\int_{0}^{t}\sigma\left\Vert F\right\Vert _{L^{4}({\mathbb{T}%
}^{d})}^{4}. \label{definition_B}%
\end{align}

\subsubsection{Bounds for the Hoff functionals}

We recall that%
\[
\left\{
\begin{array}
[c]{l}%
\mu\operatorname{curl}u=\Delta^{-1}\operatorname{curl}(\rho\dot{u}%
)+\Delta^{-1}\operatorname{curl}\left(  \operatorname{div}(\omega_{\delta}%
\ast(\mathcal{E}\nabla(\omega_{\delta}\ast u)))\right)  ,\\
F=\Delta^{-1}\operatorname{div}\left(  \rho\dot{u}\right)  +\Delta
^{-1}\operatorname{div}\left(  \operatorname{div}(\omega_{\delta}%
\ast(\mathcal{E}\nabla(\omega_{\delta}\ast u)))\right)
\end{array}
\right.
\]
and therefore, using also $\left(  \text{\ref{bound_H2}}\right)  $ we have
that%
\begin{align}
&  \int_{0}^{t}\left\Vert \nabla u\right\Vert _{L^{3}(\mathbb{T}^{d})}^{3}%
\leq\frac{C}{\left(  2\mu+\lambda\right)  ^{3}}\int_{0}^{t}\left\Vert P\left(
\rho\right)  -P\left(  M\right)  \right\Vert _{L^{3}(\mathbb{T}^{d})}%
^{3}+\frac{C}{\mu^{3}}\int_{0}^{t}\left\Vert \mu\operatorname{curl}%
u\right\Vert _{L^{3}(\mathbb{T}^{d})}^{3}+\frac{C}{\left(  2\mu+\lambda
\right)  ^{3}}\int_{0}^{t}\left\Vert F\right\Vert _{L^{3}(\mathbb{T}^{d})}%
^{3}\\
&  \leq\frac{C}{\left(  2\mu+\lambda\right)  ^{2}}\int_{\mathbb{T}^{d}}%
H_{2}\left(  \rho_{0}/M\right)  +C\int_{0}^{t}\left\Vert \Delta^{-1}%
\operatorname{div}(\rho\dot{u})\right\Vert _{L^{3}(\mathbb{T}^{d})}^{3}%
+C\int_{0}^{t}\left\Vert \Delta^{-1}\operatorname{curl}(\rho\dot
{u})\right\Vert _{L^{3}(\mathbb{T}^{d})}^{3}\\
&  \text{ \ \ \ \ \ \ \ \ \ \ \ \ \ \ \ \ \ \ \ \ \ \ \ \ \ \ \ \ \ \ \ }%
+C\left\Vert \mathcal{E-I}\right\Vert _{L^{\infty}((0,t)\times\mathbb{T}^{d}%
)}^{3}\max\left\{  \frac{1}{\left(  2\mu+\lambda\right)  ^{3}},\frac{1}%
{\mu^{3}}\right\}  \int_{0}^{t}\left\Vert \nabla u\right\Vert _{L^{3}%
(\mathbb{T}^{d})}^{3} \label{grad_u_L3_1}%
\end{align}
where $C$ is a generic constant that depends only on the dimension
$d\in\left\{  2,3\right\}  $ whose exact value can change from one line to the
other. Let us observe that using interpolation and Sobolev imbedding
inequalities, we obtain that:%

\begin{align}
&  \int_{0}^{t}\left\Vert \Delta^{-1}\operatorname{div}(\rho\dot
{u})\right\Vert _{L^{3}(\mathbb{T}^{d})}^{3}+\int_{0}^{t}\left\Vert
\Delta^{-1}\operatorname{curl}(\rho\dot{u})\right\Vert _{L^{3}(\mathbb{T}%
^{d})}^{3}\nonumber\\
&  \leq\int_{0}^{t}\left\Vert \Delta^{-1}\operatorname{div}(\rho\dot
{u})\right\Vert _{L^{2}(\mathbb{T}^{d})}\left\Vert \Delta^{-1}%
\operatorname{div}(\rho\dot{u})\right\Vert _{L^{4}(\mathbb{T}^{d})}^{2}%
+\int_{0}^{t}\left\Vert \Delta^{-1}\operatorname{curl}(\rho\dot{u})\right\Vert
_{L^{2}(\mathbb{T}^{d})}\left\Vert \Delta^{-1}\operatorname{curl}(\rho\dot
{u})\right\Vert _{L^{4}(\mathbb{T}^{d})}^{2}\nonumber\\
&  \leq C\left(  \sup_{t>0}\left\Vert \Delta^{-1}\operatorname{div}(\rho
\dot{u})\right\Vert _{L^{2}(\mathbb{T}^{d})}+\sup_{t}\left\Vert \Delta
^{-1}\operatorname{curl}(\rho\dot{u})\right\Vert _{L^{2}(\mathbb{T}^{d}%
)}\right)  \sup_{t>0}\left\Vert \sqrt{\rho\left(  t\right)  }\right\Vert
_{L^{4d/(4-d)}(\mathbb{T}^{d})}^{2}\int_{0}^{t}\left\Vert \sqrt{\rho}\dot
{u}\right\Vert _{L^{2}(\mathbb{T}^{d})}^{2}. \label{grad_u_L3_2}%
\end{align}

In order to close the estimate, we need the following:

\begin{lemma}
\label{technical_lemma} For all $\rho\geq0$ we have that:
\begin{equation}
\left(  P\left(  \rho\right)  -P\left(  M\right)  \right)  ^{2}\leq2\gamma
M^{\gamma}H_{1}\left(  \rho/M\right)  +\left(  2\gamma-1\right)  H_{2}\left(
\rho/M\right)  . \label{ineg2}%
\end{equation}

\end{lemma}

\noindent\textit{Proof of Lemma} \ref{technical_lemma}: Consider
$g:(0,+\infty)\rightarrow\mathbb{R}$ defined by:
\[
g\left(  \rho\right)  =\alpha\frac{H_{1}\left(  \rho/M\right)  }{\rho}%
+\beta\frac{H_{2}\left(  \rho/M\right)  }{\rho}-\frac{\left(  P\left(
\rho\right)  -P\left(  M\right)  \right)  ^{2}}{\rho}%
\]
and observe that%
\[
g\left(  M\right)  =0.
\]
We have that%
\[
\rho^{2}g^{\prime}\left(  \rho\right)  =\left(  \alpha-2\gamma M^{\gamma
}\right)  \left(  \rho^{\gamma}-M^{\gamma}\right)  +\beta\left\vert
\rho^{\gamma}-M^{\gamma}\right\vert \left(  \rho^{\gamma}-M^{\gamma}\right)
-\left(  2\gamma-1\right)  \left(  \rho^{\gamma}-M^{\gamma}\right)  ^{2}.
\]
For $\alpha=2\gamma M^{\gamma}$ and $\beta=2\gamma-1$ we see that%
\[
g^{\prime}\left(  \rho\right)  \leq0\text{ on }\rho\in\lbrack0,M]\text{ and
}g^{\prime}\left(  \rho\right)  \geq0\text{ if }\rho\geq M.
\]
This ends the proof of Lemma \ref{technical_lemma}.

Using \eqref{ineg2} in Lemma \ref{technical_lemma} we have that%
\[
\left\Vert P\left(  \rho\right)  -P\left(  M\right)  \right\Vert _{L^{\infty
}((0,T)\times L^{2}(\mathbb{T}^{d}))}\leq2\gamma M^{\gamma}E_{0}+\left(
2\gamma-1\right)  B\left(  t\right)
\]
Also, remark that owing to $\gamma\geq d/(4-d)$ we infer that%
\begin{align}
\left\Vert \sqrt{\rho}\right\Vert _{L^{4d/((4-d)}(\mathbb{T}^{d})}^{4}  &
=\left\Vert \rho\right\Vert _{L^{2d/(4-d)}(\mathbb{T}^{d})}^{2}\leq\left\Vert
P(\rho)\right\Vert _{L^{2}(\mathbb{T}^{d})}^{2/\gamma}\nonumber\\
&  \leq\left\{  M^{\gamma}+\int_{\mathbb{T}^{d}}\left\{  2\gamma M^{\gamma
}H_{1}\left(  \rho/M\right)  +\left(  2\gamma-1\right)  H_{2}\left(
\rho/M\right)  \right\}  \right\}  ^{\frac{2}{\gamma}}\nonumber\\
&  \leq C_{\mu,\lambda,\gamma,M}(M+B\left(  t\right)  ). \label{inegalitate}%
\end{align}
where $C_{\mu,\lambda,\gamma,M}$ is a generic constant that depends only on
the $\mu,\lambda,\gamma,M$ \ and whose exact value can change from one line to
the other. It is of course, different from the generic constant $C$ appearing
in $\left(  \text{\ref{grad_u_L3_1}}\right)  $ and $\left(
\text{\ref{grad_u_L3_2}}\right)  $ that only depends on the dimension.
Combining $\left(  \text{\ref{grad_u_L3_1}}\right)  $ with $\left(
\text{\ref{inegalitate}}\right)  $ we obtain:%
\[
\sup_{t>0}\left\Vert \Delta^{-1}\operatorname{div}(\rho\dot{u})\right\Vert
_{L^{2}(\mathbb{T}^{d})}+\sup_{t}\left\Vert \Delta^{-1}\operatorname{curl}%
(\rho\dot{u})\right\Vert _{L^{2}(\mathbb{T}^{d})}\leq C_{\mu,\lambda,\gamma
,M}\left(  E_{0}+B\left(  t\right)  +\sqrt{A_{1}\left(  t\right)  }\right)  .
\]

Thus, using $\left(  \text{\ref{Hyp4}}\right)  $, the last term of the RHS can
be absorbed into the LHS thus giving%
\begin{equation}
\int_{0}^{t}\left\Vert \nabla u\right\Vert _{L^{3}(\mathbb{T}^{d})}^{3}%
\leq\frac{C}{\left(  2\mu+\lambda\right)  ^{2}}\int_{\mathbb{T}^{d}}%
H_{2}\left(  \rho_{0}/M\right)  +C_{\mu,\lambda,\gamma,M}\left(
E_{0}+B\left(  t\right)  +\sqrt{A_{1}\left(  t\right)  }\right)  A_{1}\left(
t\right)  . \label{grad_u_L3}%
\end{equation}

Similarly, using $\left(  \text{\ref{Bound_H3}}\right)  $ we have that%
\begin{align}
&  \int_{0}^{t}\sigma\left\Vert \nabla u\right\Vert _{L^{4}(\mathbb{T}^{d}%
)}^{4}\nonumber\\
&  \leq\frac{C}{\left(  2\mu+\lambda\right)  ^{4}}\int_{0}^{t}\sigma\left\Vert
P\left(  \rho\right)  -P\left(  M\right)  \right\Vert _{L^{4}(\mathbb{T}^{d}%
)}^{4}+\frac{C}{\mu^{4}}\int_{0}^{t}\sigma\left\Vert \mu\operatorname{curl}%
u\right\Vert _{L^{4}(\mathbb{T}^{d})}^{4}+\frac{C}{\left(  2\mu+\lambda
\right)  ^{4}}\int_{0}^{t}\sigma\left\Vert F\right\Vert _{L^{4}(\mathbb{T}%
^{d})}^{4}\nonumber\\
&  \leq\frac{3M^{2\gamma}E_{0}}{\left(  2\mu+\lambda\right)  ^{3}}+\frac
{1}{\left(  2\mu+\lambda\right)  ^{3}}\int_{\mathbb{T}^{d}}H_{2}\left(
\rho_{0}/M\right)  +\frac{C}{\left(  2\mu+\lambda\right)  ^{3}}\int_{0}%
^{1}\left\Vert P\left(  \rho\right)  -P\left(  M\right)  \right\Vert
_{L^{3}(\mathbb{T}^{d})}^{3}\nonumber\\
&  +\frac{C}{\mu^{4}}\int_{0}^{t}\sigma\left\Vert \mu\operatorname{curl}%
u\right\Vert _{L^{4}(\mathbb{T}^{d})}^{4}+\frac{C}{\left(  2\mu+\lambda
\right)  ^{4}}\int_{0}^{t}\sigma\left\Vert F\right\Vert _{L^{4}(\mathbb{T}%
^{d})}^{4}. \label{grad_l4_0}%
\end{align}
Next, we see that%
\begin{align}
&  \frac{C}{\mu^{4}}\int_{0}^{t}\sigma\left\Vert \mu\operatorname{curl}%
u\right\Vert _{L^{4}(\mathbb{T}^{d})}^{4}+\frac{C}{\left(  2\mu+\lambda
\right)  ^{4}}\int_{0}^{t}\sigma\left\Vert F\right\Vert _{L^{4}(\mathbb{T}%
^{d})}^{4}\nonumber\\
&  \leq\frac{C}{\mu^{4}}\int_{0}^{t}\sigma\left\Vert \sqrt{\rho}\right\Vert
_{L^{4d/(4-d)}(\mathbb{T}^{d})}^{4}\left\Vert \sqrt{\rho}\dot{u}\right\Vert
_{L^{2}(\mathbb{T}^{d})}^{4}+C\left\Vert \mathcal{E-I}\right\Vert _{L^{\infty
}((0,t)\times\mathbb{T}^{d})}^{4}\max\left\{  \frac{1}{\left(  2\mu
+\lambda\right)  ^{4}},\frac{1}{\mu^{4}}\right\}  \int_{0}^{t}\sigma\left\Vert
\nabla u\right\Vert _{L^{4}(\mathbb{T}^{d})}^{4}\nonumber\\
&  \leq\frac{C}{\mu^{4}}\sup_{t}\left\Vert \sqrt{\rho}\right\Vert
_{L^{4d/(4-d)}(\mathbb{T}^{d})}^{4}\sup_{t}\sigma\left\Vert \sqrt{\rho}\dot
{u}\right\Vert _{L^{2}(\mathbb{T}^{d})}^{2}\int_{0}^{t}\left\Vert \sqrt{\rho
}\dot{u}\right\Vert _{L^{2}(\mathbb{T}^{d})}^{2}\nonumber\\
&  \text{
\ \ \ \ \ \ \ \ \ \ \ \ \ \ \ \ \ \ \ \ \ \ \ \ \ \ \ \ \ \ \ \ \ \ \ \ \ \ \ \ \ \ \ \ \ \ \ \ \ \ \ }%
+C\left\Vert \mathcal{E-I}\right\Vert _{L^{\infty}((0,t)\times\mathbb{T}^{d}%
)}^{4}\max\left\{  \frac{1}{\left(  2\mu+\lambda\right)  ^{4}},\frac{1}%
{\mu^{4}}\right\}  \int_{0}^{t}\sigma\left\Vert \nabla u\right\Vert
_{L^{4}(\mathbb{T}^{d})}^{4}\nonumber\\
&  \leq\frac{C}{\mu^{4}}\sup_{t}\left\Vert \sqrt{\rho}\right\Vert
_{L^{4d/(4-d)}(\mathbb{T}^{d})}^{4}A_{1}\left(  t\right)  A_{2}\left(
t\right)  +C\left\Vert \mathcal{E-I}\right\Vert _{L^{\infty}((0,t)\times
\mathbb{T}^{d})}^{4}\max\left\{  \frac{1}{\left(  2\mu+\lambda\right)  ^{4}%
},\frac{1}{\mu^{4}}\right\}  \int_{0}^{t}\sigma\left\Vert \nabla u\right\Vert
_{L^{4}(\mathbb{T}^{d})}^{4} \label{grad_l4}%
\end{align}
Of course, under Hypothesis $\left(  \text{\ref{Hyp4}}\right)  $ we see that
the last term from the RHS above inequality can be absorbed into the LHS. From
the above estimtes, it transpires that $\left\Vert P\left(  \rho\right)
-P\left(  M\right)  \right\Vert _{L^{3}((0,1)\times\mathbb{T}^{d})}$ verifies
the same estimate $\left(  \text{\ref{grad_u_L3}}\right)  $. Combining the
estimates $\left(  \text{\ref{inegalitate}}\right)  $, $\left(
\text{\ref{grad_l4_0}}\right)  $ and $\left(  \text{\ref{grad_l4}}\right)  $
we end up with:
\begin{align*}
\int_{0}^{t}\sigma\left\Vert \nabla u\right\Vert _{L^{4}(\mathbb{T}^{d})}^{4}
&  \leq\frac{M^{2\gamma}E_{0}}{\left(  2\mu+\lambda\right)  ^{3}}+\frac
{C}{\left(  2\mu+\lambda\right)  ^{3}}\int_{\mathbb{T}^{d}}H_{2}\left(
\rho_{0}/M\right) \\
&  +C_{\mu,\lambda,\gamma,M}(M+B\left(  t\right)  )\left(  E_{0}+B\left(
t\right)  +\sqrt{A_{1}\left(  t\right)  }+A_{2}\left(  t\right)  \right)
A_{1}\left(  t\right) \\
&  +C\left\Vert \mathcal{E-I}\right\Vert _{L^{\infty}((0,t)\times
\mathbb{T}^{d})}^{4}\max\left\{  \frac{1}{\left(  2\mu+\lambda\right)  ^{4}%
},\frac{1}{\mu^{4}}\right\}  \int_{0}^{t}\sigma\left\Vert \nabla u\right\Vert
_{L^{4}(\mathbb{T}^{d})}^{4},
\end{align*}
where, we recall that here, $C$ depends only on the dimension. Thus, under the
hypothesis $\left(  \text{\ref{Hyp4}}\right)  $ we obtain that
\begin{align}
&  \int_{0}^{t}\left\Vert \nabla u\right\Vert _{L^{3}(\mathbb{T}^{d})}%
^{3}+\int_{0}^{t}\sigma\left\Vert \nabla u\right\Vert _{L^{4}(\mathbb{T}^{d}%
)}^{4}\nonumber\\
&  \leq C_{\mu,\lambda,\gamma,M}\left(  E_{0}+\int_{\mathbb{T}^{d}}%
H_{2}\left(  \rho_{0}/M\right)  \right)  +C_{\mu,\lambda,\gamma,M}\left(
M+B\left(  t\right)  \right)  \left(  E_{0}+B\left(  t\right)  +\sqrt
{A_{1}\left(  t\right)  }+A_{2}\left(  t\right)  \right)  A_{1}(t).
\label{grad_u(L3+L4)}%
\end{align}

\begin{remark}
In this section we used two generic constants, $C$ and $C_{\mu,\lambda
,\gamma,M}$ in order to emphasize that the anisotropic amplitude is small only
with respect $\mu$ and $2\mu+\lambda,$ see Hypothesis $\left(
\text{\ref{Hyp4}}\right)  $. In the estimates that follow in the next
sections, $C$ will denote a generic constant that depends on the parameters of
the problem $C=C\left(  \mu,\lambda,\gamma,M\right)  $, the exact value of
which can change from a line to another.
\end{remark}

\paragraph{Terms appearing in the RHS of $\left(  \ref{ref_def_A1}\right)  $.}

Recall that%
\begin{align}
A_{1}\left(  t\right)   &  =\frac{\mu}{2}\int_{\mathbb{T}^{d}}\left\vert
\partial_{k}u_{0}^{i}\right\vert ^{2}+\frac{\mu+\lambda}{2}\int_{\mathbb{T}%
^{d}}\left\vert \operatorname{div}u_{0}\right\vert ^{2}+\int_{\mathbb{T}^{d}%
}\varepsilon_{ijkl}\partial_{\ell}u_{0,\delta}^{k}\partial_{j}u_{0,\delta}%
^{i}\nonumber\\
&  \int_{\mathbb{T}^{d}}P\left(  \rho\left(  t\right)  \right)
\operatorname{div}u\left(  t\right)  -\int_{\mathbb{T}^{d}}P\left(
\rho\left(  0\right)  \right)  \operatorname{div}u\left(  0\right) \nonumber\\
&  -\mu\int_{0}^{t}\int_{\mathbb{T}^{d}}\partial_{k}u^{i}\partial_{k}u^{\ell
}\partial_{\ell}u^{i}+\frac{\mu}{2}\int_{0}^{t}\int_{\mathbb{T}^{d}}\left\vert
\partial_{k}u^{i}\right\vert ^{2}\operatorname{div}u\nonumber\\
&  -\left(  \mu+\lambda\right)  \int_{0}^{t}\int_{\mathbb{T}^{d}%
}\operatorname{div}u\partial_{i}u^{\ell}\partial_{\ell}u^{i}+\frac{\mu
+\lambda}{2}\int_{0}^{t}\int_{\mathbb{T}^{d}}(\operatorname{div}%
u)^{3}\nonumber\\
&  +\frac{1}{2}\int_{0}^{t}\int_{\mathbb{T}^{d}}\left\{  \partial
_{t}\varepsilon_{ijkl}+\partial_{q}(\varepsilon_{ijkl}u^{q})\right\}
\partial_{j}u_{\delta}^{i}\partial_{k}u_{\delta}^{\ell}\nonumber\\
&  -\int_{0}^{t}\int_{\mathbb{T}^{d}}\varepsilon_{ijkl}\partial_{\ell
}u_{\delta}^{k}\omega_{\delta}\ast(\partial_{j}u^{q}\partial_{q}u^{i}%
)-\int_{0}^{t}\int_{\mathbb{T}^{d}}\varepsilon_{ijkl}\partial_{\ell}u_{\delta
}^{k}\left[  u^{q},\omega_{\delta}\right]  \partial_{qj}^{2}u^{i}\nonumber\\
&  +\int_{0}^{t}\int_{\mathbb{T}^{d}}\rho P^{\prime}\left(  \rho\right)
\partial_{\ell}u^{k}\partial_{k}u^{\ell}+\int_{\mathbb{T}^{d}}\left(  \rho
P^{\prime}\left(  \rho\right)  -P\left(  \rho\right)  \right)
(\operatorname{div}u)^{2}.
\end{align}
First, using $\left(  \text{\ref{technical_lemma}}\right)  $ we have that%
\begin{align*}
\int_{\mathbb{T}^{d}}P\left(  \rho\left(  t\right)  \right)
\operatorname{div}u\left(  t\right)   &  =\int_{\mathbb{T}^{d}}\left(
P\left(  \rho\left(  t\right)  \right)  -P\left(  M\right)  \right)
\operatorname{div}u\left(  t\right) \\
&  \leq C\left(  \eta\right)  \int_{\mathbb{T}^{d}}\left(  P\left(
\rho\left(  t\right)  \right)  -P\left(  M\right)  \right)  ^{2}+\eta
\int_{\mathbb{T}^{d}}\left\vert \operatorname{div}u\right\vert ^{2}\left(
t\right) \\
&  \leq C\left(  \eta\right)  \left\{  \int_{\mathbb{T}^{d}}H_{1}\left(
\rho\right)  +H_{2}\left(  \rho\right)  \right\}  +\eta\int_{\mathbb{T}^{d}%
}\left\vert \operatorname{div}u\right\vert ^{2}\left(  t\right)  ,
\end{align*}
where $\eta$ will be chosen later. Using the last estimate, we obtain that%
\begin{align}
&  \int_{\mathbb{T}^{d}}P\left(  \rho\left(  t\right)  \right)
\operatorname{div}u\left(  t\right)  -\int_{\mathbb{T}^{d}}P\left(
\rho\left(  0\right)  \right)  \operatorname{div}u\left(  0\right) \nonumber\\
&  +\frac{\mu}{2}\int_{\mathbb{T}^{d}}\left\vert \partial_{k}u_{0}%
^{i}\right\vert ^{2}+\frac{\mu+\lambda}{2}\int_{\mathbb{T}^{d}}\left\vert
\operatorname{div}u_{0}\right\vert ^{2}+\frac{1}{2}\int_{\mathbb{T}^{d}%
}\varepsilon_{ijkl}\partial_{\ell}u_{0,\delta}^{k}\partial_{j}u_{0,\delta}%
^{i}+\frac{1}{2}\int_{\mathbb{T}^{d}}\partial_{t}\varepsilon_{ijkl}%
\partial_{j}u_{\delta}^{i}\partial_{k}u_{\delta}^{\ell}\nonumber\\
&  \leq C\left(  \eta\right)  \left\{  E_{0}+\left\Vert \nabla u_{0}%
\right\Vert _{L^{2}(\mathbb{T}^{d})}+B\left(  t\right)  \right\}  +\eta
\int_{\mathbb{T}^{d}}\left\vert \nabla u\right\vert ^{2}\left(  t\right)
\label{d=2H1A_1_1}%
\end{align}

Using $\left(  \text{\ref{grad_u_L3}}\right)  $ and hypothesis $\left(
\text{\ref{Hyp3}}\right)  $ we infer that%
\begin{align}
&  \int_{0}^{t}\int_{\mathbb{T}^{d}}\rho P^{\prime}\left(  \rho\right)
\partial_{\ell}u^{k}\partial_{k}u^{\ell}=\gamma P\left(  M\right)  \int
_{0}^{t}\int_{\mathbb{T}^{d}}\partial_{\ell}u^{k}\partial_{k}u^{\ell}%
+\gamma\int_{0}^{t}\int_{\mathbb{T}^{d}}\left(  P\left(  \rho\right)
-P\left(  M\right)  \right)  \partial_{\ell}u^{k}\partial_{k}u^{\ell
}\nonumber\\
&  \leq CE_{0}+\int_{0}^{t}\int_{\mathbb{T}^{d}}\left(  P\left(  \rho\right)
-P\left(  M\right)  \right)  ^{3}+\int_{0}^{t}\int_{\mathbb{T}^{d}}\left\vert
\nabla u\right\vert ^{3}. \label{d=2H1A_1_2}%
\end{align}
The term $\int_{0}^{t}\int_{\mathbb{T}^{d}}\left(  \rho P^{\prime}\left(
\rho\right)  -P\left(  \rho\right)  \right)  (\operatorname{div}u)^{2}$ is
treated similarly.

Using Proposition \ref{Prop_ren1} in order to trat the comutator term, we have
that:
\begin{align}
&  -\mu\int_{0}^{t}\int_{\mathbb{T}^{d}}\partial_{k}u^{i}\partial_{k}u^{\ell
}\partial_{\ell}u^{i}+\frac{\mu}{2}\int_{0}^{t}\int_{\mathbb{T}^{d}}\left\vert
\partial_{k}u^{i}\right\vert ^{2}\operatorname{div}u-\left(  \mu
+\lambda\right)  \int_{0}^{t}\int_{\mathbb{T}^{d}}\operatorname{div}%
u\partial_{i}u^{\ell}\partial_{\ell}u^{i}+\frac{\mu+\lambda}{2}\int_{0}%
^{t}\int_{\mathbb{T}^{d}}(\operatorname{div}u)^{3}\nonumber\\
&  -\int_{0}^{t}\int_{\mathbb{T}^{d}}\varepsilon_{ijkl}\partial_{\ell
}u_{\delta}^{k}\omega_{\delta}\ast(\partial_{j}u^{q}\partial_{q}u^{i}%
)-\int_{0}^{t}\int_{\mathbb{T}^{d}}\varepsilon_{ijkl}\partial_{\ell}u_{\delta
}^{k}\left[  u^{q},\omega_{\delta}\right]  \partial_{qj}^{2}u^{i}\leq
C\int_{\mathbb{T}^{d}}\left\Vert \nabla u\right\Vert _{L^{3}(\mathbb{T}^{d}%
)}^{3}. \label{A_1_3}%
\end{align}

Using Poincar\'{e}'s inequality, we obtain that%
\begin{align}
&  \frac{1}{2}\int_{\mathbb{T}^{d}}\left\{  \partial_{t}\varepsilon
_{ijkl}+\partial_{q}(\varepsilon_{ijkl}u^{q})\right\}  \partial_{j}u_{\delta
}^{i}\partial_{k}u_{\delta}^{\ell}\nonumber\\
&  \hskip1cm\leq CE_{0}+C\int_{0}^{t}\int\left\Vert \nabla u\right\Vert
_{L^{3}(\mathbb{T}^{d})}^{3}+\int_{0}^{t}\left\vert \int_{\mathbb{T}^{d}}%
u^{q}\left(  \tau\right)  \right\vert \left\vert \int_{\mathbb{T}^{d}}%
\partial_{q}\varepsilon_{ijkl}\partial_{j}u_{\delta}^{i}\partial_{k}u_{\delta
}^{\ell}\right\vert
\end{align}
In order to treat the mean of $u^{q}$ we follow the ideea of P.--L. Lions's
\cite{lions1996mathematical2}. Let us recall, after verifying that $\rho^{2}$
is controlled by $B\left(  t\right)  $, that
\[
\left\Vert \rho\left(  t\right)  \left(  u^{q}\left(  t\right)  -\int
_{\mathbb{T}^{d}}u^{q}\left(  t\right)  \right)  \right\Vert _{L^{1}%
(\mathbb{T}^{d})}\leq C\left\Vert \rho\left(  t\right)  \right\Vert
_{L^{2}(\mathbb{T}^{d})}\left\Vert \nabla u\left(  t\right)  \right\Vert
_{L^{2}(\mathbb{T}^{d})}.
\]
Thus, we have that%
\[
\left\vert \int_{\mathbb{T}^{d}}\left(  \rho\left(  t,x\right)  \int
_{\mathbb{T}^{d}}u^{q}\left(  t,y\right)  dy-\rho\left(  t,x\right)
u^{q}\left(  t,x\right)  \right)  dx\right\vert \leq\left\Vert \rho\left(
t\right)  \left(  u^{q}\left(  t\right)  -\int_{\mathbb{T}^{d}}u^{q}\left(
t\right)  \right)  \right\Vert _{L^{1}(\mathbb{T}^{d})}%
\]
from which it follows that%
\begin{equation}
M\left\vert \int_{\mathbb{T}^{d}}u^{q}\left(  t,y\right)  dy\right\vert \leq
M^{\frac{1}{2}}E_{0}^{\frac{1}{2}}+C\left\Vert \rho\left(  t\right)
\right\Vert _{L^{2}(\mathbb{T}^{d})}\left\Vert \nabla u\left(  t\right)
\right\Vert _{L^{2}(\mathbb{T}^{d})}. \label{moyenne_u}%
\end{equation}
Consequently
\begin{align*}
&  \int_{0}^{t}\left\vert \int_{\mathbb{T}^{d}}u^{q}\left(  \tau\right)
\right\vert \left\vert \int_{\mathbb{T}^{d}}\partial_{q}\varepsilon
_{ijkl}\partial_{j}u_{\delta}^{i}\partial_{k}u_{\delta}^{\ell}\right\vert \\
&  \leq\left\Vert \partial_{q}\varepsilon_{ijkl}\right\Vert _{L^{\infty
}((0,t)\times\mathbb{T}^{d})}\left(  M^{-\frac{1}{2}}E_{0}^{\frac{3}{2}}%
+C\int_{0}^{t}\left\Vert \nabla u\left(  s\right)  \right\Vert _{L^{2}%
(\mathbb{T}^{d})}^{3}ds\right) \\
&  \leq\left\Vert \partial_{q}\varepsilon_{ijkl}\right\Vert _{L^{\infty
}((0,t)\times\mathbb{T}^{d})}\left(  M^{-\frac{1}{2}}E_{0}^{\frac{3}{2}}%
+C\int_{0}^{t}\left\Vert \nabla u\left(  s\right)  \right\Vert _{L^{3}%
(\mathbb{T}^{d})}^{3}ds\right)  .
\end{align*}
Taking $\eta$ sufficiently small and summing up $\left(
\text{\ref{d=2H1A_1_1}}\right)  $, $\left(  \text{\ref{d=2H1A_1_2}}\right)  $
we obtain
\begin{equation}
A_{1}\left(  t\right)  \leq C\left(  E_{0}+\int_{\mathbb{T}^{d}}H_{2}\left(
\rho_{0}/M\right)  +\left\Vert \nabla u_{0}\right\Vert _{L^{2}({\mathbb{T}%
}^{d})}^{2}\right)  +C_{\mu,\lambda,\gamma,M}\left(  E_{0}+B\left(  t\right)
+\sqrt{A_{1}\left(  t\right)  }\right)  A_{1}\left(  t\right)  .
\label{Bound_A1}%
\end{equation}

\paragraph{Terms appearing in the RHS of $\left(  \text{\ref{ref_def_A2}%
}\right)  $.}

We recall the definition of $A_{2}\left(  t\right)  $ which is%
\begin{align*}
&  \sigma\left(  t\right)  \int_{\mathbb{T}^{d}}\frac{\rho\left(  t\right)
\left\vert \dot{u}\left(  t\right)  \right\vert ^{2}}{2}+\mu\int
_{\mathbb{T}^{d}}\sigma\left\vert \partial_{k}\dot{u}^{i}\right\vert
^{2}+\left(  \mu+\lambda\right)  \int_{\mathbb{T}^{d}}\sigma\left\vert
\operatorname{div}\dot{u}\right\vert ^{2}+\int_{\mathbb{T}^{d}}\sigma
\varepsilon_{ijkl}\partial_{\ell}\dot{u}_{\delta}^{k}\partial_{j}\dot
{u}_{\delta}^{i}\\
&  =\int_{0}^{1}\int_{\mathbb{T}^{d}}\sigma\frac{\rho\left\vert \dot
{u}\right\vert ^{2}}{2}\\
&  +\mu\int_{0}^{t}\int_{\mathbb{T}^{d}}\sigma\partial_{k}u^{q}\partial
_{q}u^{i}\partial_{k}\dot{u}^{i}+\mu\int_{0}^{t}\int_{\mathbb{T}^{d}}%
\sigma\partial_{k}u^{q}\partial_{k}u^{i}\partial_{q}\dot{u}^{i}-\mu\int
_{0}^{t}\int_{\mathbb{T}^{d}}\sigma\operatorname{div}u\partial_{k}%
u^{i}\partial_{k}\dot{u}^{i}\\
&  +\left(  \mu+\lambda\right)  \int_{0}^{t}\int_{\mathbb{T}^{d}}%
\sigma\partial_{\ell}u^{q}\partial_{q}u^{\ell}\operatorname{div}\dot
{u}+\left(  \mu+\lambda\right)  \int_{0}^{t}\int_{\mathbb{T}^{d}}%
\sigma\partial_{i}u^{q}\partial_{q}\dot{u}^{i}\operatorname{div}u-\left(
\mu+\lambda\right)  \int_{0}^{t}\int_{\mathbb{T}^{d}}\sigma\left\vert
\operatorname{div}u\right\vert ^{2}\operatorname{div}\dot{u}\\
&  -\int_{0}^{t}\int_{\mathbb{T}^{d}}\sigma\left(  \partial_{t}\varepsilon
_{ijkl}+\partial_{q}\left(  u^{q}\varepsilon_{ijkl}\right)  \right)
\partial_{\ell}u_{\delta}^{k}\partial_{j}\dot{u}_{\delta}^{i}-\int_{0}^{t}%
\int_{\mathbb{T}^{d}}\sigma\varepsilon_{ijkl}(\omega_{\delta}\ast
(\partial_{\ell}u^{q}\partial_{q}u^{k}))\partial_{j}\dot{u}_{\delta}^{i}%
-\int_{0}^{t}\int_{\mathbb{T}^{d}}\sigma\varepsilon_{ijkl}\partial_{\ell
}u_{\delta}^{k}\omega_{\delta}\ast(\partial_{j}u^{q}\partial_{q}\dot{u}^{i})\\
&  +\int_{0}^{t}\int_{\mathbb{T}^{d}}\sigma\varepsilon_{ijkl}(\left[
u^{q},\omega_{\delta}\ast\right]  \partial_{\ell q}^{2}u^{k})\partial_{j}%
\dot{u}_{\delta}^{i}+\int_{0}^{t}\int_{\mathbb{T}^{d}}\sigma\varepsilon
_{ijkl}\partial_{\ell}u_{\delta}^{k}(\left[  u^{q},\omega_{\delta}\ast\right]
\partial_{jq}^{2}\dot{u}^{i})\\
&  -\int_{0}^{t}\int_{\mathbb{T}^{d}}\sigma\left\{  P\left(  \rho\right)
\partial_{j}u^{k}\partial_{k}\dot{u}^{j}+\left(  \rho P^{\prime}\left(
\rho\right)  -P\left(  \rho\right)  \right)  \operatorname{div}%
u\operatorname{div}\dot{u}\right\}  .
\end{align*}
We observe that%
\begin{align*}
&  \mu\int_{0}^{t}\int_{\mathbb{T}^{d}}\sigma\partial_{k}u^{q}\partial
_{q}u^{i}\partial_{k}\dot{u}^{i}+\mu\int_{0}^{t}\int_{\mathbb{T}^{d}}%
\sigma\partial_{k}u^{q}\partial_{k}u^{i}\partial_{q}\dot{u}^{i}-\mu\int
_{0}^{t}\int_{\mathbb{T}^{d}}\sigma\operatorname{div}u\partial_{k}%
u^{i}\partial_{k}\dot{u}^{i}\\
&  +\left(  \mu+\lambda\right)  \int_{0}^{t}\int_{\mathbb{T}^{d}}%
\sigma\partial_{\ell}u^{q}\partial_{q}u^{\ell}\operatorname{div}\dot
{u}+\left(  \mu+\lambda\right)  \int_{0}^{t}\int_{\mathbb{T}^{d}}%
\sigma\partial_{i}u^{q}\partial_{q}\dot{u}^{i}\operatorname{div}u-\left(
\mu+\lambda\right)  \int_{0}^{t}\int_{\mathbb{T}^{d}}\sigma\left\vert
\operatorname{div}u\right\vert ^{2}\operatorname{div}\dot{u}\\
&  +\int_{0}^{t}\int_{\mathbb{T}^{d}}\sigma\varepsilon_{ijkl}(\omega_{\delta
}\ast(\partial_{\ell}u^{q}\partial_{q}u^{k}))\partial_{j}\dot{u}_{\delta}%
^{i}-\int_{0}^{t}\int_{\mathbb{T}^{d}}\sigma\varepsilon_{ijkl}\partial_{\ell
}u_{\delta}^{k}\omega_{\delta}\ast(\partial_{j}u^{q}\partial_{q}\dot{u}^{i})\\
&  +\int_{0}^{t}\int_{\mathbb{T}^{d}}\sigma\varepsilon_{ijkl}(\left[
u^{q},\omega_{\delta}\ast\right]  \partial_{\ell q}^{2}u^{k})\partial_{j}%
\dot{u}_{\delta}^{i}+\int_{0}^{t}\int_{\mathbb{T}^{d}}\sigma\varepsilon
_{ijkl}\partial_{\ell}u_{\delta}^{k}(\left[  u^{q},\omega_{\delta}\ast\right]
\partial_{jq}^{2}\dot{u}^{i})\\
&  \int_{0}^{t}\int_{\mathbb{T}^{d}}\sigma\left\{  P\left(  \rho\right)
\partial_{j}u^{k}\partial_{k}\dot{u}^{j}+\left(  \rho P^{\prime}\left(
\rho\right)  -P\left(  \rho\right)  \right)  \operatorname{div}%
u\operatorname{div}\dot{u}\right\} \\
&  \leq P\left(  M\right)  E_{0}+\frac{C\left(  \eta\right)  }{2}\int_{0}%
^{t}\int_{\mathbb{T}^{d}}\sigma\left\vert P\left(  \rho\right)  -P\left(
M\right)  \right\vert ^{4}+C\left(  \eta\right)  \int_{0}^{t}\int
_{\mathbb{T}^{d}}\sigma\left\vert \nabla u\right\vert ^{4}+\eta\int_{0}%
^{t}\int_{\mathbb{T}^{d}}\sigma\left\vert \nabla\dot{u}\right\vert ^{2}.
\end{align*}
Moreover, using the Poincar\'{e} inequality along with $\left(
\text{\ref{moyenne_u}}\right)  $ we get that
\[
-\int_{0}^{t}\int_{\mathbb{T}^{d}}\sigma\left(  \partial_{t}\varepsilon
_{ijkl}+\partial_{q}\left(  u^{q}\varepsilon_{ijkl}\right)  \right)
\partial_{\ell}u_{\delta}^{k}\partial_{j}\dot{u}_{\delta}^{i}\leq C\left(
\eta\right)  E_{0}+\eta\int_{0}^{t}\int_{\mathbb{T}^{d}}\sigma\left\vert
\nabla\dot{u}\right\vert ^{2}.
\]
We thus see that for $\eta$ sufficiently small, we obtain that%
\begin{equation}
A_{2}\left(  t\right)  \leq C\left(  E_{0}+\int_{\mathbb{T}^{d}}H_{2}\left(
\rho_{0}/M\right)  +\left\Vert \nabla u_{0}\right\Vert _{L^{2}(\mathbb{T}%
^{d})}^{2}\right)  +C\left(  M+B\left(  t\right)  \right)  \left(
E_{0}+B\left(  t\right)  +\sqrt{A_{1}\left(  t\right)  }+A_{2}\left(
t\right)  \right)  A_{1}(t). \label{d=2H1_A2_final}%
\end{equation}

\paragraph{The bootstrap argument.}

Suppose that initially, $\left(  \rho,u\right)  $ are defined
on a time interval $[0,T)$. Assume that
\[
E\left(  \rho_{0}/M,u_{0}\right)  +\int_{\mathbb{T}^{d}}H_{2}\left(  \rho
_{0}/M\right)  +\left\Vert \nabla u_{0}\right\Vert _{L^{2}(\mathbb{T}^{d}%
)}^{2}\leq c_{0}%
\]
for some $c_{0}$ to be fixed later. We want to show that there exists a
constant $C=C\left(  \mu,\lambda,\gamma,M,E_{0}\right)  $ depending on
$\mu,\lambda,\gamma,M,E_{0}$ such that for $c_{0}$ sufficiently small%
\[
\forall t\in\lbrack0,T^{\ast}):\max_{[0,T)}\left\{  A_{1}\left(  t\right)
+A_{2}\left(  t\right)  +B\left(  t\right)  \right\}  \leq C\left(
\mu,\lambda,\gamma,M,E_{0}\right)  c_{0}.
\]
We will use a bootstrap argument. Recall that in the previous two sections, we
showed, see $\left(  \text{\ref{Bound_A1}}\right)  $ and $\left(
\text{\ref{d=2H1_A2_final}}\right)  $, that there exists a constant $C$ such
that%
\begin{align*}
A_{1}\left(  t\right)  +A_{2}\left(  t\right)   &  \leq C\left(  E_{0}%
+\int_{\mathbb{T}^{d}}H_{2}\left(  \rho_{0}/M\right)  +\left\Vert \nabla
u_{0}\right\Vert _{L^{2}}^{2}\right) \\
&  +C\left(  M+B\left(  t\right)  \right)  \left(  E_{0}+B\left(  t\right)
+\sqrt{A_{1}\left(  t\right)  }+A_{2}\left(  t\right)  \right)  A_{1}(t).
\end{align*}
We observe that
\[
B\left(  t\right)  \leq C\left(  E_{0}+\int_{\mathbb{T}^{d}}H_{2}\left(
\rho_{0}/M\right)  \right)  +C\left(  M+B\left(  t\right)  \right)  \left(
\sqrt{A_{1}\left(  t\right)  }+A_{2}\left(  t\right)  \right)  A_{1}\left(
t\right)  .
\]
We thus obtain that there exists a constant $\tilde{C}=\tilde{C}\left(
\mu,\lambda,\gamma,M,E_{0}\right)  $ such that
\begin{align}
A_{1}\left(  t\right)  +A_{2}\left(  t\right)  +B\left(  t\right)   &
\leq\tilde{C}\left(  E_{0}+\int_{\mathbb{T}^{2}}H_{2}\left(  \rho
_{0}/M\right)  +\left\Vert \nabla u_{0}\right\Vert _{L^{2}}^{2}\right)
\nonumber\\
&  +\tilde{C}\left(  M+B\left(  t\right)  \right)  \left(  E_{0}+B\left(
t\right)  +\sqrt{A_{1}\left(  t\right)  }+A_{2}\left(  t\right)  \right)
A_{1}\left(  t\right)  . \label{estimate_de_final}%
\end{align}

Let us introduce $T^{\ast}\in(0,T]$ such that%
\[
\max_{t\in\lbrack0,T^{\ast})}\left\{  \frac{1}{2}A_{1}\left(  t\right)
+A_{2}\left(  t\right)  +B\left(  t\right)  \right\}  \leq2\tilde{C}c_{0}.
\]
If $c_{0}$ is sufficiently small such that%
\[
\tilde{C}(M+2\tilde{C}c_{0})E_{0}\leq\frac{1}{2}%
\]
then the last two inequalities we observe that%
\[
\frac{1}{2}A_{1}\left(  t\right)  +A_{2}\left(  t\right)  +B\left(  t\right)
\leq\tilde{C}c_{0}+\tilde{C}\left(  M+2\tilde{C}c_{0}\right)  \left(
2\sqrt{\tilde{C}c_{0}}+4\tilde{C}c_{0}\right)  4\tilde{C}c_{0}%
\]
and thus for $c_{0}$ sufficiently small we have that
\[
\sqrt{c_{0}}\tilde{C}\left(  M+2\tilde{C}c_{0}\right)  \left(  \sqrt
{2\tilde{C}}+2\tilde{C}\sqrt{c_{0}}\right)  2\tilde{C}\leq\frac{1}{2}\tilde
{C},
\]
such that by a bootstrap argument we obtain that $T^{\ast}=T$.

\subsection{Construction of a sequence of approximate solution to $\left(
\text{\ref{system}}\right)  \label{Construction}$}

This entire part of the proof, from the local to the global existence of
strong solutions for the approximate system, can be achieved by repeating the
same arguments as in the work of R. Danchin and P.B. Mucha
\cite{DanchinMucha2019}, Section $3.$ It is for this reason that we only
briefly recall the strategy from \cite{DanchinMucha2019} and we comment on
what is different in our case. The form of the approximate systems $\left(
\text{\ref{ANS_delta}}\right)  $ is inspired from the one proposed by the
first author and P.--E. Jabin in \cite{bresch2018global} which we recopy here
for the reader's convenience:
\begin{equation}
\left\{
\begin{array}
[c]{l}%
\rho_{t}+\operatorname{div}\left(  \rho u\right)  =0,\\
\left(  \rho u\right)  _{t}+\operatorname{div}\left(  \rho u\otimes u\right)
+\nabla\rho^{\gamma}=\mu\Delta u+\left(  \mu+\lambda\right)  \nabla
\operatorname{div}u+\omega_{\delta}\ast\operatorname{div}(\mathcal{E(}%
\nabla\omega_{\delta}\ast u)).
\end{array}
\right.  \label{ANS_delta_1}%
\end{equation}
It differs from the classical (isotropic) Navier-Stokes system by a smooth
term. Therefore, we argue that the classical results regarding existence of
local solutions, see for instance Theorem $3.1.$ from \cite{DanchinMucha2019}
remain true in our case with a proof that is essentially the same. Thus, we
have that

\begin{theorem}
\label{ThmExistLocal2}Let $\rho_{0}\in W^{1,p}\left(  \mathbb{T}^{d}\right)  $
and $u_{0}\in W^{2-\frac{2}{p},2}\left(  \mathbb{T}^{d}\right)  $ for some
$p>d$ with $d\geq2.$ Assume that $\rho_{0}>0.$ Then there exists $T_{\ast}>0$
depending only on the norms of the data and on $\inf_{\mathbb{T}^{d}}\rho_{0}$
such that $\left(  \text{\ref{ANS_delta_1}}\right)  $ supplemented with data
$\rho_{0}$ and $u_{0}$ has a unique solution $\left(  \rho,u\right)  $ on the
time interval $\left[  0,T^{\ast}\right]  $, satisfying%
\[
u\in W^{1,p}\left(  \left(  0,T_{\ast}\right)  ;L^{p}\left(  \mathbb{T}%
^{d}\right)  \right)  \cap L^{p}\left(  \left(  0,T_{\ast}\right)
;W^{2,p}\left(  \mathbb{T}^{d}\right)  \right)  \text{ and }\rho\in
\mathcal{C}([0,T^{\ast}];W^{1,p}(\mathbb{T}^{d})).
\]

\end{theorem}

Consider
\[
\rho_{0}\in L^{2\gamma}\left(  \mathbb{T}^{d}\right)  \text{ and }u_{0}%
\in(H^{1}\left(  \mathbb{T}^{d}\right)  )^{d},
\]
with%
\[
\int_{\mathbb{T}^{d}}\rho_{0}\left(  x\right)  dx=M.
\]
Also, consider $\omega_{\delta}=\frac{1}{\delta^{d}}\omega\left(  \frac{\cdot
}{\delta}\right)  $ with $\omega$ a smooth, nonnegative, radial function
compactly supported in the unit ball centered at the origin and with integral
equal to $1$. To this end, for all $\delta\in(0,M)$ there exists $\xi_{\delta
}>\delta$ such that
\[
\tilde{\rho}_{0}^{\delta}\left(  x\right)  =\min\left\{  \xi_{\delta},\rho
_{0}\left(  x\right)  +\delta\right\}  \text{ and }\int_{\mathbb{T}^{d}}%
\tilde{\rho}_{0}^{\delta}\left(  x\right)  dx=M.
\]
We consider%
\[
\rho_{0}^{\delta}\left(  x\right)  =\omega_{\delta}\ast\tilde{\rho}%
_{0}^{\delta}\text{ and }u_{0}^{\delta}=\omega_{\delta}\ast u_{0}.
\]
Observe that for a subsequence that, which by slightly abusing the notation we
still denote by the index $\delta$, it holds that:%
\[
\int_{\mathbb{T}^{d}}\rho_{0}^{\delta}\left(  x\right)  dx=M\text{ and }%
\lim_{\delta\rightarrow0}\left\Vert \rho_{0}^{\delta}-\rho_{0}\right\Vert
_{L^{2\gamma}(\mathbb{T}^{d})}+\lim_{\delta\rightarrow0}\left\Vert
u_{0}^{\delta}-u_{0}\right\Vert _{H^{1}(\mathbb{T}^{d})}=0.
\]

Consider $\left(  \rho^{\delta},u^{\delta}\right)  _{\delta}$ the sequence of
solutions for the Cauchy problem associated to system $\left(
\text{\ref{ANS_delta}}\right)  $ with initial
\[
\left\{
\begin{array}
[c]{c}%
\rho_{|t=0}=\rho_{0}^{\delta},\\
u_{|t=0}=u_{0}^{\delta},
\end{array}
\right.
\]
the existence of which is granted by Theorem \ref{ThmExistLocal2}. A priori,
each of $\left(  \rho^{\delta},u^{\delta}\right)  $ is defined on its own
maximal time interval $[0,T_{\delta})$ with $T_{\delta}\in(0,\infty]$. On
these time intervals the solution has enough regularity such that the
computations performed above make sense and, as a consequence, we have that
$\left(  \rho^{\delta},u^{\delta}\right)  $ have
\[
E\left(  \rho^{\delta}/M,u^{\delta}\right)  ,A_{1}^{\delta}\left(  t\right)
,A_{2}^{\delta}\left(  t\right)  ,B^{\delta}\left(  t\right)
\]
bounded independently w.r.t. $\delta$ where $A_{1}^{\delta}\left(  t\right)
,A_{2}^{\delta}\left(  t\right)  ,B^{\delta}\left(  t\right)  $ are the
expressions defined in $\left(  \text{\ref{A1_Hoff}}\right)  $, $\left(
\text{\ref{A2_Hoff}}\right)  $ respectively in $\left(
\text{\ref{definition_B}}\right)  $ with $\left(  \rho^{\delta},u^{\delta
}\right)  $ instead of $\left(  \rho,u\right)  $. From here on, the argument
leading to the conclusion that $T^{\delta}=+\infty$ continues as in
\cite{DanchinMucha2019} due to the fact that the term $\omega_{\delta}%
\ast\operatorname{div}(\mathcal{E(}\nabla\omega_{\delta}\ast u))$ is regular.

In particular, taking also into account $\left(  \text{\ref{moyenne_u}%
}\right)  $ we get that for all $T>0$ we have that%
\begin{equation}
\left\Vert \rho^{\delta}\right\Vert _{L^{\infty}((0,T);L^{2\gamma}%
(\mathbb{T}^{d}))\cap L^{3\gamma}((0,T)\times\mathbb{T}^{d})}+\left\Vert
u^{\delta}\right\Vert _{L^{3}((0,T);W^{1,3}(\mathbb{T}^{d}))}+\left\Vert
\sqrt{\sigma}\nabla\dot{u}^{\delta}\right\Vert _{L^{2}((0,T)\times
\mathbb{T}^{d})}\leq C\left(  \mu,\lambda,\gamma,\mathcal{E}\right)  .
\label{Hoff_uniform_bounds}%
\end{equation}

\textit{Strong convergence of the sequence }$\left(  \mathit{u}^{\delta
}\right)  _{\delta>0}$\textit{.} We want to prove that the uniform bounds
verified by the solutions constructed above imply that up to a subsequence%
\[
\lim_{\delta\rightarrow0}u^{\delta}=u\text{ strongly in }L^{2}((\tfrac{1}%
{n},T)\times\mathbb{T}^{d}),
\]
for all $n\in\mathbb{N}^{\ast}$. This a consequence of the fact that the
second Hoff functional is uniformly bounded w.r.t. $\delta>0$ which implies
that for all $T>0$ and all $\delta>0$
\[
\int_{0}^{T}\int_{\mathbb{T}^{d}}\sigma\left(  t\right)  \left\Vert \nabla
\dot{u}^{\delta}\right\Vert _{L^{2}(\mathbb{T}^{d})}^{2}\leq c,
\]
for some $c$. This implies that for any $n\in\mathbb{N}^{\ast}$
\[
\int_{1/n}^{T}\int_{\mathbb{T}^{d}}\left\Vert \nabla\dot{u}^{\delta
}\right\Vert _{L^{2}(\mathbb{T}^{d})}^{2}\leq nc.
\]
We remark that since $d\in\left\{  2,3\right\}  $ we have that
\[
\left\Vert u^{\delta}\right\Vert _{L^{3}((0,T);L^{p}(\mathbb{T}^{d}%
))}+\left\Vert \nabla u^{\delta}\right\Vert _{L^{3}((0,T);L^{p}(\mathbb{T}%
^{d}))}%
\]
is uniformly bounded for any $p<\infty$. Taking into account that
$\partial_{t}u^{\delta}=\dot{u}^{\delta}-u^{\delta}\cdot\nabla u^{\delta},$
(the mean value of $\dot{u}^{\delta}$ is controlled exactly as in
\ref{moyenne_u}) we obtain that for all $\eta\in(0,1]$ we have
\begin{equation}
\partial_{t}u^{\delta}\text{ is uniformly bounded in }L^{\frac{3}{2}%
}((1/n,T);L^{3-\eta}(\mathbb{T}^{d})). \label{bound_on_time_derivative}%
\end{equation}
For any $n$, by the Aubin-Lions theorem $\left(  u^{\delta}\right)
_{\delta>0}$ converges strongly in $L^{2}\left(  \left(  \frac{1}{n},T\right)
\times\mathbb{T}^{d}\right)  $ while applying a Cantor's diagonal type process
provides us with a subsequence $\left(  u^{\delta}\right)  _{\delta>0}$
converging for any $n$ in $L^{2}\left(  \left(  \frac{1}{n},T\right)
\times\mathbb{T}^{d}\right)  $.

\begin{remark}
We use $\delta$ as upperscript to designate the sequence of approximate
solutions $\left(  \rho^{\delta},u^{\delta}\right)  _{\delta>0}$. This should
not to be confused with the lower-script notation used in the previous section
which denoted the regularisation with $\omega_{\delta}$. In order to avoid
possible confusion, all along this section we explicitely write $\omega
_{\delta}\ast(\cdot)$.
\end{remark}

\subsection{Stability of solutions for $\left(  \text{\ref{ANS_delta_1}%
}\right)  $}

\label{weakstab}

In this section we show that from the sequence $\left(  \rho^{\delta
},u^{\delta}\right)  _{\delta}$ constructed in the previous section and
verifying uniformly the Hoff-type estimates $\left(
\text{\ref{Hoff_uniform_bounds}}\right)  $, one can extract a subsequence
converging weakly towards a solution of the system $\left(
\text{\ref{system_pert}}\right)  $. We recall that this is not trivial given
the fact that the pressure is a nonlinear function of the density. From the
estimates we have gathered so far for $\left(  \rho^{\delta},u^{\delta
}\right)  _{\delta}$ we infer the existence of $\left(  \rho,u,\overline
{\rho^{\gamma}}\right)  $ such that modulo an extraction of a subsequence:%

\begin{equation}
\left\{
\begin{array}
[c]{l}%
\rho^{\delta}\rightharpoonup\rho\text{ weakly in }L^{3\gamma}\left(
(0,T)\times\mathbb{T}^{d}\right)  \text{ and weakly-}\ast\text{ in }L^{\infty
}(\left(  0,T\right)  ;L^{2\gamma}\left(  \mathbb{T}^{d}\right)  ),\\
\rho^{\delta}\rightarrow\rho\text{ strongly in }\mathcal{C}^{0}([0,T];L_{weak}%
^{2\gamma}),\\
(\rho^{\delta})^{\gamma}\rightharpoonup\overline{\rho^{\gamma}}\text{ weakly
in }L^{3}\left(  (0,T)\times\mathbb{T}^{d}\right)  \text{ and weakly-}%
\ast\text{ in }L^{\infty}(\left(  0,T\right)  ;L^{2}\left(  \mathbb{T}%
^{d}\right)  ),\\
(\rho^{\delta})^{\gamma}\rightharpoonup\overline{\rho^{\gamma}}\text{ strongly
in }\mathcal{C}^{0}([0,T];L_{weak}^{2}),\\
u^{\delta}\rightharpoonup u\text{ weakly in }\left(  L^{3}(0,T;W^{1,3}%
(\mathbb{T}^{d}))\right)  ^{d},\\
u^{\delta}\rightarrow u\text{ strongly in }\left(  L^{2}((\frac{1}{n}%
,T)\times{\mathbb{T}}^{d}))\right)  ^{d},
\end{array}
\right.  \label{uniform}%
\end{equation}
and, in view of $\left(  \text{\ref{bound_on_time_derivative}}\right)  $ for
any $n\in\mathbb{N}$ and for any $p\in\lbrack2,3):$
\begin{equation}
\partial_{t}u\in L^{\frac{3}{2}}((1/n,T);L^{p}(\mathbb{T}^{d})).
\label{estimate_time_derivative}%
\end{equation}
See the Appendix for a definition for $\mathcal{C}^{0}([0,T];L_{weak}^{2})$.
In particular, since we work in dimension $d\in\left\{  2,3\right\}  $, since
$u\in\left(  L^{3}((0,T);W^{1,3}(\mathbb{T}^{d}))\right)  ^{d}$ and owing to
$\left(  \text{\ref{estimate_time_derivative}}\right)  $ we obtain that for
all $n\in\mathbb{N}$ we have:%
\begin{equation}
\left\{
\begin{array}
[c]{l}%
\rho\in L_{t}^{\infty}L_{x}^{2},\text{ }\rho u\in L^{3}((0,T);L^{\frac{12}{7}%
}(\mathbb{T}^{d})),\text{ }\rho u\otimes u\in L^{\frac{3}{2}}((0,T)\times
\mathbb{T}^{d}),\\
\partial_{t}u\in L^{\frac{3}{2}}((1/n,T);L^{\frac{12}{5}}(\mathbb{T}%
^{d})),\partial_{t}\left\vert u\right\vert ^{2}\in L^{1}((1/n,T);L^{2}%
(\mathbb{T}^{d}))\text{, }\nabla\left\vert u\right\vert ^{2}\in L^{\frac{3}%
{2}}((0,T);L^{\frac{12}{5}}(\mathbb{T}^{d})).\text{ }%
\end{array}
\right.  \label{information_for_the_limit}%
\end{equation}
All the relations of $\left(  \text{\ref{uniform}}\right)  $ are applications
of classical results from functional analysis. The second and fourth relations
are obtained using a weak variant of the Arzel\`{a}-Ascoli see Theorem
\ref{Arzela-Ascoli-weak} from the appendix, a proof of which can be found, for
instance in Vrabie \cite{Vrabie2003Semigroups}. We will give details for the
proof of the fourth relation of $\left(  \text{\ref{uniform}}\right)  $ since
the second one follows from similar arguments. Since the sequence $\left(
(\rho^{\delta})^{\gamma}\right)  _{\delta>0}$ is bounded in $L^{\infty
}((0,T);L^{2}\left(  \mathbb{T}^{d}\right)  )$ the second condition from
Theorem \ref{Arzela-Ascoli-weak} obviously holds true. It remains to prove
that $\left(  (\rho^{\delta})^{\gamma}\right)  _{\delta>0}$ is weakly
equicontinious on $\left[  0,T\right]  $. Fix $\varepsilon>0$ and a $w\in
L^{2}\left(  \mathbb{T}^{d}\right)  $ and consider $\tilde{w}\in
\mathcal{C}_{per}^{\infty}\left(  \mathbb{R}^{d}\right)  $. We see that there
exists a numerical constant $C$ independent of $\delta$ such that%
\begin{align*}
&  \left\langle (\rho^{\delta})^{\gamma}\left(  t\right)  -(\rho^{\delta
})^{\gamma}\left(  s\right)  ,w\right\rangle \\
&  =\left\langle (\rho^{\delta})^{\gamma}\left(  t\right)  -(\rho^{\delta
})^{\gamma}\left(  s\right)  ,w-\tilde{w}\right\rangle +\left\langle
(\rho^{\delta})^{\gamma}\left(  t\right)  -(\rho^{\delta})^{\gamma}\left(
s\right)  ,\tilde{w}\right\rangle \\
&  \leq2\left\Vert (\rho^{\delta})^{\gamma}\right\Vert _{L^{\infty
}((0,T);L^{2}(\mathbb{T}^{d}))}\left\Vert w-\tilde{w}\right\Vert
_{L^{2}(\mathbb{T}^{d})}+C\int_{s}^{t}\left\Vert (\rho^{\delta})^{\gamma
}\right\Vert _{L^{3}(\mathbb{T}^{d})}\left\Vert u^{\delta}\right\Vert
_{W^{1,3}(\mathbb{T}^{d})}\left\Vert \nabla\tilde{w}\right\Vert _{L^{3}%
(\mathbb{T}^{d})}.\\
&  \leq2\left\Vert (\rho^{\delta})^{\gamma}\right\Vert _{L^{\infty
}((0,T);L^{2}(\mathbb{T}^{d}))}\left\Vert w-\tilde{w}\right\Vert
_{L^{2}(\mathbb{T}^{d})}+C\left\vert t-s\right\vert ^{\frac{1}{3}}\left\Vert
(\rho^{\delta})^{\gamma}\right\Vert _{L^{3}((0,T)\times\mathbb{T}^{d}%
)}\left\Vert u^{\delta}\right\Vert _{L^{3}((0,T);W^{1,3}(\mathbb{T}^{d}%
))}\left\Vert \nabla\tilde{w}\right\Vert _{L^{3}(\mathbb{T}^{d})}.
\end{align*}
Using $\left(  \text{\ref{Hoff_uniform_bounds}}\right)  $, the first term can
be made arbitrarily small because we can approximate $L^{2}\left(
\mathbb{T}^{d}\right)  $ functions with smooth periodic functions while the
second term can be made arbitrarily small provided $\left\vert t-s\right\vert
$ is chosen appropriately

Using the strong convergences of $\rho^{\delta}\rightarrow\rho$ in
$\mathcal{C}^{0}([0,T];L_{weak}^{2\gamma})$, of $(\rho^{\delta})^{\gamma
}\rightharpoonup\overline{\rho^{\gamma}}$ in $\mathcal{C}^{0}([0,T];L_{weak}%
^{2})$ and the fact that $\rho_{0}^{\delta}\rightarrow\rho_{0}$ strongly in
$L^{2\gamma}\left(  \mathbb{T}^{d}\right)  $ we recover that%
\begin{equation}
\left\{
\begin{array}
[c]{l}%
\lim\limits_{t\rightarrow0}\int_{\mathbb{T}^{d}}\rho\left(  t,x\right)
\psi\left(  x\right)  dx=\int_{\mathbb{T}^{d}}\rho_{0}\left(  x\right)
\psi\left(  x\right)  dx\text{ and that}\\
\lim\limits_{t\rightarrow0}\int_{\mathbb{T}^{d}}\overline{\rho^{\gamma}%
}\left(  t,x\right)  \psi\left(  x\right)  dx=\int_{\mathbb{T}^{d}}\rho
_{0}^{\gamma}\left(  x\right)  \psi\left(  x\right)  dx,\text{ for all }%
\psi\in C_{per}^{\infty}\text{.}%
\end{array}
\right.  \label{weak_compact_0}%
\end{equation}

It is by now well-understood that the assumptions $\left(  \text{\ref{uniform}%
}\right)  $ are sufficient in order to conclude that%
\begin{equation}
\left\{
\begin{array}
[c]{l}%
\partial_{t}\rho+\operatorname{div}\left(  \rho u\right)  =0,\\
\partial_{t}(\rho u)+\operatorname{div}(\rho u\otimes u)-\mu\Delta u-\left(
\mu+\lambda\right)  \nabla\operatorname{div}u+\nabla\overline{\rho^{\gamma}%
}=\operatorname{div}\left(  \mathcal{E(}\nabla u)\right)  ,
\end{array}
\right.  \label{limit_not_yet_ident}%
\end{equation}
Moreover,
\begin{equation}
\rho\in C\left(  [0,T);L^{p}\left(  \mathbb{T}^{3}\right)  \right)  \text{ for
all }1\leq p<2\gamma, \label{continuity_rho}%
\end{equation}
see for instance Lemma $6.15$, page $312$ of \cite{NoSt}. In the same manner,
using the fact that
\[
\partial_{t}\overline{\rho^{\gamma}}+\operatorname{div}\left(  \overline
{\rho^{\gamma}}u\right)  +\left(  \gamma-1\right)  \overline{\rho^{\gamma
}\operatorname{div}u}=0,
\]
which comes from passing to the limit in the equations verified by
$(\rho^{\delta})^{\gamma}$ we obtain that:
\begin{equation}
\overline{\rho^{\gamma}}\in C\left(  [0,T);L^{p}\left(  \mathbb{T}^{3}\right)
\right)  \text{ for }1\leq p<2. \label{continuity_rho_gamma}%
\end{equation}
Of course, in order to finish the proof we must show that the function
$\overline{\rho^{\gamma}}$ coincides with the function $\rho^{\gamma}$. To
this end we will essentially mimic the proof from \cite{BrBu1} which consists
of taking the difference between the limit of the energy equations with the
energy equation of the limiting system and "multiplying" it with an
appropriate quantity that yields a "conservative" identity.

Using once more $\left(  \text{\ref{uniform}}\right)  $ we obtain the
existence of positive functions $\overline{\nabla u:\nabla u},\overline
{(\operatorname{div}u)^{2}},$ $\overline{\mathcal{E(}\nabla u):\nabla u}$ $\in
L^{\frac{3}{2}}\left(  (0,T)\times\mathbb{T}^{d}\right)  $ such that up to a
subsequence we have%
\begin{equation}
\left\{
\begin{array}
[c]{l}%
\mathcal{\nabla}u^{\delta}:\nabla u^{\delta}\rightharpoonup\overline{\nabla
u:\nabla u}\text{ in }L^{\frac{3}{2}}\left(  (0,T)\times\mathbb{T}^{d}\right)
\text{ and }\nabla u:\nabla u\leq\overline{\nabla u:\nabla u},\\
(\operatorname{div}u^{\delta})^{2}\rightharpoonup\overline{(\operatorname{div}%
u)^{2}}\text{ in }L^{\frac{3}{2}}\left(  (0,T)\times\mathbb{T}^{d}\right)
\text{ and }(\operatorname{div}u)^{2}\leq\overline{(\operatorname{div}u)^{2}%
},\\
\mathcal{E(\nabla(}\omega_{\delta}\ast u^{\delta})):\nabla(\omega_{\delta}\ast
u^{\delta})\rightharpoonup\overline{\mathcal{E(}\nabla u):\nabla u}\text{ in
}L^{\frac{3}{2}}\left(  (0,T)\times\mathbb{T}^{d}\right)  \text{ }\\
\text{
\ \ \ \ \ \ \ \ \ \ \ \ \ \ \ \ \ \ \ \ \ \ \ \ \ \ \ \ \ \ \ \ \ \ \ \ \ \ \ \ \ \ \ and
}\mathcal{E(}\nabla u):\nabla u\leq\overline{\mathcal{E(}\nabla u):\nabla
u}\text{ .}%
\end{array}
\right.  \label{weak_conv2}%
\end{equation}
It is in the proof of the last property, that the we need to regularize a
positive definite operator and the assumption made in Remark \ref{important}.
See the remark to see that simple change of shear viscosity that may be made
to satisfy such property starting with a viscosity tensor satisfying
Hypothesis \eqref{Hyp1}--\eqref{Hyp4}.

\smallskip

\noindent\textit{Lower semi-continuity.} Indeed, for any $\phi\in
L^{3}{((0,T)\times\mathbb{T}^{d})}$ with $\phi\geq0$, denoting $u^{\delta
}=\left(  u^{\delta,1},\cdots,u^{\delta,d}\right)  $, we have that%
\begin{align*}
0  &  \leq\int_{0}^{T}\int_{\mathbb{T}^{d}}\mathcal{E(}\nabla(\omega_{\delta
}\ast u^{\delta})-\nabla u):(\nabla(\omega_{\delta}\ast u^{\delta})-\nabla
u)\phi\\
&  =\int_{0}^{T}\int_{\mathbb{T}^{d}}\varepsilon_{ijkl}(\partial_{\ell}%
\omega_{\delta}\ast u^{\delta,k}-\partial_{\ell}u^{k})(\partial_{j}%
\omega_{\delta}\ast u^{\delta,i}-\partial_{j}u^{i})\phi\\
&  =\int_{0}^{T}\int_{\mathbb{T}^{d}}\varepsilon_{ijkl}\partial_{\ell}%
\omega_{\delta}\ast u^{\delta,k}\partial_{j}\omega_{\delta}\ast u^{\delta
,i}\phi-\int_{0}^{T}\int_{\mathbb{T}^{d}}\varepsilon_{ijkl}\partial_{\ell
}\omega_{\delta}\ast u^{\delta,k}\partial_{j}u^{i}\phi-\int_{0}^{T}%
\int_{\mathbb{T}^{d}}\varepsilon_{ijkl}\partial_{\ell}u^{k}\partial_{j}%
\omega_{\delta}\ast u^{\delta,i}\phi\\
&  +\int_{0}^{T}\int_{\mathbb{T}^{d}}\varepsilon_{ijkl}\partial_{\ell}%
u^{k}\partial_{j}u^{i}\phi.
\end{align*}
We obviously have%
\[
\lim_{\delta\rightarrow0}\int_{0}^{T}\int_{\mathbb{T}^{d}}\omega_{\delta}%
\ast(\varepsilon_{ijkl}\partial_{j}u^{i}\phi)\partial_{\ell}u^{\delta,k}%
=\int_{0}^{T}\int_{\mathbb{T}^{d}}\varepsilon_{ijkl}\partial_{\ell}%
u^{k}\partial_{j}u^{i}\phi
\]
and the same for the other similar term. Thus we obtain that%
\begin{align*}
0  &  \leq\lim_{\delta\rightarrow0}\int_{0}^{T}\int_{\mathbb{T}^{d}%
}\varepsilon_{ijkl}\partial_{\ell}\omega_{\delta}\ast u^{\delta,k}\partial
_{j}\omega_{\delta}\ast u^{\delta,i}\phi-\int_{0}^{T}\int_{\mathbb{T}^{d}%
}\varepsilon_{ijkl}\partial_{\ell}u^{k}\partial_{j}u^{i}\phi\\
&  =\int_{0}^{T}\int_{\mathbb{T}^{d}}\left(  \overline{\mathcal{E(}\nabla
u):\nabla u}-\mathcal{E(}\nabla u):\nabla u\right)  \phi.
\end{align*}

\noindent\textit{Energy identities and conclusion.} On the one hand, for any
$\delta>0$ , the regularity of $\left(  \rho^{\delta},u^{\delta}\right)  $,
see Theorem \ref{ThmExistLocal2} and the remark that follows, we may write the
following energy equation:%
\begin{align}
\frac{1}{2}\frac{\partial}{\partial t}\left\{  \rho^{\delta}\left\vert
u^{\delta}\right\vert ^{2}+\frac{(\rho^{\delta})^{\delta}}{\gamma-1}\right\}
+\operatorname{div}\left(  \left(  \rho^{\delta}\left\vert u^{\delta
}\right\vert ^{2}+\frac{(\rho^{\delta})^{\gamma}}{\gamma-1}\right)  u^{\delta
}\right)  +\mu\nabla u^{\delta}  &  :\nabla u^{\delta}+\left(  \mu
+\lambda\right)  (\operatorname{div}u^{\delta})^{2}\nonumber\\
-\mu\Delta\frac{\left\vert u^{\delta}\right\vert ^{2}}{2}-\left(  \mu
+\lambda\right)  \operatorname{div}\left(  u^{\delta}\operatorname{div}%
u^{\delta}\right)  -\omega_{\delta}\ast\operatorname{div}\left(
\mathcal{E}\nabla\omega_{\delta}\ast u^{\delta}\right)  u^{\delta}  &  =0,
\label{equ_delta}%
\end{align}
Let us observe that for all $\phi\in\mathcal{C}_{c}\left(  (0,T);\mathcal{C}%
_{per}^{\infty}\left(  \mathbb{R}^{d}\right)  \right)  $ we have that%
\begin{align*}
-\int_{\mathbb{T}^{d}}\omega_{\delta}\ast\operatorname{div}\left(
\mathcal{E}\nabla\omega_{\delta}\ast u^{\delta}\right)  u^{\delta}\phi &
=\int_{\mathbb{T}^{d}}\mathcal{E}\nabla(\omega_{\delta}\ast u^{\delta}%
):\omega_{\delta}\ast\nabla(u^{\delta}\phi)\\
&  =\int_{\mathbb{T}^{d}}\mathcal{E}\nabla(\omega_{\delta}\ast u^{\delta
}):\omega_{\delta}\ast(\nabla u^{\delta}\phi)+\int_{\mathbb{T}^{d}}%
\mathcal{E}\nabla(\omega_{\delta}\ast u^{\delta}):\omega_{\delta}%
\ast(u^{\delta}\otimes\nabla\phi).
\end{align*}
Owing to the fact that there exits some $n$ such that
\[
\operatorname*{Supp}\phi\left(  \cdot,\cdot\right)  \subset(1/n,T)\times
\mathbb{R}^{d},
\]
that $u^{\delta}\rightarrow u$ strongly in $L^{2}((1/n,T)\times{\mathbb{T}%
}^{d}))^{d}$ and that $\nabla(\omega_{\delta}\ast u^{\delta})\rightarrow\nabla
u$ weakly in $L^{3}((0,T)\times{\mathbb{T}}^{d}))^{d\times d}$ we obtain that%
\[
\lim_{\delta\rightarrow0}\int_{0}^{T}\int_{\mathbb{T}^{d}}\mathcal{E(}%
\nabla(\omega_{\delta}\ast u^{\delta})):\omega_{\delta}\ast(u^{\delta}%
\otimes\nabla\phi)=\int_{0}^{T}\int_{\mathbb{T}^{d}}\mathcal{E(}\nabla
u):(u\otimes\nabla\phi).
\]
Next, we observe that%
\begin{align*}
\int_{0}^{T}\int_{\mathbb{T}^{d}}\mathcal{E}\nabla(\omega_{\delta}\ast
u^{\delta})  &  :\omega_{\delta}\ast(\nabla u^{\delta}\phi)\\
&  =\int_{0}^{T}\int_{\mathbb{T}^{d}}\mathcal{E}\nabla(\omega_{\delta}\ast
u^{\delta}):\nabla(\omega_{\delta}\ast u^{\delta})\phi+\int_{0}^{T}%
\int_{\mathbb{T}^{d}}\mathcal{E}\nabla(\omega_{\delta}\ast u^{\delta}):\left[
\omega_{\delta}\ast,\phi\right]  \nabla u^{\delta}.
\end{align*}
Now, for any $j,q\in\overline{1,d}$ one has%
\begin{align*}
\left[  \omega_{\delta}\ast,\phi\right]  \partial_{j}u^{\delta,q}\left(
t,x\right)   &  =\left(  \omega_{\delta}\ast\left(  \phi\partial_{j}%
u^{\delta,q}\right)  -\phi\omega_{\delta}\ast\partial_{j}u^{\delta,q}\right)
\left(  t,x\right) \\
&  =\int_{\mathbb{T}^{d}}\left(  \phi\left(  t,x-y\right)  -\phi\left(
t,x\right)  \right)  \partial_{j}u^{\delta,q}\left(  t,x-y\right)
\omega_{\delta}\left(  y\right)  dy\\
&  =\int_{\mathbb{T}^{d}}\left(  \phi\left(  t,x-\delta z\right)  -\phi\left(
t,x\right)  \right)  \partial_{j}u^{\delta,q}\left(  t,x-\delta z\right)
\omega\left(  z\right)  dz.
\end{align*}
Thus%
\[
\left\vert \int_{0}^{T}\int_{\mathbb{T}^{d}}\mathcal{E}\nabla(\omega_{\delta
}\ast u^{\delta}):\left[  \omega_{\delta}\ast,\phi\right]  \nabla u^{\delta
}\right\vert \leq\delta\max_{i,j,k,l}\left\Vert \varepsilon_{ijkl}\right\Vert
_{L^{\infty}(\mathbb{T}^{d})}\left\Vert \nabla u^{\delta}\right\Vert
_{L^{2}(\mathbb{T}^{d})}^{2}\left\Vert \nabla\phi\right\Vert _{L^{\infty}%
}\underset{\delta\rightarrow0}{\rightarrow}0.
\]
Moreover, using the information of relation $\left(  \text{\ref{weak_conv2}%
}\right)  $ we may pass to the limit in $\left(  \text{\ref{equ_delta}%
}\right)  $ such as to obtain
\begin{gather}
\frac{1}{2}\frac{\partial}{\partial t}\left\{  \rho\left\vert u\right\vert
^{2}+\frac{\overline{\rho^{\gamma}}}{\gamma-1}\right\}  +\operatorname{div}%
\left(  \left(  \rho\left\vert u\right\vert ^{2}+\frac{\overline{\rho^{\gamma
}}}{\gamma-1}\right)  u\right)  +\mu\overline{\nabla u:\nabla u}+\left(
\mu+\lambda\right)  \overline{(\operatorname{div}u)^{2}}+\overline
{\mathcal{E}\left(  \nabla u\right)  :\nabla u}\nonumber\\
-\mu\Delta\frac{\left\vert u\right\vert ^{2}}{2}-\left(  \mu+\lambda\right)
\operatorname{div}\left(  u\operatorname{div}u\right)  -\operatorname{div}%
\left(  u\mathcal{E}(\nabla u)\right)  =0. \label{passage_limite}%
\end{gather}

On the other hand, let us observe that system $\left(
\text{\ref{limit_not_yet_ident}}\right)  $ can be put under the form
\begin{equation}
\left\{
\begin{array}
[c]{l}%
\partial_{t}\rho+\operatorname{div}\left(  \rho u\right)  =0,\\
\partial_{t}(\rho u)+\operatorname{div}(\rho u\otimes u)-\operatorname{div}%
\left(  \mathcal{E}\nabla u\right)  +\nabla\rho^{\gamma}=\nabla(\rho^{\gamma
}-\overline{\rho^{\gamma}}).
\end{array}
\right.  \label{NS_neidentificat}%
\end{equation}
We observe that in view of the estimates $\left(
\text{\ref{information_for_the_limit}}\right)  $ and by a density argument,
the weak formulations for the transport equation holds true for any test
function $\psi$ such that%
\[
\partial_{t}\psi\in L^{1}((0,T);L^{2}(\mathbb{T}^{d}))\text{, }\nabla\psi\in
L^{\frac{3}{2}}((0,T);L^{\frac{12}{5}}(\mathbb{T}^{d}))
\]
while the momentum equation holds true for vector fields with coefficients
$\psi$ such that:%
\[
\partial_{t}\psi\in L^{\frac{3}{2}}((0,T);L^{\frac{12}{5}}(\mathbb{T}%
^{d}))\text{ and }\nabla\psi\in L^{3}((0,T)\times\mathbb{T}^{d}).
\]
The same estimates $\left(  \text{\ref{information_for_the_limit}}\right)  $
show that for any $\phi\in\mathcal{C}_{c}^{\infty}\left(  (0,T);\mathcal{C}%
_{per}^{\infty}\left(  \mathbb{R}^{d}\right)  \right)  $ one can use
$\left\vert u\right\vert ^{2}\phi/2$ as a test function in the first equation
of $\left(  \text{\ref{NS_neidentificat}}\right)  $ while one can also use
$u\phi$ as a test function in the second equation of $\left(
\text{\ref{NS_neidentificat}}\right)  $. By doing so, summing up the two
relations that result, taking into account the chain-rule respectively the
derivation rule of products of functions in Sobolev spaces, Proposition
\ref{Prop_ren1} and finally the fact that $\phi$ is chosen arbitrarily, we
obtain that:
\begin{gather}
\frac{1}{2}\frac{\partial}{\partial t}\left\{  \rho\left\vert u\right\vert
^{2}+\frac{\rho^{\gamma}}{\gamma-1}\right\}  +\operatorname{div}\left(
\left(  \rho\left\vert u\right\vert ^{2}+\frac{\rho^{\gamma}}{\gamma
-1}\right)  u\right)  +\mu\nabla u:\nabla u+\left(  \mu+\lambda\right)
(\operatorname{div}u)^{2}+\mathcal{E}\left(  \nabla u\right)  :\nabla
u\nonumber\\
-\mu\Delta\frac{\left\vert u\right\vert ^{2}}{2}-\left(  \mu+\lambda\right)
\operatorname{div}\left(  u\operatorname{div}u\right)  -\operatorname{div}%
\left(  u\mathcal{E(\nabla}u)\right)  =\operatorname{div}\left(  u\left(
\rho^{\gamma}-\overline{\rho^{\gamma}}\right)  \right)  -\left(  \rho^{\gamma
}-\overline{\rho^{\gamma}}\right)  \operatorname{div}u,
\label{renorm_sans_eps}%
\end{gather}
which holds true in $\mathcal{D}_{t,x}^{\prime}\left(  (0,T)\times
\mathbb{R}_{per}^{d}\right)  $. Next, we take the difference between $\left(
\text{\ref{renorm_sans_eps}}\right)  $ and $\left(  \text{\ref{passage_limite}%
}\right)  $, we multiply it with $\gamma-1$ in order to obtain that%
\begin{equation}
\partial_{t}\Theta+\operatorname{div}\left(  \Theta u\right)  +\left(
\gamma-1\right)  \Theta\operatorname{div}u=-\left(  \gamma-1\right)  \Xi\text{
in }\mathcal{D}_{t,x}^{\prime}\left(  (0,T)\times\mathbb{R}_{per}^{d}\right)
, \label{difference1}%
\end{equation}
where
\begin{align*}
&  \Theta\overset{not.}{=}\overline{\rho^{\gamma}}-\rho^{\gamma},\\
&  \Xi\overset{not.}{=}\left(  \mu\overline{\nabla u:\nabla u}+\left(
\mu+\lambda\right)  \overline{(\operatorname{div}u)^{2}}+\overline
{\mathcal{E}\left(  \nabla u\right)  :\nabla u}\right) \\
&  \text{ \ \ \ \ \ }-\left(  \mu\nabla u:\nabla u+\left(  \mu+\lambda\right)
(\operatorname{div}u)^{2}+\mathcal{E}\left(  \nabla u\right)  :\nabla
u\right)  .
\end{align*}
Obviously,
\[
\Theta,\Xi\geq0.
\]
We regularize the previous equation with the help of a sequence of
approximations of the identity $\omega_{\varepsilon}:$
\begin{equation}
\partial_{t}\omega_{\varepsilon}\ast\Theta+\operatorname{div}\left(
\omega_{\varepsilon}\ast\Theta u\right)  +\left(  \gamma-1\right)
\omega_{\varepsilon}\ast(\Theta\operatorname{div}u)=r_{\varepsilon}\left(
\Theta,u\right)  -\left(  \gamma-1\right)  \omega_{\varepsilon}\ast\Xi,
\label{reg}%
\end{equation}
see the notations introduced in $\left(  \text{\ref{notation_approx}}\right)
$ and $\left(  \text{\ref{def_reminder}}\right)  $. Since the time derivative
$\partial_{t}\omega_{\varepsilon}\ast\Theta$ belongs to some Lebesgue space,
we may multiply relation $\left(  \text{\ref{reg}}\right)  $ with $\frac
{1}{\gamma}(h+\omega_{\varepsilon}\ast\Theta)^{\frac{1}{\gamma}-1}$ where
$h>0$ is a fixed positive constant and apply the chain rule. We end up with%
\begin{align}
&  \partial_{t}\left(  h+\omega_{\varepsilon}\ast\Theta\right)  ^{\frac
{1}{\gamma}}+\operatorname{div}\left(  \left(  h+\omega_{\varepsilon}%
\ast\Theta\right)  ^{\frac{1}{\gamma}}u\right)  +(h+\omega_{\varepsilon}%
\ast\Theta)^{\frac{1}{\gamma}-1}[\left(  \frac{1}{\gamma}-1\right)
\omega_{\varepsilon}\ast\Theta-h]\operatorname{div}u\nonumber\\
&  +\left(  1-\frac{1}{\gamma}\right)  (h+\omega_{\varepsilon}\ast
\Theta)^{\frac{1}{\gamma}-1}\omega_{\varepsilon}\ast(\Theta\operatorname{div}%
u)\nonumber\\
&  =\frac{1}{\gamma}(h+\omega_{\varepsilon}\ast\Theta)^{\frac{1}{\gamma}%
-1}r_{\varepsilon}\left(  \Theta,u\right)  -\frac{1}{\gamma}(h+\omega
_{\varepsilon}\ast\Theta)^{\frac{1}{\gamma}-1}\left(  \gamma-1\right)
\omega_{\varepsilon}\ast\Xi, \label{derivative}%
\end{align}
Owing to $\left(  \text{\ref{continuity_rho}}\right)  $ and $\left(
\text{\ref{continuity_rho_gamma}}\right)  $ the application $t\rightarrow
\int_{\mathbb{T}^{d}}\left(  h+\omega_{\varepsilon}\ast\Theta\right)
^{\frac{1}{\gamma}}$ is continuous and since by integrating w.r.t. space in
$\left(  \text{\ref{derivative}}\right)  $, its distributional time derivative
belongs to some Lebesgue space, we deduce that it is absolutely continuous and
that the distributional derivative coincides with the derivative a.e.. We may
thus write that for any $t\in\left(  0,T\right)  $ we have that%
\begin{align*}
&  \int_{\mathbb{T}^{d}}\left(  h+\omega_{\varepsilon}\ast\Theta\right)
^{\frac{1}{\gamma}}\left(  t\right) \\
&  =\int_{\mathbb{T}^{d}}\left(  h+\omega_{\varepsilon}\ast\Theta\right)
^{\frac{1}{\gamma}}\left(  0\right)  -\int_{0}^{t}\int_{\mathbb{T}^{d}}\left(
\frac{1}{\gamma}-1\right)  (h+\omega_{\varepsilon}\ast\Theta)^{\frac{1}%
{\gamma}-1}\left[  \omega_{\varepsilon},\operatorname{div}u\right]
\Theta+\int_{0}^{t}\int_{\mathbb{T}^{d}}(h+\omega_{\varepsilon}\ast
\Theta)^{\frac{1}{\gamma}-1}h\operatorname{div}u\\
&  +\int_{0}^{t}\int_{\mathbb{T}^{d}}\left[  \frac{1}{\gamma}(h+\omega
_{\varepsilon}\ast\Theta)^{\frac{1}{\gamma}-1}r_{\varepsilon}\left(
\Theta,u\right)  -\frac{1}{\gamma}(h+\omega_{\varepsilon}\ast\Theta)^{\frac
{1}{\gamma}-1}\left(  \gamma-1\right)  \omega_{\varepsilon}\ast\Xi\right] \\
&  \leq\int_{\mathbb{T}^{d}}\left(  h+\omega_{\varepsilon}\ast\Theta\right)
^{\frac{1}{\gamma}}\left(  0\right)  -\int_{0}^{t}\left(  \frac{1}{\gamma
}-1\right)  (h+\omega_{\varepsilon}\ast\Theta)^{\frac{1}{\gamma}-1}\left[
\omega_{\varepsilon},\operatorname{div}u\right]  \Theta\\
&  +\int_{0}^{t}\int_{\mathbb{T}^{d}}\frac{1}{\gamma}(h+\omega_{\varepsilon
}\ast\Theta)^{\frac{1}{\gamma}-1}r_{\varepsilon}\left(  \Theta,u\right)  ,
\end{align*}
where we used the positivity of $\Xi$. Using Proposition \ref{Prop_ren1}, we
obtain that
\[
\left[  \omega_{\varepsilon},\operatorname{div}u\right]  \Theta\text{ and
}r_{\varepsilon}\left(  \Theta,u\right)  \rightarrow0\text{ in }L^{1}\left(
\left(  0,T\right)  \times\mathbb{T}^{d}\right)  .
\]
Notice that since $\gamma>1$ along with $\omega_{\varepsilon}\ast\Theta\geq0,$
we also have that%
\[
(h+\omega_{\varepsilon}\ast\Theta)^{1/\gamma-1}\leq h^{1/\gamma-1}.
\]
Taking into account the last observations, by making $\varepsilon\rightarrow0$
we get that
\[
\int_{\mathbb{T}^{d}}\left(  h+\Theta\right)  ^{\frac{1}{\gamma}}\left(
t\right)  \leq\int_{\mathbb{T}^{d}}\left(  h+\Theta\right)  ^{\frac{1}{\gamma
}}\left(  0\right)  +h^{1/\gamma}\int_{0}^{t}\int_{{\mathbb{T}}^{d}%
}|\operatorname{div}\,u|.
\]
Letting $h$ go to zero and using that at initial time $\Theta\left(  0\right)
$ is $0,$ shows that%
\[
\int_{\mathbb{T}^{d}}\Theta^{\frac{1}{\gamma}}\left(  t\right)  =\int
_{\mathbb{T}^{d}}\left(  \overline{\rho^{\gamma}}\left(  t\right)
-\rho^{\gamma}\left(  t\right)  \right)  ^{\frac{1}{\gamma}}=0,
\]
from which follows the conclusion that%
\[
\overline{\rho^{\gamma}}=\rho^{\gamma}\text{ a.e. on }\left(  0,T\right)
\times\mathbb{T}^{d}.
\]
This ends the proof of Theorem \ref{MainTheorem}.

\section{Appendix}

\subsection{Appendix A: tool box}

In this section, we gather some classical results that are used throught the text.

\begin{lemma}
[Fourier Multipliers]Consider $m:\mathbb{R}^{d}\backslash\left\{  0\right\}
\rightarrow\mathbb{R}$ a function verifying
\[
\left\vert \partial^{\alpha}m\left(  \xi\right)  \right\vert \leq c_{\alpha
}\left\vert \xi\right\vert ^{-\alpha}%
\]
for all $\alpha\in\mathbb{N}^{d}$ with $\left\vert \alpha\right\vert \leq
d+1$. Then, for all $p\in\left(  1,\infty\right)  $, there exists $C_{p}$ such
that for any $u\in L^{p}$
\[
\left\Vert \mathcal{F}^{-1}\left(  m\left(  \xi\right)  \mathcal{F}\left(
u-\int_{\mathbb{T}^{d}}u\right)  \right)  \right\Vert _{L^{p}(\mathbb{T}^{d}%
)}\leq C_{p}\left\Vert u\right\Vert _{L^{p}(\mathbb{T}^{d})}.
\]

\end{lemma}

\begin{lemma}
[Sobolev's Inequality]For all $p\in(1,d)$ and $u\in L^{p}\cap D^{1,p^{\ast}}$
we have that%
\[
\left\Vert u-\int_{\mathbb{T}^{d}}u\right\Vert _{L^{p}(\mathbb{T}^{d})}%
\leq\left\Vert \nabla u\right\Vert _{L^{p^{\ast}}(\mathbb{T}^{d})},
\]
where $1/p+1/d=1/p^{\ast}$.
\end{lemma}

\bigskip Let $X$ be a Banach space. We consider $\mathcal{C}^{0}(\left[
0,T\right]  ;X_{weak})$ the space of continuous functions from $\left[
0,T\right]  $ to $X$ endowed with the weak topology:
\[
\mathcal{C}^{0}(\left[  0,T\right]  ;X_{weak})=\left\{
\begin{array}
[c]{l}%
f:\left[  0,T\right]  \rightarrow X\text{ such that }\\
\forall w\in X^{\prime}\text{, }t\rightarrow\left\langle w,f\left(  t\right)
\right\rangle _{X^{\prime}\times X}\text{ is continuous}%
\end{array}
\right\}
\]

\begin{definition}
A subset $\mathcal{F}$ of $\mathcal{C}^{0}(\left[  0,T\right]  ;X_{weak})$ is
called weakly equicontinious on $\left[  0,T\right]  $ if for all $w\in
X^{\prime}$ and for all $\varepsilon>0$ there exists a $\delta=\delta\left(
w,\varepsilon\right)  >0$ such that for all $t,s\in\left[  0,T\right]  $%
\[
\left\vert t-s\right\vert \leq\delta\Rightarrow\left\vert \left\langle
f\left(  t\right)  -f\left(  s\right)  ,w\right\rangle \right\vert
\leq\varepsilon.
\]

\end{definition}

The following theorem is a version of the Arzel\`{a}-Ascoli theorem:

\begin{theorem}
\label{Arzela-Ascoli-weak}Let $X$ be a reflexive Banach space. A subset
$\mathcal{F}$ of $\mathcal{C}^{0}(\left[  0,T\right]  ;X_{weak})$ endowed with
the uniformly weak topology is sequentially relatively compact if and only if

\begin{itemize}
\item $\mathcal{F}$ is weakly equicontinious on $\left[  0,T\right]  $.

\item There exists $D\subset\left[  0,T\right]  $ dense such that for all
$t\in D$ the set $\mathcal{F}\left(  t\right)  :=\left\{  f\left(  t\right)
:f\in\mathcal{F}\right\}  $ is bounded in $X$.
\end{itemize}
\end{theorem}

For a proof of a slightly more general result see Theorem $A.3.1.$ from Vrabie
\cite{Vrabie2003Semigroups}, page $302$.

Let $g\in L^{q}((0,T);L^{p}(\mathbb{T}^{d}))$ with $p,q\geq1$, introduce a new
function
\begin{equation}
g_{\delta}\left(  x\right)  =g\ast\omega_{\delta}(x)\qquad\hbox{ with }\qquad
\omega_{\delta}\left(  x\right)  =\frac{1}{\delta^{d}}\omega(\frac{x}{\delta})
\label{notation_approx}%
\end{equation}
with $\omega$ a smooth, nonnegative, even function compactly supported in the
unit ball centered at the origin and with integral equal to 1. We recall the
following classical analysis result
\[
\lim_{\delta\rightarrow0}\left\Vert g_{\delta}-g\right\Vert _{L^{q}%
(0,T;L^{p}({\mathbb{T}}^{d}))}=0.
\]
Next let us recall the following commutator estimate which was obtained for
the first time by DiPerna and Lions:

\begin{proposition}
\label{Prop_ren1}Consider $\beta\in(1,\infty)$ and $\left(  a,b\right)  $ such
that $a\in L^{\beta}\left(  \left(  0,T\right)  \times\mathbb{T}^{d}\right)  $
and $b,\nabla b\in L^{p}\left(  \left(  0,T\right)  \times\mathbb{T}%
^{d}\right)  $ where $\frac{1}{s}=\frac{1}{\beta}+\frac{1}{p}\leq1$. Then, we
have%
\[
\lim_{\delta\rightarrow0}r_{\delta}^{k}\left(  a,b\right)  =0\text{ in }%
L^{s}\left(  \left(  0,T\right)  \times\mathbb{T}^{d}\right)  ,
\]
for $k\in\left\{  1,2\right\}  $ where
\begin{equation}
r_{\delta}^{1}\left(  a,b\right)  =b\partial_{i}a_{\delta}-\left(
b\partial_{i}a\right)  _{\delta}\text{ and }r_{\delta}^{2}\left(  a,b\right)
=\partial_{i}\left(  a_{\delta}b\right)  -\partial_{i}\left(  \left(
ab\right)  _{\delta}\right)  .\text{ } \label{def_reminder}%
\end{equation}
Moreover, the following commutator estimates hold true%
\begin{align}
\left\Vert b\partial_{i}a_{\delta}-\left(  b\partial_{i}a\right)  _{\delta
}\right\Vert _{L_{t}^{s}L_{x}^{s}}  &  \leq\left\Vert \nabla b\right\Vert
_{L_{t}^{p}L_{x}^{p}}\left\Vert a\right\Vert _{L_{t}^{\beta}L_{x}^{\beta}%
}\label{comut_estim1}\\
\left\Vert \partial_{i}\left(  a_{\delta}b\right)  -\partial_{i}\left(
\left(  ab\right)  _{\delta}\right)  \right\Vert _{L_{t}^{s}L_{x}^{s}}  &
\leq\left\Vert \nabla b\right\Vert _{L_{t}^{p}L_{x}^{p}}\left\Vert
a\right\Vert _{L_{t}^{\beta}L_{x}^{\beta}} \label{comut_estim2}%
\end{align}
where $b\partial_{i}a$ should be understood as%
\[
b\partial_{i}a=\partial_{i}\left(  ab\right)  -a\partial_{i}b.
\]

\end{proposition}

Whenever we have a \textit{regular solution} for the transport equation
\begin{equation}
\partial_{t}\rho+\operatorname{div}\left(  \rho u\right)  =0,
\label{transport_eq1}%
\end{equation}
then, multiplying the former equation with $b^{\prime}\left(  \rho\right)  $
gives
\begin{equation}
\partial_{t}b\left(  \rho\right)  +\operatorname{div}\left(  b\left(
\rho\right)  u\right)  +\left\{  \rho b^{\prime}\left(  \rho\right)  -b\left(
\rho\right)  \right\}  \operatorname{div}u=0. \label{renorm1}%
\end{equation}
The following proposition gives us a framework for justifying this
computations when $\rho$ is just a Lebesgue function.

\begin{proposition}
\label{Prop_ren2}Consider $2\leq\beta<\infty$ and $\lambda_{0},\lambda_{1}$
such that $\lambda_{0}<1$ and $-1\leq\lambda_{1}\leq\beta/2-1$. Also, consider
$\rho\in L^{\beta}\left(  \left(  0,T\right)  \times\mathbb{T}^{3}\right)  $,
$\rho\geq0$ a.e. and $u,\nabla u\in L^{2}\left(  \left(  0,T\right)
\times\mathbb{T}^{3}\right)  $ verifying the transport equation $\left(
\text{\ref{transport_eq1}}\right)  $ in the sense of distributions. Then, for
any function $b\in C^{0}\left(  [0,\infty)\right)  \cap C^{1}\left(  \left(
0,\infty\right)  \right)  $ such that%
\[
\left\{
\begin{array}
[c]{l}%
b^{\prime}\left(  t\right)  \leq ct^{-\lambda_{0}}\text{ for }t\in(0,1],\\
\left\vert b^{\prime}\left(  t\right)  \right\vert \leq ct^{\lambda_{1}}\text{
for }t\geq1.
\end{array}
\right.
\]
Then, equation $\left(  \text{\ref{renorm1}}\right)  $ holds in the sense of distributions.
\end{proposition}

\noindent\ The proof of the above results follow by adapting in a
straightforward manner lemmas $6.7.$ and $6.9$ from the book of A.
Novotn\'{y}- I.Stra\v{s}kraba \cite{NoSt} pages 304--308.

\subsection{Appendix B: detailed computations for the Hoff functionals}

\subsubsection{Hoff's first energy functional}

The momentum equation reads:%
\[
\rho\dot{u}-\mu\Delta u-(\mu+\lambda)\nabla\operatorname{div}%
u-\operatorname{div}\omega_{\delta}\ast\mathcal{E}\left(  \nabla u_{\delta
}\right)  +\nabla P\left(  \rho\right)  =\rho f.
\]
where
\[
\dot{u}=\partial_{t}u+u\nabla u.
\]
We multiply the above equation with $\dot{u}$ and integrate. Owing to the
hypothesis%
\[
\varepsilon_{ijkl}a_{ij}b_{kl}=\varepsilon_{ijkl}a_{kl}b_{ij}%
\]
we can write that%

\begin{align*}
\left\langle \operatorname{div}\omega_{\delta}\ast\mathcal{E}\left(
\omega_{\delta}\ast\nabla u\right)  ,\dot{u}\right\rangle  &  =-\int
_{\mathbb{T}^{d}}\partial_{j}(\varepsilon_{ijkl}\partial_{\ell}u_{\delta}%
^{k})\dot{u}_{\delta}^{i}\\
&  =\int_{\mathbb{T}^{d}}\varepsilon_{ijkl}\partial_{\ell}u_{\delta}%
^{k}\partial_{j}\partial_{t}u_{\delta}^{i}+\int_{\mathbb{T}^{d}}%
\varepsilon_{ijkl}\partial_{\ell}u_{\delta}^{k}\omega_{\delta}\ast
(\partial_{j}u^{q}\partial_{q}u^{i})+\int_{\mathbb{T}^{d}}\varepsilon
_{ijkl}\partial_{\ell}u_{\delta}^{k}\omega_{\delta}\ast(u^{q}\partial_{qj}%
^{2}u^{i})\\
&  =\frac{1}{2}\left\{  \int_{\mathbb{T}^{d}}\varepsilon_{ijkl}\partial_{\ell
}u_{\delta}^{k}\partial_{j}\partial_{t}u_{\delta}^{i}+\int_{\mathbb{T}^{d}%
}\varepsilon_{ijkl}\partial_{t}\partial_{\ell}u_{\delta}^{k}\partial
_{j}u_{\delta}^{i}\right\}  +\int_{\mathbb{T}^{d}}\varepsilon_{ijkl}%
\partial_{\ell}u_{\delta}^{k}\omega_{\delta}\ast(\partial_{j}u^{q}\partial
_{q}u^{i})\\
&  +\int_{\mathbb{T}^{d}}\varepsilon_{ijkl}\partial_{\ell}u_{\delta}^{k}%
u^{q}\partial_{qj}^{2}u_{\delta}^{i}+\int_{\mathbb{T}^{d}}\varepsilon
_{ijkl}\partial_{\ell}u_{\delta}^{k}\left[  u^{q},\omega_{\delta}\right]
\partial_{qj}^{2}u^{i}\\
&  =\frac{1}{2}\frac{d}{dt}\int_{\mathbb{T}^{d}}\varepsilon_{ijkl}%
\partial_{\ell}u_{\delta}^{k}\partial_{j}u_{\delta}^{i}-\frac{1}{2}%
\int_{\mathbb{T}^{d}}\partial_{t}\varepsilon_{ijkl}\partial_{\ell}u_{\delta
}^{k}\partial_{j}u_{\delta}^{i}+\int_{\mathbb{T}^{d}}\varepsilon
_{ijkl}\partial_{\ell}u_{\delta}^{k}\omega_{\delta}\ast(\partial_{j}%
u^{q}\partial_{q}u^{i})\\
&  +\frac{1}{2}\left\{  \int_{\mathbb{T}^{d}}\varepsilon_{ijkl}\partial
_{j}u_{\delta}^{i}u^{q}\partial_{q\ell}^{2}u_{\delta}^{k}+\int_{\mathbb{T}%
^{d}}\varepsilon_{ijkl}\partial_{\ell}u^{k}u^{q}\partial_{qj}^{2}u_{\delta
}^{i}\right\}  +\int_{\mathbb{T}^{d}}\varepsilon_{ijkl}\partial_{\ell
}u_{\delta}^{k}\left[  u^{q},\omega_{\delta}\right]  \partial_{qj}^{2}u^{i}\\
&  =\frac{1}{2}\frac{d}{dt}\int_{\mathbb{T}^{d}}\varepsilon_{ijkl}%
\partial_{\ell}u^{k}\partial_{j}u^{i}+\int_{\mathbb{T}^{d}}\varepsilon
_{ijkl}\partial_{\ell}u_{\delta}^{k}\omega_{\delta}\ast(\partial_{j}%
u^{q}\partial_{q}u^{i})-\frac{1}{2}\int_{\mathbb{T}^{d}}\left\{  \partial
_{t}\varepsilon_{ijkl}+\partial_{q}(\varepsilon_{ijkl}u^{q})\right\}
\partial_{j}u_{\delta}^{i}\partial_{\ell}u_{\delta}^{k}\\
&  +\int_{\mathbb{T}^{d}}\varepsilon_{ijkl}\partial_{\ell}u_{\delta}%
^{k}\left[  u^{q},\omega_{\delta}\right]  \partial_{qj}^{2}u^{i}.
\end{align*}
Similar computations show that%
\begin{align*}
-\left\langle \left(  \mu\Delta+\left(  \mu+\lambda\right)  \nabla
\operatorname{div}\right)  u,\dot{u}\right\rangle  &  =\frac{1}{2}\frac{d}%
{dt}\left\{  \mu\int_{\mathbb{T}^{d}}\left\vert \partial_{k}u^{i}\right\vert
^{2}+\left(  \mu+\lambda\right)  \int_{\mathbb{T}^{d}}\left\vert
\operatorname{div}u\right\vert ^{2}\right\} \\
&  +\mu\int_{\mathbb{T}^{d}}\partial_{k}u^{i}\partial_{k}u^{\ell}%
\partial_{\ell}u^{i}-\frac{\mu}{2}\int_{\mathbb{T}^{d}}\left\vert \partial
_{k}u^{i}\right\vert ^{2}\operatorname{div}u\\
&  +\left(  \mu+\lambda\right)  \int_{\mathbb{T}^{d}}\operatorname{div}%
u\partial_{i}u^{\ell}\partial_{\ell}u^{i}-\frac{\mu+\lambda}{2}\int
_{\mathbb{T}^{d}}(\operatorname{div}u)^{3}%
\end{align*}
Next, we treat the pressure term%
\begin{align*}
\int_{\mathbb{T}^{d}}\dot{u}\nabla P\left(  \rho\right)   &  =-\int
_{\mathbb{T}^{d}}P\left(  \rho\right)  \operatorname{div}\dot{u}=-\frac{d}%
{dt}\left\{  \int_{\mathbb{T}^{d}}P\left(  \rho\right)  \operatorname{div}%
u\right\}  +\int_{\mathbb{T}^{d}}\partial_{t}P\left(  \rho\right)
\operatorname{div}u-\int_{\mathbb{T}^{d}}P\left(  \rho\right)
\operatorname{div}\left(  u\nabla u\right) \\
&  =-\frac{d}{dt}\left\{  \int_{\mathbb{T}^{d}}P\left(  \rho\right)
\operatorname{div}u\right\}  +\int_{\mathbb{T}^{d}}\partial_{t}P\left(
\rho\right)  \operatorname{div}u-\int_{\mathbb{T}^{d}}P\left(  \rho\right)
\partial_{\ell}u^{k}\partial_{k}u^{\ell}-\int_{\mathbb{T}^{d}}P\left(
\rho\right)  u^{k}\partial_{k\ell}^{2}u^{\ell}\\
&  =-\frac{d}{dt}\left\{  \int_{\mathbb{T}^{d}}P\left(  \rho\right)
\operatorname{div}u\right\}  +\int_{\mathbb{T}^{d}}\partial_{t}P\left(
\rho\right)  \operatorname{div}u-\int_{\mathbb{T}^{d}}P\left(  \rho\right)
\partial_{\ell}u^{k}\partial_{k}u^{\ell}+\int_{\mathbb{T}^{d}}\partial
_{k}(P\left(  \rho\right)  u^{k})\partial_{\ell}u^{\ell}\\
&  =-\frac{d}{dt}\left\{  \int_{\mathbb{T}^{d}}P\left(  \rho\right)
\operatorname{div}u\right\}  +\int_{\mathbb{T}^{d}}\left(  \partial
_{t}P\left(  \rho\right)  +\operatorname{div}\left(  P(\rho)u\right)  \right)
\operatorname{div}u-\int_{\mathbb{T}^{d}}P\left(  \rho\right)  \partial_{\ell
}u^{k}\partial_{k}u^{\ell}\\
&  =-\frac{d}{dt}\left\{  \int_{\mathbb{T}^{d}}P\left(  \rho\right)
\operatorname{div}u\right\}  +\int_{\mathbb{T}^{d}}\left(  P\left(
\rho\right)  -\rho P^{\prime}\left(  \rho\right)  \right)  (\operatorname{div}%
u)^{2}-\int_{\mathbb{T}^{d}}P\left(  \rho\right)  \partial_{\ell}u^{k}%
\partial_{k}u^{\ell}.
\end{align*}
Putting together all the above computations, we end up with%

\begin{align}
&  \frac{1}{2}\frac{d}{dt}\left\{  \mu\int_{\mathbb{T}^{d}}\left\vert
\partial_{k}u^{i}\right\vert ^{2}+\left(  \mu+\lambda\right)  \int
_{\mathbb{T}^{d}}\left\vert \operatorname{div}u\right\vert ^{2}+\int
_{\mathbb{T}^{d}}\varepsilon_{ijkl}\partial_{\ell}u_{\delta}^{k}\partial
_{j}u_{\delta}^{i}-\int_{\mathbb{T}^{d}}P\left(  \rho\right)
\operatorname{div}u\right\}  +\int_{\mathbb{T}^{d}}\rho\left\vert \dot
{u}\right\vert ^{2}\nonumber\\
&  =-\mu\int_{\mathbb{T}^{d}}\partial_{k}u^{i}\partial_{k}u^{\ell}%
\partial_{\ell}u^{i}+\frac{\mu}{2}\int_{\mathbb{T}^{d}}\left\vert \partial
_{k}u^{i}\right\vert ^{2}\operatorname{div}u\nonumber\\
&  -\left(  \mu+\lambda\right)  \int_{\mathbb{T}^{d}}\operatorname{div}%
u\partial_{i}u^{\ell}\partial_{\ell}u^{i}+\frac{\mu+\lambda}{2}\int
_{\mathbb{T}^{d}}(\operatorname{div}u)^{3}\nonumber\\
&  +\frac{1}{2}\int_{\mathbb{T}^{d}}\left\{  \partial_{t}\varepsilon
_{ijkl}+\partial_{q}(\varepsilon_{ijkl}u^{q})\right\}  \partial_{j}u_{\delta
}^{i}\partial_{k}u_{\delta}^{\ell}-\int_{\mathbb{T}^{d}}\varepsilon
_{ijkl}\partial_{\ell}u_{\delta}^{k}\omega_{\delta}\ast(\partial_{j}%
u^{q}\partial_{q}u^{i})-\int_{\mathbb{T}^{d}}\varepsilon_{ijkl}\partial_{\ell
}u_{\delta}^{k}\left[  u^{q},\omega_{\delta}\right]  \partial_{qj}^{2}%
u^{i}.\nonumber\\
&  +\int_{0}^{t}\int_{\mathbb{T}^{d}}\rho P^{\prime}\left(  \rho\right)
\partial_{\ell}u^{k}\partial_{k}u^{\ell}+\int_{\mathbb{T}^{d}}\left(  \rho
P^{\prime}\left(  \rho\right)  -P\left(  \rho\right)  \right)
(\operatorname{div}u)^{2}+\int_{0}^{t}\int_{\mathbb{T}^{d}}\rho\dot{u}f.
\label{2}%
\end{align}

\subsubsection{Hoff's second energy functional}

The idea leading to the construction of this second functional is to apply to
the momentum equation the material time derivative $\partial_{t}%
\cdot+\operatorname{div}\left(  u\cdot\right)  $, multiply with $\dot{u}$ and
integrate. The detailed computations are presented below. First, we obviously
have that
\[
\int_{\mathbb{T}^{d}}\left(  \partial_{t}(\rho\dot{u}^{j})+\partial_{k}%
(u^{k}\rho\dot{u}^{j})\right)  \dot{u}^{j}=\frac{d}{dt}\int_{\mathbb{T}^{d}%
}\frac{\rho\left\vert \dot{u}\right\vert ^{2}}{2}%
\]
Next, let us deal with the pressure term. First of all, owing to the density
equation we write that%
\[
\partial_{t}P\left(  \rho\right)  +\operatorname{div}\left(  P\left(
\rho\right)  u\right)  +\left(  \rho P^{\prime}\left(  \rho\right)  -P\left(
\rho\right)  \right)  \operatorname{div}u=0
\]
which implies that for all $j\in\overline{1,d}$ it holds true that%
\[
\partial_{t}\partial_{j}P\left(  \rho\right)  +\operatorname{div}\left(
\partial_{j}P\left(  \rho\right)  u\right)  +\operatorname{div}\left(
P\left(  \rho\right)  \partial_{j}u\right)  +\partial_{j}\left\{  \left(  \rho
P^{\prime}\left(  \rho\right)  -P\left(  \rho\right)  \right)
\operatorname{div}u\right\}  =0.
\]
We use this relation in order to infer
\begin{align*}
\int_{\mathbb{T}^{d}}\left(  \partial_{t}\partial_{j}P\left(  \rho\right)
+\partial_{k}(u^{k}\partial_{j}P\left(  \rho\right)  )\right)  \dot{u}^{j}  &
=-\int_{\mathbb{T}^{d}}\left\{  \operatorname{div}\left(  P\left(
\rho\right)  \partial_{j}u\right)  +\partial_{j}\left\{  \left(  \rho
P^{\prime}\left(  \rho\right)  -P\left(  \rho\right)  \right)
\operatorname{div}u\right\}  \right\}  \dot{u}^{j}\\
&  =\int_{\mathbb{T}^{d}}\left\{  P\left(  \rho\right)  \partial_{j}%
u^{k}\partial_{k}\dot{u}^{j}+\left(  \rho P^{\prime}\left(  \rho\right)
-P\left(  \rho\right)  \right)  \operatorname{div}u\operatorname{div}\dot
{u}\right\}  .
\end{align*}
Finally, let us treat the dissipative term. We observe that%
\begin{align*}
&  -\left\langle \partial_{t}\operatorname{div}\omega_{\delta}\ast
\mathcal{E}\left(  \nabla u_{\delta}\right)  +\operatorname{div}\left(
u\operatorname{div}\omega_{\delta}\ast\mathcal{E}\left(  \nabla u_{\delta
}\right)  \right)  ,\dot{u}\right\rangle \\
&  -\int_{\mathbb{T}^{d}}\partial_{j}\left(  \partial_{t}\varepsilon
_{ijkl}\partial_{\ell}u_{\delta}^{k}\right)  \dot{u}_{\delta}^{i}%
-\int_{\mathbb{T}^{d}}\partial_{j}\left(  \varepsilon_{ijkl}\partial
_{t}\partial_{\ell}u_{\delta}^{k}\right)  \dot{u}_{\delta}^{i}-\int
_{\mathbb{T}^{d}}\partial_{q}\left(  u^{q}\omega_{\delta}\ast\partial
_{j}\left(  \varepsilon_{ijkl}\partial_{\ell}u_{\delta}^{k}\right)  \right)
\dot{u}^{i}\\
& \\
&  =\int_{\mathbb{T}^{d}}\partial_{t}\varepsilon_{ijkl}\partial_{\ell
}u_{\delta}^{k}\partial_{j}\dot{u}_{\delta}^{i}+\int_{\mathbb{T}^{d}%
}\varepsilon_{ijkl}\partial_{t}\partial_{\ell}u_{\delta}^{k}\partial_{j}%
\dot{u}_{\delta}^{i}+\int_{\mathbb{T}^{d}}\partial_{j}\left(  \varepsilon
_{ijkl}\partial_{\ell}u_{\delta}^{k}\right)  \omega_{\delta}\ast(u^{q}%
\partial_{q}\dot{u}^{i})\\
& \\
&  =\int_{\mathbb{T}^{d}}\partial_{t}\varepsilon_{ijkl}\partial_{\ell
}u_{\delta}^{k}\partial_{j}\dot{u}_{\delta}^{i}\\
&  +\int_{\mathbb{T}^{d}}\varepsilon_{ijkl}\partial_{\ell}(\partial
_{t}u_{\delta}^{k}+\omega_{\delta}\ast(u^{q}\partial_{q}u^{k}))\partial
_{j}\dot{u}_{\delta}^{i}-\int_{\mathbb{T}^{d}}\varepsilon_{ijkl}\partial
_{\ell}(\omega_{\delta}\ast(u^{q}\partial_{q}u^{k}))\partial_{j}\dot
{u}_{\delta}^{i}\\
&  +\int_{\mathbb{T}^{d}}\partial_{j}\left(  \varepsilon_{ijkl}\partial_{\ell
}u_{\delta}^{k}\right)  \omega_{\delta}\ast(u^{q}\partial_{q}\dot{u}^{i})\\
& \\
&  =\int_{\mathbb{T}^{d}}\partial_{t}\varepsilon_{ijkl}\partial_{\ell
}u_{\delta}^{k}\partial_{j}\dot{u}_{\delta}^{i}+\int_{\mathbb{T}^{d}%
}\varepsilon_{ijkl}\partial_{\ell}\dot{u}_{\delta}^{k}\partial_{j}\dot
{u}_{\delta}^{i}\\
&  -\int_{\mathbb{T}^{d}}\varepsilon_{ijkl}(\omega_{\delta}\ast(\partial
_{\ell}u^{q}\partial_{q}u^{k}))\partial_{j}\dot{u}_{\delta}^{i}-\int
_{\mathbb{T}^{d}}\varepsilon_{ijkl}(\omega_{\delta}\ast(u^{q}\partial_{\ell
q}^{2}u^{k}))\partial_{j}\dot{u}_{\delta}^{i}\\
&  -\int_{\mathbb{T}^{d}}\varepsilon_{ijkl}\partial_{\ell}u_{\delta}^{k}%
\omega_{\delta}\ast(\partial_{j}u^{q}\partial_{q}\dot{u}^{i})-\int
_{\mathbb{T}^{d}}\varepsilon_{ijkl}\partial_{\ell}u_{\delta}^{k}\omega
_{\delta}\ast(u^{q}\partial_{jq}^{2}\dot{u}^{i}).
\end{align*}
Again, integrating by parts leads to the following identity:
\begin{align*}
&  -\left\langle \partial_{t}\operatorname{div}\omega_{\delta}\ast
\mathcal{E}\left(  \nabla u_{\delta}\right)  +\operatorname{div}\left(
u\operatorname{div}\omega_{\delta}\ast\mathcal{E}\left(  \nabla u_{\delta
}\right)  \right)  ,\dot{u}\right\rangle \\
& \\
&  =\int_{\mathbb{T}^{d}}\partial_{t}\varepsilon_{ijkl}\partial_{\ell
}u_{\delta}^{k}\partial_{j}\dot{u}_{\delta}^{i}+\int_{\mathbb{T}^{d}%
}\varepsilon_{ijkl}\partial_{\ell}\dot{u}_{\delta}^{k}\partial_{j}\dot
{u}_{\delta}^{i}\\
&  -\int_{\mathbb{T}^{d}}\varepsilon_{ijkl}(\omega_{\delta}\ast(\partial
_{\ell}u^{q}\partial_{q}u^{k}))\partial_{j}\dot{u}_{\delta}^{i}-\int
_{\mathbb{T}^{d}}\varepsilon_{ijkl}\partial_{\ell}u_{\delta}^{k}\omega
_{\delta}\ast(\partial_{j}u^{q}\partial_{q}\dot{u}^{i})\\
&  -\int_{\mathbb{T}^{d}}\varepsilon_{ijkl}(\left[  u^{q},\omega_{\delta}%
\ast\right]  \partial_{\ell q}^{2}u^{k})\partial_{j}\dot{u}_{\delta}^{i}%
-\int_{\mathbb{T}^{d}}\varepsilon_{ijkl}u^{q}\partial_{\ell q}^{2}u_{\delta
}^{k}\partial_{j}\dot{u}_{\delta}^{i}\\
&  -\int_{\mathbb{T}^{d}}\varepsilon_{ijkl}\partial_{\ell}u_{\delta}%
^{k}(\left[  u^{q},\omega_{\delta}\ast\right]  \partial_{jq}^{2}\dot{u}%
^{i})-\int_{\mathbb{T}^{d}}\varepsilon_{ijkl}\partial_{\ell}u_{\delta}%
^{k}u^{q}\partial_{jq}^{2}\dot{u}_{\delta}^{i}\\
& \\
&  =\int_{\mathbb{T}^{d}}\left(  \partial_{t}\varepsilon_{ijkl}+\partial
_{q}\left(  u^{q}\varepsilon_{ijkl}\right)  \right)  \partial_{\ell}u_{\delta
}^{k}\partial_{j}\dot{u}_{\delta}^{i}+\int_{\mathbb{T}^{d}}\varepsilon
_{ijkl}\partial_{\ell}\dot{u}_{\delta}^{k}\partial_{j}\dot{u}_{\delta}^{i}\\
&  -\int_{\mathbb{T}^{d}}\varepsilon_{ijkl}(\omega_{\delta}\ast(\partial
_{\ell}u^{q}\partial_{q}u^{k}))\partial_{j}\dot{u}_{\delta}^{i}-\int
_{\mathbb{T}^{d}}\varepsilon_{ijkl}\partial_{\ell}u_{\delta}^{k}\omega
_{\delta}\ast(\partial_{j}u^{q}\partial_{q}\dot{u}^{i})\\
&  -\int_{\mathbb{T}^{d}}\varepsilon_{ijkl}(\left[  u^{q},\omega_{\delta}%
\ast\right]  \partial_{\ell q}^{2}u^{k})\partial_{j}\dot{u}_{\delta}^{i}%
-\int_{\mathbb{T}^{d}}\varepsilon_{ijkl}\partial_{\ell}u_{\delta}^{k}(\left[
u^{q},\omega_{\delta}\ast\right]  \partial_{jq}^{2}\dot{u}^{i})
\end{align*}

Similar computations lead to the identity%
\begin{align*}
&  -\left\langle \partial_{t}\left(  \mu\Delta u+\left(  \mu+\lambda\right)
\nabla\operatorname{div}u\right)  +\operatorname{div}\left(  u\left(
\mu\Delta u+\left(  \mu+\lambda\right)  \nabla\operatorname{div}u\right)
\right)  ,\dot{u}\right\rangle \\
&  =\mu\int_{\mathbb{T}^{d}}\left\vert \partial_{k}\dot{u}^{i}\right\vert
^{2}+\left(  \mu+\lambda\right)  \int_{\mathbb{T}^{d}}\left\vert
\operatorname{div}\dot{u}\right\vert ^{2}\\
&  -\mu\int_{\mathbb{T}^{d}}\partial_{k}u^{q}\partial_{q}u^{i}\partial_{k}%
\dot{u}^{i}-\mu\int_{\mathbb{T}^{d}}\partial_{k}u^{q}\partial_{k}u^{i}%
\partial_{q}\dot{u}^{i}+\mu\int_{\mathbb{T}^{d}}\operatorname{div}%
u\partial_{k}u^{i}\partial_{k}\dot{u}^{i}\\
&  -\left(  \mu+\lambda\right)  \int_{\mathbb{T}^{d}}\partial_{\ell}%
u^{q}\partial_{q}u^{\ell}\operatorname{div}\dot{u}-\left(  \mu+\lambda\right)
\int_{\mathbb{T}^{d}}\partial_{i}u^{q}\partial_{q}\dot{u}^{i}%
\operatorname{div}u+\left(  \mu+\lambda\right)  \int_{\mathbb{T}^{d}%
}\left\vert \operatorname{div}u\right\vert ^{2}\operatorname{div}\dot{u}%
\end{align*}

Putting together all the above computations, we end up with%
\begin{align*}
&  \frac{d}{dt}\int_{\mathbb{T}^{d}}\frac{\rho\left\vert \dot{u}\right\vert
^{2}}{2}+\mu\int_{\mathbb{T}^{d}}\left\vert \partial_{k}\dot{u}^{i}\right\vert
^{2}+\left(  \mu+\lambda\right)  \int_{\mathbb{T}^{d}}\left\vert
\operatorname{div}\dot{u}\right\vert ^{2}+\int_{\mathbb{T}^{d}}\varepsilon
_{ijkl}\partial_{\ell}\dot{u}_{\delta}^{k}\partial_{j}\dot{u}_{\delta}^{i}\\
&  =\mu\int_{\mathbb{T}^{d}}\partial_{k}u^{q}\partial_{q}u^{i}\partial_{k}%
\dot{u}^{i}+\mu\int_{\mathbb{T}^{d}}\partial_{k}u^{q}\partial_{k}u^{i}%
\partial_{q}\dot{u}^{i}-\mu\int_{\mathbb{T}^{d}}\operatorname{div}%
u\partial_{k}u^{i}\partial_{k}\dot{u}^{i}\\
&  +\left(  \mu+\lambda\right)  \int_{\mathbb{T}^{d}}\partial_{\ell}%
u^{q}\partial_{q}u^{\ell}\operatorname{div}\dot{u}+\left(  \mu+\lambda\right)
\int_{\mathbb{T}^{d}}\partial_{i}u^{q}\partial_{q}\dot{u}^{i}%
\operatorname{div}u-\left(  \mu+\lambda\right)  \int_{\mathbb{T}^{d}%
}\left\vert \operatorname{div}u\right\vert ^{2}\operatorname{div}\dot{u}\\
&  -\int_{\mathbb{T}^{d}}\left(  \partial_{t}\varepsilon_{ijkl}+\partial
_{q}\left(  u^{q}\varepsilon_{ijkl}\right)  \right)  \partial_{\ell}u_{\delta
}^{k}\partial_{j}\dot{u}_{\delta}^{i}\\
&  -\int_{\mathbb{T}^{d}}\varepsilon_{ijkl}(\omega_{\delta}\ast(\partial
_{\ell}u^{q}\partial_{q}u^{k}))\partial_{j}\dot{u}_{\delta}^{i}-\int
_{\mathbb{T}^{d}}\varepsilon_{ijkl}\partial_{\ell}u_{\delta}^{k}\omega
_{\delta}\ast(\partial_{j}u^{q}\partial_{q}\dot{u}^{i})\\
&  +\int_{\mathbb{T}^{d}}\varepsilon_{ijkl}(\left[  u^{q},\omega_{\delta}%
\ast\right]  \partial_{\ell q}^{2}u^{k})\partial_{j}\dot{u}_{\delta}^{i}%
+\int_{\mathbb{T}^{d}}\varepsilon_{ijkl}\partial_{\ell}u_{\delta}^{k}(\left[
u^{q},\omega_{\delta}\ast\right]  \partial_{jq}^{2}\dot{u}^{i})\\
&  -\int_{\mathbb{T}^{d}}\left\{  P\left(  \rho\right)  \partial_{j}%
u^{k}\partial_{k}\dot{u}^{j}+\left(  \rho P^{\prime}\left(  \rho\right)
-P\left(  \rho\right)  \right)  \operatorname{div}u\operatorname{div}\dot
{u}\right\}  .
\end{align*}
We multiply the above with $\sigma\left(  t\right)  $ such that we obtain%
\begin{align*}
&  \sigma\left(  t\right)  \int_{\mathbb{T}^{d}}\frac{\rho\left(  t\right)
\left\vert \dot{u}\left(  t\right)  \right\vert ^{2}}{2}+\mu\int
_{\mathbb{T}^{d}}\sigma\left\vert \partial_{k}\dot{u}^{i}\right\vert
^{2}+\left(  \mu+\lambda\right)  \int_{\mathbb{T}^{d}}\sigma\left\vert
\operatorname{div}\dot{u}\right\vert ^{2}+\int_{\mathbb{T}^{d}}\sigma
\varepsilon_{ijkl}\partial_{\ell}\dot{u}_{\delta}^{k}\partial_{j}\dot
{u}_{\delta}^{i}\\
&  =\int_{0}^{1}\int_{\mathbb{T}^{d}}\sigma\frac{\rho\left\vert \dot
{u}\right\vert ^{2}}{2}\\
&  +\mu\int_{0}^{t}\int_{\mathbb{T}^{d}}\sigma\partial_{k}u^{q}\partial
_{q}u^{i}\partial_{k}\dot{u}^{i}+\mu\int_{0}^{t}\int_{\mathbb{T}^{d}}%
\sigma\partial_{k}u^{q}\partial_{k}u^{i}\partial_{q}\dot{u}^{i}-\mu\int
_{0}^{t}\int_{\mathbb{T}^{d}}\sigma\operatorname{div}u\partial_{k}%
u^{i}\partial_{k}\dot{u}^{i}\\
&  +\left(  \mu+\lambda\right)  \int_{0}^{t}\int_{\mathbb{T}^{d}}%
\sigma\partial_{\ell}u^{q}\partial_{q}u^{\ell}\operatorname{div}\dot
{u}+\left(  \mu+\lambda\right)  \int_{0}^{t}\int_{\mathbb{T}^{d}}%
\sigma\partial_{i}u^{q}\partial_{q}\dot{u}^{i}\operatorname{div}u-\left(
\mu+\lambda\right)  \int_{0}^{t}\int_{\mathbb{T}^{d}}\sigma\left\vert
\operatorname{div}u\right\vert ^{2}\operatorname{div}\dot{u}\\
&  -\int_{0}^{t}\int_{\mathbb{T}^{d}}\sigma\left(  \partial_{t}\varepsilon
_{ijkl}+\partial_{q}\left(  u^{q}\varepsilon_{ijkl}\right)  \right)
\partial_{\ell}u_{\delta}^{k}\partial_{j}\dot{u}_{\delta}^{i}-\int_{0}^{t}%
\int_{\mathbb{T}^{d}}\sigma\varepsilon_{ijkl}(\omega_{\delta}\ast
(\partial_{\ell}u^{q}\partial_{q}u^{k}))\partial_{j}\dot{u}_{\delta}^{i}%
-\int_{0}^{t}\int_{\mathbb{T}^{d}}\sigma\varepsilon_{ijkl}\partial_{\ell
}u_{\delta}^{k}\omega_{\delta}\ast(\partial_{j}u^{q}\partial_{q}\dot{u}^{i})\\
&  +\int_{0}^{t}\int_{\mathbb{T}^{d}}\sigma\varepsilon_{ijkl}(\left[
u^{q},\omega_{\delta}\ast\right]  \partial_{\ell q}^{2}u^{k})\partial_{j}%
\dot{u}_{\delta}^{i}+\int_{0}^{t}\int_{\mathbb{T}^{d}}\sigma\varepsilon
_{ijkl}\partial_{\ell}u_{\delta}^{k}(\left[  u^{q},\omega_{\delta}\ast\right]
\partial_{jq}^{2}\dot{u}^{i})\\
&  -\int_{0}^{t}\int_{\mathbb{T}^{d}}\sigma\left\{  P\left(  \rho\right)
\partial_{j}u^{k}\partial_{k}\dot{u}^{j}+\left(  \rho P^{\prime}\left(
\rho\right)  -P\left(  \rho\right)  \right)  \operatorname{div}%
u\operatorname{div}\dot{u}\right\}  .
\end{align*}

\bigskip

\noindent\textbf{Acknowledgments.} D. Bresch and C. Burtea are supported by
the SingFlows project, grant ANR-18-CE40-0027. D. Bresch is also supported by
the Fraise project, grant ANR-16-CE06-0011 of the French National Research
Agency (ANR). C. Burtea is also supported by the CRISIS project, grant
ANR-20-CE40-0020-01. D. Bresch acknowledges the support by the National
Science Foundation while the author participated in a program hosted by the
Mathematical Sciences Research Institute (MSRI) in Berkeley, California,
during the Spring 2021 semester.

\bigskip

\bibliographystyle{alpha}
\bibliography{BibliografieGeneralaCompletaDidi}

\newcommand{\etalchar}[1]{$^{#1}$}
\begin{thebibliography}{MRRS19}

\bibitem[BB20]{BrBu1}
D~Bresch and C.~Burtea.
\newblock Global existence of weak solutions for the anisotropic compressible
  stokes system.
\newblock {\em Annales de l'Institut Henri Poincar{\'e} C, Analyse non
  lin{\'e}aire}, 37(6):1271--1297, 2020.

\bibitem[BB21]{BrBu2}
D~Bresch and C.~Burtea.
\newblock Weak solutions for the stationary anisotropic and nonlocal
  compressible navier-stokes system.
\newblock {\em Journal de Math{\'e}matiques Pures et Appliqu{\'e}es},
  146:183--217, 2021.

\bibitem[BH11]{bresch2011multi}
D~Bresch and X.~Huang.
\newblock A multi-fluid compressible system as the limit of weak solutions of
  the isentropic compressible navier--stokes equations.
\newblock {\em Archive for rational mechanics and analysis}, 201(2):647--680,
  2011.

\bibitem[BJ18]{bresch2018global}
D.~Bresch and P.-E. Jabin.
\newblock Global existence of weak solutions for compressible navier--stokes
  equations: thermodynamically unstable pressure and anisotropic viscous stress
  tensor.
\newblock {\em Annals of Mathematics}, 188(2):577--684, 2018.

\bibitem[BJW21]{bresch2021compressible}
D.~Bresch, P.-E. Jabin, and F.~Wang.
\newblock Compressible {N}avier--{S}tokes equations with heterogeneous pressure
  laws.
\newblock {\em Nonlinearity}, 34(6):4115, 2021.

\bibitem[Dan00]{danchin2000global}
R.~Danchin.
\newblock Global existence in critical spaces for compressible navier-stokes
  equations.
\newblock {\em Inventiones Mathematicae}, 141(3):579--614, 2000.

\bibitem[Dan10]{danchin2010solvability}
R.~Danchin.
\newblock On the solvability of the compressible navier--stokes system in
  bounded domains.
\newblock {\em Nonlinearity}, 23(2):383, 2010.

\bibitem[Des97]{desjardins1997regularity}
B.~Desjardins.
\newblock Regularity of weak solutions of the compressible isentropic
  {N}avier-{S}tokes equations.
\newblock {\em Communications in Partial Differential Equations},
  22(5-6):977--1008, 1997.

\bibitem[DM19]{DanchinMucha2019}
R.~Danchin and P.B. Mucha.
\newblock Compressible navier-stokes equations with ripped density.
\newblock {\em arXiv preprint arXiv:1903.09396}, 2019.

\bibitem[Fei01]{feireisl2001compactness}
E.~Feireisl.
\newblock On compactness of solutions to the compressible isentropic
  navier-stokes equations when the density is not square integrable.
\newblock {\em Commentationes Mathematicae Universitatis Carolinae},
  42(1):83--98, 2001.

\bibitem[FNP01]{feireisl2001existence}
E.~Feireisl, A.~Novotn{\`y}, and H.~Petzeltov{\'a}.
\newblock On the existence of globally defined weak solutions to the
  navier—stokes equations.
\newblock {\em Journal of Mathematical Fluid Mechanics}, 3(4):358--392, 2001.

\bibitem[Hof87]{hoff1987global}
D.~Hoff.
\newblock Global existence for 1d, compressible, isentropic {N}avier-{S}tokes
  equations with large initial data.
\newblock {\em Transactions of the American Mathematical Society},
  303(1):169--181, 1987.

\bibitem[Hof95a]{hoff1995global}
D.~Hoff.
\newblock Global solutions of the {N}avier-{S}tokes equations for
  multidimensional compressible flow with discontinuous initial data.
\newblock {\em Journal of Differential Equations}, 120(1):215--254, 1995.

\bibitem[Hof95b]{hoff1995strong}
D.~Hoff.
\newblock Strong convergence to global solutions for multidimensional flows of
  compressible, viscous fluids with polytropic equations of state and
  discontinuous initial data.
\newblock {\em Archive for rational mechanics and analysis}, 132(1):1--14,
  1995.

\bibitem[Hof02]{hoff2002dynamics}
D.~Hoff.
\newblock Dynamics of singularity surfaces for compressible, viscous flows in
  two space dimensions.
\newblock {\em Communications on Pure and Applied Mathematics: A Journal Issued
  by the Courant Institute of Mathematical Sciences}, 55(11):1365--1407, 2002.

\bibitem[HS85]{HS85}
D.~Hoff and J.~Smoller.
\newblock Solutions in the large for certain nonlinear parabolic systems.
\newblock {\em Annales de l'Institut Henri Poincare (C) Non Linear Analysis},
  2(3):213--235, 1985.

\bibitem[HS08]{hoff2008lagrangean}
D.~Hoff and M.M. Santos.
\newblock Lagrangean structure and propagation of singularities in
  multidimensional compressible flow.
\newblock {\em Archive for rational mechanics and analysis}, 188(3):509--543,
  2008.

\bibitem[Lio96]{lions1996mathematical2}
P.-L. Lions.
\newblock {\em Mathematical Topics in Fluid Mechanics: Volume 2: Compressible
  Models}, volume~2.
\newblock Oxford University Press on Demand, 1996.

\bibitem[MN{\etalchar{+}}80]{matsumura1980initial}
A.~Matsumura, T.~Nishida, et~al.
\newblock The initial value problem for the equations of motion of viscous and
  heat-conductive gases.
\newblock {\em Journal of Mathematics of Kyoto University}, 20(1):67--104,
  1980.

\bibitem[MRRS19]{merle2019implosion}
F.~Merle, P.~Rapha\"el, I.~Rodnianski, and J.~Szeftel.
\newblock On the implosion of a three dimensional compressible fluid.
\newblock {\em arXiv preprint arXiv:1912.11009}, 2019.

\bibitem[Nas62]{nash1962probleme}
J.~Nash.
\newblock Le probl{\`e}me de cauchy pour les {\'e}quations diff{\'e}rentielles
  d'un fluide g{\'e}n{\'e}ral.
\newblock {\em Bulletin de la Soci{\'e}t{\'e} Math{\'e}matique de France},
  90:487--497, 1962.

\bibitem[NS04]{NoSt}
A.~Novotn\'y and I.~Straskraba.
\newblock {\em Introduction to the mathematical theory of compressible flow}.
\newblock Oxford University Press on Demand, 2004.

\bibitem[Sol80]{solonnikov1980solvability}
V.A. Solonnikov.
\newblock Solvability of the initial-boundary-value problem for the equations
  of motion of a viscous compressible fluid.
\newblock {\em Journal of Soviet Mathematics}, 14(2):1120--1133, 1980.

\bibitem[Vai94]{vaigant1994example}
V.A. Vaigant.
\newblock An example of the nonexistence with respect to time of the global
  solution of navier-stokes equations for a compressible viscous barotropic
  fluid.
\newblock {\em Doklady Akademii Nauk}, 339(2):155--156, 1994.

\bibitem[Vra03]{Vrabie2003Semigroups}
I.I. Vrabie.
\newblock {\em Co-semigroups and applications}.
\newblock Number 191 in {N}orth-{H}olland {M}athematics {S}tudies. Elsevier,
  2003.

\end{thebibliography}

\end{document}